\documentclass[10pt,a4paper,final]{amsart}
\usepackage{amsmath,amsfonts,amssymb,dsfont,graphicx,enumerate,verbatim,esint,mathtools,mathrsfs,enumitem,ifpdf}
\usepackage[msc-links]{amsrefs}
\usepackage{kbordermatrix}
\usepackage{xargs} 
\usepackage[title]{appendix}

\ifpdf
  \usepackage[final,expansion=alltext,protrusion=alltext]{microtype} 
\fi

\usepackage{tikz}
\usetikzlibrary{matrix,arrows}
\usepackage{tikz-cd}

\tikzset{
  commutative diagrams/.cd,
  arrow style=tikz,
  diagrams={>=stealth},
  shift up/.style={
    to path={([yshift=#1]\tikztostart.east) -- ([yshift=#1]\tikztotarget.west) \tikztonodes}},
  mathdouble/.style={-,double equal sign distance}
}

\def\clap#1{\hbox to 0pt{\hss#1\hss}} 

\def\mathclap{\mathpalette\mathclapinternal} 
\def\mathrlap{\mathpalette\mathrlapinternal}
\def\mathclap{\mathpalette\mathclapinternal}

\def\mathrlapinternal#1#2{%
  \rlap{$\mathsurround=0pt#1{#2}$}} 
\def\mathclapinternal#1#2{%
  \clap{$\mathsurround=0pt#1{#2}$}}

\theoremstyle{plain}

\newtheorem*{Th*}{Theorem}

\newtheorem*{Cor*}{Corollary}

\theoremstyle{definition}

\theoremstyle{remark}

\numberwithin{equation}{section}

\makeatletter
\def\Set@Scallop[#1]#2#3{{#1}\Parens{#2}{#3}}
\newcommand\DeclareScalableOperator[2]{%
  \expandafter\def\csname#1\endcsname{\@ifnextchar[{{#2}\Set@Scallop}{{#2}\Set@Scallop[{}]}}
}
\makeatother

\def\ev{{\bar 0}}
\def\odd{{\bar 1}}
\newcommand{\defi}{\coloneqq}     
\def\csvr#1{#1_{\ev,\reals}}
\def\csvrst#1{#1^*_{\smash{\ev,\reals}}}
\def\csv#1{#1,\csvr#1} 
\def\csvst#1{#1^*,\csvrst#1} 

\newcommand{\Fa}{For all }
\newcommand{\fa}{for all }

\newcommand{\fs}{for some }
\newcommand\mathfa[1][{}]{\quad\text{\fa{#1} }}

\newcommand{\scth}{such that }
\newcommand{\AND}{and}

\newcommand\mathtxt[1]{\quad\text{{#1}}\quad}

\newcommand{\nd}{\mathtxt\AND}

\newcommand\cf{\emph{cf.}~}
\newcommand\vq{\emph{v.}~}
\newcommand\eg{\emph{e.g.}~}
\newcommand\ie{\emph{i.e.}~}

\newcommand\via{\emph{via}~}
\newcommand\loccit{\emph{loc.~cit.}}
\newcommand\opcit{\emph{op.~cit.}}

\newcommand\etc{\emph{etc.}}
\newcommand\aposteriori{\emph{a posteriori}}

\newcommand\vphi{\varphi}
\newcommand\vrho{\varrho}

\newcommand\eps{\varepsilon}
\newcommand\nats{\mathbb{N}}
\newcommand\ints{\mathbb{Z}}

\newcommand\reals{\mathbb{R}}
\newcommand\cplxs{\mathbb{C}}

\newcommand\vvoid{\varnothing}
\newcommand\sle{\leqslant}
\newcommand\sge{\geqslant}

\DeclareMathOperator\Ad{\mathrm{Ad}}
\DeclareMathOperator\ad{\mathrm{ad}}
\DeclareMathOperator\GL{\mathrm{GL}}

\DeclareMathOperator\id{\mathrm{id}}

\DeclareMathOperator\diag{\mathrm{diag}}
\DeclareMathOperator\ddiv{\mathrm{div}}

\makeatletter
\newcommand\Size[7][1]{
                                 \ifx#20%
                                        \def\r@l{}\def\r@m{}\def\r@r{}%
                                 \else%
                                    \ifx#21%
                                           \def\r@l{\bigl}\def\r@r{\bigr}\def\r@m{\bigm}%
                                    \else%
                                           \ifx#22%
                                                 \def\r@l{\Bigl}\def\r@r{\Bigr}\def\r@m{\Bigm}%
                                            \else%
                                                 \ifx#23%
                                                        \def\r@l{\biggl}\def\r@r{\biggr}\def\r@m{\biggm}%
                                                  \else
                                                        \ifx#24%
                                                        \def\r@l{\Biggl}\def\r@r{\Biggr}\def\r@m{\Biggm}%
                                                        \fi%
                                                  \fi%
                                            \fi%
                                      \fi%
                                 \fi%
                                 \ifx#10%
                                       \def\r@m{}%
                                 \fi%
                                 \r@l#3{#4}\r@m#5{#6}\r@r#7%
}%

\makeatother

\newcommand\Set[3]{
                                 \Size{#1}{\{}{#2}{|}{#3}{\}}%
}%
\newcommand\Scp[3]{
                                 \Size{#1}{(}{#2}{|}{#3}{)}%
}%
\newcommand\Dual[3]{
                                 \Size[0]{#1}{\langle}{#2}{,}{#3}{\rangle}%
}%
\newcommand\Parens[2]{
  \Size[0]{#1}{(}{#2}{}{}{)}
}
\newcommand\Bracks[2]{
  \Size[0]{#1}{[}{#2}{}{}{]}
}

\newcommand\Norm[2]{
  \Size[0]{#1}{\lVert}{#2}{}{}{\rVert}
}
\newcommand\Abs[2]{
  \Size[0]{#1}{\lvert}{#2}{}{}{\rvert}
}

\makeatletter

\def\theoremn@me#1{\csname#1name\endcsname\ignorespaces}
  
\def\thmref#1#2{\theoremn@me{#1}\ref{#1:#2}}

\def\uc#1#2{\MakeUppercase{#1}{#2}} 

\newcommand{\DefTheorem}[2]{\newenvironmentx{#1}[2][1=\empty,2=\empty]{%
    \ignorespaces%
    \ifx##2\empty%
      \begin{#2}%
    \else%
      \begin{#2}[{\uc##2}]
    \fi%
    \ifx##1\empty%
      {}%
    \else%
      \label{#1:##1}%
    \fi%
    \ignorespaces}{\end{#2}\ignorespacesafterend}}

\DefTheorem{Th}{theorem}
\DefTheorem{Prop}{proposition}
\DefTheorem{Cor}{corollary}
\DefTheorem{Lem}{lemma}
\DefTheorem{Def}{definition}
\DefTheorem{Rem}{remark}

\makeatletter

\newif\if@smallmat
\newif\if@none
\newif\if@paren
\newif\if@brack
\newif\if@brace
\newif\if@vline
\newenvironment{Matrix}[2][1]
                                 {\ifx#20%
                                        \@smallmattrue%
                                  \else%
                                         \@smallmatfalse
                                  \fi%
                                  \ifx#11%
                                         \@nonefalse\@parentrue\@brackfalse\@bracefalse\@vlinefalse%
                                  \else%
                                       \ifx#12%
                                            \@nonefalse\@parenfalse\@bracktrue\@bracefalse\@vlinefalse%
                                        \else%
                                            \ifx#13%
                                                 \@nonefalse\@parenfalse\@brackfalse\@bracetrue\@vlinefalse%
                                            \else%
                                                 \ifx#14%
                                                       \@nonefalse\@parenfalse\@brackfalse\@bracefalse\@vlinetrue
                                                 \else%
                                                       \ifx#15%
                                                             \@nonefalse\@parenfalse\@brackfalse\@bracefalse\@vlinefalse%
                                                       \else%
                                                             \@nonetrue\@parenfalse\@brackfalse\@bracefalse\@vlinefalse%
                                                       \fi%
                                                 \fi%
                                            \fi%
                                        \fi%
                                   \fi%
                                   \if@smallmat%
                                        \if@none%
                                             \begin{smallmatrix}%
                                        \else%
                                            \if@paren%
                                                  \bigl(\begin{smallmatrix}%
                                            \else%
                                                  \if@brack%
                                                          \bigl[\begin{smallmatrix}%
                                                  \else%
                                                          \if@brace%
                                                               \bigl\{\begin{smallmatrix}%
                                                          \else%
                                                               \if@vline%
                                                                    \bigl\lvert\begin{smallmatrix}%
                                                                \else%
                                                                    \bigl\lVert\begin{smallmatrix}%
                                                                \fi%
                                                          \fi%
                                                  \fi%
                                            \fi%
                                        \fi%
                                   \else%
                                        \if@none%
                                             \begin{matrix}%
                                        \else%
                                            \if@paren%
                                                  \begin{pmatrix}%
                                            \else%
                                                  \if@brack%
                                                          \begin{bmatrix}%
                                                  \else%
                                                          \if@brace%
                                                               \begin{Bmatrix}%
                                                          \else%
                                                               \if@vline%
                                                                    \begin{vmatrix}%
                                                                \else%
                                                                    \begin{Vmatrix}%
                                                                \fi%
                                                          \fi%
                                                  \fi%
                                            \fi%
                                        \fi%
                                   \fi}%
                                  {\if@smallmat%
                                        \if@none%
                                             \end{smallmatrix}%
                                        \else%
                                            \if@paren%
                                                  \end{smallmatrix}\bigr)%
                                            \else%
                                                  \if@brack%
                                                          \end{smallmatrix}\bigr]%
                                                  \else%
                                                          \if@brace%
                                                               \end{smallmatrix}\bigr\}%
                                                          \else%
                                                               \if@vline%
                                                                    \end{smallmatrix}\bigr\rvert%
                                                                \else%
                                                                    \end{smallmatrix}\bigr\rVert%
                                                                \fi%
                                                          \fi%
                                                  \fi%
                                            \fi%
                                         \fi%
                                   \else%
                                        \if@none%
                                             \end{matrix}%
                                        \else%
                                            \if@paren%
                                                  \end{pmatrix}%
                                            \else%
                                                  \if@brack%
                                                          \end{bmatrix}%
                                                  \else%
                                                          \if@brace%
                                                               \end{Bmatrix}%
                                                          \else%
                                                               \if@vline%
                                                                    \end{vmatrix}%
                                                                \else%
                                                                    \end{Vmatrix}%
                                                                \fi%
                                                          \fi%
                                                  \fi%
                                            \fi%
                                        \fi%
                                   \fi}%

\newcount\dotcnt\newdimen\deltay
\def\Ddot#1#2(#3,#4,#5,#6){\deltay=#6\setbox1=\hbox to0pt{\smash{\dotcnt=1
\kern#3\loop\raise\dotcnt\deltay\hbox to0pt{\hss#2}\kern#5\ifnum\dotcnt<#1
\advance\dotcnt 1\repeat}\hss}\setbox2=\vtop{\box1}\ht2=#4\box2}

\def\trimatclip(#1,#2,#3,#4)[#5,#6,#7,#8]{%
  \left(%
  \begin{tikzpicture}[baseline=(current bounding box.center)]
    \clip (#5,#6) rectangle (#7,#8);
    \matrix (m) [ampersand replacement=\&,matrix of math nodes,nodes in empty cells]{
      {#2}\&{#4}\& \&{#4}\\
      {#1}\\
      \& \& \&{#4}\\
      {#1}\& \& {#1} \& {#3}\\
      } ;
    \draw[thick,dotted] (m-1-1) -- (m-4-4);
    \draw[thick,dotted] (m-1-2) -- (m-1-4);
    \draw[thick,dotted] (m-1-4) -- (m-3-4);   
    \draw[thick,dotted] (m-1-2) -- (m-3-4);
    \draw[thick,dotted] (m-2-1) -- (m-4-3);
    \draw[thick,dotted] (m-2-1) -- (m-4-1);
    \draw[thick,dotted] (m-4-1) -- (m-4-3);
  \end{tikzpicture}%
  \right)
}

\def\trimat(#1,#2,#3,#4){%
  \trimatclip(#1,#2,#3,#4)[-0.7,-0.9,0.7,0.78]
}

\def\trismallmat(#1,#2,#3,#4){%
  \left(\begin{Matrix}[0]0
    {#2}\Ddot6$\cdot$(4pt,-1pt,2.4pt,-2.4pt)&{#4}\Ddot3$\cdot$(3.5pt,-0.5pt,3pt,0pt)\Ddot3$\cdot$(3pt,-0.5pt,3pt,-3pt)&&{#4}\Ddot3$\cdot$(-0.5pt,-2pt,0pt,-2.4pt)\\
    {#1}\Ddot3$\cdot$(-0.5pt,-2.5pt,0pt,-2pt)\Ddot3$\cdot$(4pt,-0.5pt,2.4pt,-2.4pt)&\\
    &&&{#4}\\
    {#1}\Ddot3$\cdot$(3.5pt,0pt,3pt,0pt)&&{#1}&{#3}
  \end{Matrix}\right)
}

\makeatother

\def\ger{\mathfrak}

\newcommand\CategoryTypeface{\mathbf}
\def\cat{\CategoryTypeface}

\newcommand\SheafTypeface{\mathcal}
\def\sh{\SheafTypeface}

\DeclareMathOperator\Sets{\cat{Sets}}

\DeclareScalableOperator{Sh}{\cat{Sh}}

\newcommand\Db{\sh Db}

\newcommand{\lBr}{[\kern-.65ex[}
\newcommand{\rBr}{]\kern-.65ex]}

\DeclareMathOperator\supp{\mathrm{supp}}
\DeclareMathOperator\str{\mathrm{str}}
\DeclareMathOperator\tr{\mathrm{tr}}

\DeclareMathOperator\Herm{\mathrm{Herm}}
\DeclareMathOperator\LT{\mathscr L}

\def\DSO{\DeclareScalableOperator}

\DeclareScalableOperator{Ber}{\mathrm{Ber}}
\DeclareScalableOperator{Vol}{\mathrm{Vol}}
\DeclareScalableOperator{Form}{\Omega}
\DeclareScalableOperator{Ind}{\mathrm{Ind}}
\DeclareScalableOperator{Coind}{\mathrm{Coind}}
\DeclareScalableOperator{Der}{\mathrm{Der}}
\DeclareScalableOperator{Maps}{\mathrm{Maps}}
\DeclareScalableOperator{Ct}{\mathcal{C}} 
\DeclareScalableOperator{GCt}{\underline{\mathcal{C}}} 
\DeclareScalableOperator{Lp}{\mathbf{L}} 
\DeclareScalableOperator{Hom}{\mathrm{Hom}} 
\DeclareScalableOperator{GHom}{\underline{\mathrm{Hom}}} 
\DeclareScalableOperator{GEnd}{\underline{\mathrm{End}}} 
\DeclareScalableOperator{GPol}{\underline{\sh{P}}} 
\DeclareScalableOperator{GLin}{\underline{\sh{L}}} 
\DeclareScalableOperator{End}{\mathrm{End}} 
\DeclareScalableOperator{Aut}{\mathrm{Aut}} 
\DeclareScalableOperator{Uenv}{\mathfrak{U}} 
\DeclareScalableOperator{Cuenv}{\widehat{\mathfrak{U}}} 
\DeclareScalableOperator{ABer}{\Abs0{\mathrm{Ber}}}
\DeclareScalableOperator{Dist}{\mathcal{D}'}
\DeclareScalableOperator{Sw}{\mathscr{S}}
\DeclareScalableOperator{TDi}{\mathscr{S}'}
\DeclareScalableOperator{PW}{\mathscr{Z}}
\DeclareScalableOperator{PWFn}{\mathscr{Z}'}

\DSO{ShHom}{\sh Hom} 
\DSO{ShGHom}{\underline{\sh Hom}} 
\DSO{ShEnd}{\sh End} 
\DSO{ShGEnd}{\underline{\sh End}} 
\DSO{ShDer}{\sh Der} 
\DSO{ShGDer}{\underline{\sh Der}} 

\newcommand{\fibint}[2][\empty]{\!\sideset{_{#1}}{_{#2}}{\fint}} 
\newcommand{\tfibint}[2][\empty]{\!\sideset{_{#1}}{_{#2}}{\textstyle\fint}} 

\newcommand\Define[1]{\emph{#1}}

\newcommand{\C}{\mathbb{C}}

\renewcommand{\kbldelim}{(}
\renewcommand{\kbrdelim}{)}

\begin{document}

\title{Superbosonisation via Riesz superdistributions}

\author[Alldridge]
{Alexander Alldridge}
\address{Universit\"at zu K\"oln\\
Mathematisches Institut\\
Weyertal 86-90\\
50931 K\"oln, Germany}
\email{alldridg@math.uni-koeln.de}

\author[Shaikh]{Zain Shaikh}
\address{Universit\"at zu K\"oln\\
Mathematisches Institut\\
Weyertal 86-90\\
50931 K\"oln, Germany}
\curraddr{
  Universit\"at Paderborn\\
  Institut f\"ur Mathematik\\
  Warburger Str.~100\\
  33100 Paderborn, Germany
}

\email{zain@math.upb.de}

\begin{abstract}
  The superbosonisation identity of Littelmann--Sommers--Zirnbauer is a new tool for use in studying universality of random matrix ensembles \via supersymmetry, which is applicable to non-Gaussian invariant distributions. We give a new conceptual interpretation of this formula, linking it to harmonic superanalysis of Lie supergroups and symmetric superspaces, in particular, to a super-generalisation of the Riesz distributions. Using the super-Laplace transformation of generalised superfunctions, the theory of which we develop, we reduce the proof to the computation of the Gindikin gamma function of a Riemannian symmetric superspace, which we determine explicitly.
\end{abstract}

\thanks{This research was funded by Deutsche Forschungsgemeinschaft, grant nos.~DFG ZI 513/2-1 and SFB TR/12, and the Institutional Strategy of the University of Cologne within the German Excellence Initiative}



\maketitle

\section{Introduction}

Supersymmetry was introduced in physics as a means for formulating Bose--Fermi symmetry in quantum field theory. Since the advent of supergravity, it is usually connected to superstring theory. Although this relationship is indeed intimate and fundamental, supersymmetry is also deeply rooted in the physics of condensed matter. The so-called \emph{supersymmetry method}, developed by Efetov and Wegner \cite{efetov-susy}, has been used to great effect in the study of disordered systems, in particular in connection to the metal-insulator transition; or in other words, in the analysis of localisation and delocalisation for certain random matrix ensembles \cites{zirn-cmp,zirn-prl}. 

In connection with the physics of thin wires, the subject was well studied in the 1990s; it has recently gained new interest, since the `symmetry classes' investigated in this context \cites{zirn-rmt,az-10fold,hhz} have been found to occur as `edge modes' of certain 2D systems termed `topological insulators' (and superconductors) \cite{freed-moore}.

Mathematically, several aspects of the method beg justification. One both subtle and salient point is the transformation of certain integrals over flat superspace in high dimension $N\to\infty$, which encode statistical quantities of a given random matrix ensemble, into integrals over a curved superspace of fixed rank and dimension, more tractable to saddle-point analysis. 

This step relies on an integral transform, the so-called Hubbard--Stratonovich transformation \cite{schaefer-wegner}, which, in its traditional form, assumes a Gaussian distribution for the initial data. Any generalisation beyond this class of random matrices is challenging analytically, severely limiting the range of the method: Questions such as universality for invariant ensembles prompted the development of new, more robust and versatile tools. 

To that end, a new transform, dubbed `superbosonisation', was introduced by Efetov \emph{et al.} \cite{efetov-etal}, based on ideas of Lehmann, Saher, Sokolov, and Sommers \cite{lsss}, and of Hackenbroich and Weidenm\"uller \cite{hack-wei}. Sommers \cite{sommers-app} was the first to consider it for arbitrary probability distributions with unitary symmetry. In parallel, Guhr \cite{guhr-jphys}, Guhr \emph{et al.}~\cite{kgg} defined a generalised Hubbard--Stratonovich transformation for non-Gaussian ensembles. Following the work of Fyodorov \cite{fyodorov} in the non-supersymmetric situation, Littelmann--Sommers--Zirnbauer in their seminal paper \cite{LSZ08} extended superbosonisation to include the cases of orthogonal and symplectic symmetry. At the same time, they gave a mathematically rigorous derivation of the superbosonisation identity. Littelmann--Sommers--Zirnbauer \cite{LSZ08} also state that when both methods (the Hubbard--Stratonovich method and superbosonisation) are applicable, they are equivalent to each other; this was proved in Ref.~\cite{ksg}. 

\bigskip\noindent
In general, the superbosonisation identity holds in the unitary, orthogonal, and unitary-symplectic symmetry cases. We consider only unitary symmetry here, although neither our methods nor our results are restricted to this case. 

One considers the space $W=\cplxs^{p|q\times p|q}$ of square super-matrices and a certain subsupermanifold $\Omega$ of purely even codimension with underlying manifold 
\[
  \Omega_0=\Herm^+(p)\times\mathrm U(q),
\]
the product of the positive Hermitian $p\times p$ matrices with the unitary $q\times q$ matrices. 

Let $f$ be a superfunction defined and holomorphic on the tube domain based on $\Herm^+(p)\times\Herm(q)$. If $f$ has sufficient decay at infinity along $\Omega_0$ (\ie along $\Herm^+(p)$), then the superbosonisation identity states 
\begin{equation}\label{eq:sbos}
  \int_{\cplxs^{p|q\times n}\oplus\cplxs^{n\times p|q}}\Abs0{Dv}\,f(Q(v))=C\int_\Omega\Abs0{Dy}\,\Ber0y^nf(y),
\end{equation}
where $C$ is some finite positive constant depending only on $p$, $q$, and $n$. Here, $Q$ is the quadratic map $Q(v)=vv^*$, $\Abs0{Dv}$ is the standard Berezinian density, and $\Abs0{Dy}$ is a Berezinian density on $\Omega$, invariant under a natural transitive supergroup action. 

The precise meaning of all the quantities involved will be made clear in the course of the paper. However, let us remark that any unitary invariant superfunction on $\cplxs^{p|q\times n}\oplus\cplxs^{n\times p|q}$ may be written in the form $f(Q(v))$. A notable feature of the formula is thus that puts the `hidden supersymmetries' (from $\GL(p|q,\cplxs)$) in evidence through the invariant integral over the homogeneous superspace $\Omega$, whereas the `manifest symmetries' (from $\GL(n,\cplxs)$) only enter \via some character (namely, $\Ber0y^n$).

A remarkable special case of the identity occurs when $p=0$. Then Equation \eqref{eq:sbos} reduces to
\[
  \int_{\cplxs^{0|q\times n}\oplus\cplxs^{n\times 0|q}}\Abs0{Dv}\,f(Q(v))=C\int_{\mathrm U(q)}\Abs0{Dy}\,\det(y)^{-n}f(y),
\]
which is known as the bosonisation identity in physics. Notice that the left-hand side is a purely fermionic Berezin integral, whereas the right-hand side is purely bosonic. Formally, it turns fermions $\psi\bar\psi$ into bosons $e^{i\vphi}$. 

In this case, the identity can be proved by the use of the Schur orthogonality relations: It expresses the fact that, up to a constant factor, the $2n^{\text{th}}$ homogeneous part of $f$ equals its projection onto the Peter--Weyl component $L^2(\mathrm U(q))_n$ with spherical vector $\det^n$. The doubling of degree is explained by viewing $\mathrm U(q)$ as the symmetric space $(\mathrm U(q)\times\mathrm U(q))/\mathrm U(q)$; for instance, $L^2(\mathrm U(2))_n$ is $\End0{V_n}$, where $V_n$ is the $\mathrm U(2)$-representation on polynomials $p(z_1,z_2)$ homogeneous of degree $n$. For the case of $q=1$, the left-hand side is $f^{(n)}(0)$, up to some constant factor, and we obtain the Cauchy integral formula. 

On the other hand, if $q=0$, then $\Omega=\Herm^+(p)$, and the right-hand side 
\[
  \Dual0{T_n}f\defi\int_{\Herm^+(p)}\Abs0{Dy}\,\det(y)^nf(y)
\]
is the so-called (unweighted) \Define{Riesz distribution}. In this case, Equation \eqref{eq:sbos} is due to Ingham and Siegel \cites{ingham,siegel}. Its first use in the context that inspired superbosonisation was by Fyodorov \cite{fyodorov}. Moreover, it is known to admit a far-reaching generalisation in the framework of Euclidean Jordan algebras \cite{FK94}. 

This observation links the identity to equivariant geometry, Lie theory, and harmonic analysis; and thus forms the starting point of the present paper. Our strategy is to exploit the transitive action of a certain supergroup on $\Omega$ to compute the Laplace transform $\LT(T_n)$ of the functionals $T_n$. The corresponding transformed identity is easy to verify, since the Laplace transform of the left-hand side is straightforward to evaluate. 

Of course, since the geometry is more complicated for $q>0$, the evaluation of the right-hand side $\LT(T_n)$ becomes more delicate. Also, a theory of Laplace transforms for superdistributions had to be developed, because it was unavailable in the literature. Here, a technical difficulty is that $T_n$ is not obviously a superdistribution (although \aposteriori, Equation \eqref{eq:sbos} shows that it is), but rather a functional on a space of holomorphic superfunctions. Thus, we are obliged to study the Laplace transformation on various different spaces of generalised functions. 

Besides providing a conceptual framework in which Equation \eqref{eq:sbos} follows with relative ease, our approach also establishes a connection to analytic representation theory that previously went unnoticed. Namely, for a suitable choice of normalisation, the constant $C$ in Equation \eqref{eq:sbos} is
\[
  C=\sqrt\pi^{np}\,\Gamma_\Omega(n)^{-1},
\]
where $\Gamma_\Omega(n)$ is the evaluation at $(n,\dotsc,n)$ of the meromorphic function of $p+q$ indeterminates, known as the Gindikin $\Gamma$ function for $q=0$. In \thmref{Th}{gamma}, we explicitly determine $\Gamma_\Omega(\mathbf m)$ for any $\mathbf m=(m_1,\dotsc,m_{p+q})$, as follows: 
\[
  (2\pi)^{p(p-1)/2}\prod_{j=1}^p \Gamma(m_j-(j-1))\prod_{k=1}^q \frac{\Gamma(q-(k-1))}{\Gamma(m_{p+k}+q-(k-1))}\frac{\Gamma(m_{p+k}+k)}{\Gamma(m_{p+k}-p+k)}.
\]
This function has zeros and poles for $q>0$, whereas for $q=0$, it only has poles.

When $q=0$, this function is closely related to the $c$-function of the Riemannian symmetric space $\Omega=\Herm^+(p)$. Moreover, the renormalised Riesz distribution $R_n\defi\Gamma_\Omega(n)^{-1}T_n$ defines the unitary structure on the holomorphic discrete series representation of $\mathrm U(p,p)$ whose lowest $\mathrm U(p)$-type is the character $\det(z)^n$ \cites{FK90,rossi-vergne}.

The Gindikin $\Gamma$ function $\Gamma_\Omega$ also appears in the $b$-function equation for the relatively invariant polynomial $\det(z)$, \via
\[
  \det\Parens1{\tfrac\partial{\partial z}}\det(z)^n=n(n+1)\dotsm(n+p-1)\det(z)^{n-1}=(-1)^p\frac{\Gamma_\Omega(1-n)}{\Gamma_\Omega(-n)}\det(z)^{n-1}.
\]
Here, $\det(\tfrac\partial{\partial z})$ is the polynomial differential operator obtained by inserting the matrix $\smash{\Parens1{\tfrac\partial{\partial z_{ij}}}}$ of coordinate derivations into the determinant. 

This leads to a functional equation for $R_n$ that can be exploited to give a meromorphic extension of $R_n$ as a distribution-valued function of $n$, a fact famously applied by Rossi--Vergne in their proof of the analytic extension of the holomorphic discrete series \cite{rossi-vergne}. The implications of these connections for the representation theory of supergroups will be investigated in a forthcoming paper. 

\bigskip\noindent
Let us give a brief synopsis of the paper's contents. In Section \ref{sec:two}, we give the basic setup for the statement and proof of the superbosonisation formula, introducing the relevant supergroups and the functionals which define both sides of the equation. We give a proof of the identity in \thmref{Th}{sbos}, up to the computation of the Laplace transform $\LT(T_n)$ of the right-hand side, which is deferred to Section \ref{sec:lap-conical}. In that section, we actually discuss more generally the so-called conical superfunctions attached to $\Omega$, and determine their Laplace transforms. The main step is the explicit determination of the Gindikin $\Gamma$ function $\Gamma_\Omega$, in \thmref{Th}{gamma}.

The paper is complemented with an extensive appendix section. In Appendix \ref{app:functorpoints}, we discuss the language of (generalised) points, which will be an indispensable tool. Appendix \ref{app:super-int} covers the theory of Berezinian fibre integrals, which is used throughout Section \ref{sec:lap-conical}. Finally, the theory of generalised superfunctions and their Laplace transforms is developed in Appendix \ref{app:super-lap}. These techniques form the basis of our proof of the superbosonisation identity in \thmref{Th}{sbos}.

\bigskip\noindent
\emph{Acknowledgement.} We wish thank Martin Zirnbauer for extensive discussions and for bringing this topic to our attention. The first named author wishes to thank Jacques Faraut for his interest, and the second named author wishes to thank Bent \O{}rsted for some useful comments. Last, not least we thank the anonymous referees for their detailed suggestions, which greatly improved the article. 

\section{The superbosonisation identity}\label{sec:two}

In this section, we set up the basic framework for formulating the superbosonisation identity. We reduce the proof to a super version of the theory of the Laplace transform of generalised functions, as discussed in \ref{app:super-lap}, and the explicit computation of certain Laplace transforms, which is performed in Section \ref{sec:lap-conical}.

\subsection{Preliminaries}

In this article, we will use the machinery of supergeometry extensively. As general references the reader may, for instance, consult Refs.~\cites{fioresi-bk,QFS99,Lei80,Man88}. In this subsection, we review the basic definitions, fixing our notation and highlighting points in which we deviate slightly from the standard lore. We will introduce further notions as needed, referring the reader to Appendix \ref{app:super} and \ref{app:super-int} for detailed summaries.

\emph{Superspaces.} We will work exclusively in the category of $\cplxs$-superspaces and certain full subcategories thereof. By definition, a \Define{$\cplxs$-superspace} is a pair $X=(X_0,\sh O_X)$; here, $X_0$ is a topological space and $\sh O_X$ is a sheaf of unital supercommutative superalgebras over $\cplxs$, whose stalks $\sh O_{X,x}$ are local rings. A morphism $f:X\to Y$ is by definition a pair $(f_0,f^\sharp)$ comprising a continuous map $f_0:X_0\to Y_0$ and a sheaf map $f^\sharp:f_0^{-1}\sh O_Y\to O_X$, which is local in the sense that $f^\sharp(\ger m_{Y,f_0(x)})\subseteq\ger m_{X,x}$ for any $x$. 

Global sections $f\in\Gamma(\sh O_X)$ of $\sh O_X$ are called \Define{superfunctions}. Due to the locality condition, the \Define{value} $f(x)\defi f+\ger m_{X,x}\in\sh O_{X,x}/\ger m_{X,x}=\cplxs$ is defined for any $x$. \Define{Open subspaces} of a $\cplxs$-superspace $X$ are given by $(U,\sh O_X|_U)$, for any open $U\subseteq X_0$.

\emph{Model spaces.} We consider two types of model spaces. Firstly, when $V=V_\ev\oplus V_\odd$ is a complex super-vector space, we define 
\[
  \sh O_V\defi\mathscr H_{V_\ev}\otimes_\cplxs\textstyle\bigwedge(V_\odd)^*
\]
where $\mathscr H$ denotes the sheaf of holomorphic functions. The superspace 
\[
L(V)\defi(V_\ev,\sh O_V)
\]
is called the \Define{complex affine superspace} associated with $V$. 

Secondly, if in addition, we are given a real form $\smash{\csvr V}$ of $V_\ev$, then the pair $(\smash{\csv V})$ is called a \Define{\emph{cs} vector space}, and we define 
\[
\sh O_{\csv V}\defi\sh C^\infty_{\csvr V}\otimes_\cplxs\textstyle\bigwedge(V_\odd)^*, 
\]
where $\sh C^\infty$ denotes the sheaf of complex-valued smooth functions. Here, by a \Define{real form} $V_\reals$ of a complex vector space $V$, we mean a real subspace such that the map $V_\reals\to V$ induces an isomorphism $V_\reals\otimes\cplxs\to V$. The superspace 
\[
L(\smash{\csv V})\defi(\smash{\csvr V,\sh O_{\csv V}})
\]
is called the \Define{\emph{cs} affine superspace} associated with $(\smash{\csv V})$. (The \emph{cs} terminology is due to J.~Bernstein.)

\emph{Complex supermanifolds and \emph{cs} manifolds.} Consider now a superspace $X$ whose underlying topological space $X_0$ is Hausdorff and that admits a cover by open subspaces, which are isomorphic to open subspaces of some $L(V)$ (resp.~$L(\smash{\csv V})$), where $V$ (resp.~$(\smash{\csv V})$) may vary. Then $X$ is called a \Define{complex supermanifold} (resp.~a \Define{\emph{cs} manifold}) of dimension 
\[
  \dim X=\dim_\cplxs V_\ev|\dim_\cplxs V_\odd\quad\text{(resp.~}\dim_{cs}X=\dim_\reals V_{\ev,\reals}|\dim_\cplxs V_\odd\text{).}
\]

Both complex supermanifolds and \emph{cs} manifolds form full subcategories of the category of $\cplxs$-superspaces that admit finite products. Since \emph{cs} manifolds and complex supermanifolds belong to the same larger category of $\cplxs$-superspaces, one can easily consider morphisms between such spaces.\footnote{In passing, note that E.~Witten (Notes on supermanifolds and integration, arxiv:1209.2199) has recently advocated the study of \emph{cs} submanifolds of complex supermanifolds; the setting of functors on \emph{cs} manifolds or slightly more general $\cplxs$-superspaces seems to be well-adapted to such a study.} A remarkable fact \cite{ahw-sing}, generalising the well-known Leites theorem on morphisms, is that for a \emph{cs} manifold $S$ and complex supermanifolds $X$ and $Y$, there is a bijection
\[
  \Hom0{S,X\times Y}\to\Hom0{S,X}\times\Hom0{S,Y}
\]
whose components are given by composition with the projections $X\times Y\to X$ and $X\times Y\to Y$, respectively. Here, $X\times Y$ is the direct product of complex supermanifolds. 

\emph{Supergroups and supergroup pairs.} Group objects in the category of complex supermanifolds (resp.~\emph{cs} manifolds) are called complex Lie supergroups (resp.~\emph{cs} Lie supergroups). The category of complex Lie supergroups is equivalent to the category of complex supergroup pairs, \cf Ref.~\cite{fioresi-bk}. 

These are pairs $(\ger g,G_0)$ consisting of a complex Lie superalgebra $\ger g$ and complex Lie group $G_0$, \scth $\ger g_\ev$ is the Lie algebra of $G_0$, endowed with an action $\Ad$ of $G_0$ on $\ger g$ by Lie superalgebra automorphisms that extends the adjoint action of $G_0$ on $\ger g_\ev$, and whose derivative coincides with the restriction of the bracket of $\ger g$. Morphisms $(\ger g,G_0)\to(\ger h,H_0)$ of such pairs are given by pairs $(d\phi,\phi_0)$ of a morphism $\phi_0$ of Lie groups and a $\phi_0$-equivariant morphism $d\phi$ of Lie superalgebras. 

By assuming instead that $G_0$ be a real Lie group and that $\ger g_\ev$ is the complexification of the Lie algebra of $G_0$, and modifying all definitions accordingly, one obtains the category of \emph{cs} supergroup pairs, which is equivalent to the category of \emph{cs} Lie supergroups.

\emph{Forms of complex supergroups}. Given a Lie supergroup (complex or \emph{cs}), a \Define{closed subsupergroup} is a closed subsupermanifold that is a Lie supergroup, \scth the embedding morphism is a morphism of supergroups. 

More generally, one can consider morphisms $\vphi:H\to G$ of $\cplxs$-superspaces from a \emph{cs} supergroup $H$ to a complex Lie supergroup $G$ that are closed embeddings, \ie $\vphi_0$ is a closed topological embedding and $d\vphi$ is an injective map. If now $\vphi$ is a morphism of group objects in the sense that $\vphi$ induces group homomorphisms $\Hom0{S,H}\to\Hom0{S,G}$ for any \emph{cs} manifold $S$, then we say that $H$ is a closed \emph{cs} subsupergroup of $G$. (The sets $\Hom0{S,G}$ are groups by the remarks on products made above.)

In particular, we say that $H$ is a \Define{\emph{cs} form} of $G$ if the Lie superalgebra $\ger g$ of $G$ coincides with that of $H$, and $H_0$ is a real form of $G_0$, or equivalently, that $(\ger g,H_0)$ is a \emph{cs} supergroup pair. 

Here, by a \Define{real form} of a complex Lie group $G_0$, we mean a closed subgroup whose Lie algebra is a real form of the Lie algebra of $G_0$ (considered as a complex vector space). In view of the above remarks, to define a \emph{cs} form of $G$, it suffices to specify a real form $H_0$ of $G_0$.

\emph{$S$-valued points}. Superfunctions and morphisms are not determined by their values at points. If the notion of points is extended, one can, however, work with supergeometric objects largely as if they were ordinary manifolds. 

A point of $X$ contains the same information as a map $*\to X$ from a singleton space. Instead of the singleton space, one might actually take any other space $S$ of the same type as $X$. An \Define{$S$-valued point} is then a map $x:S\to X$, thought of as a family of points parametrised by $S$. To stress this point of view, we use the notation $x\in_SX$ to indicate an $S$-valued point of $X$. The concept directly carries over to superspaces (and indeed, to any category).

This idea is common in physics, where one treats `elements' of supermanifolds as quantities containing sufficiently general Grassmann variables. In the above terminology, these are nothing but $S$-valued points for a `superpoint' $S$, \ie a superspace over the singleton set, with a Grassmann algebra as the (constant) sheaf of superfunctions. 

Technically, it is better to allow more general auxiliary superspaces $S$. Following this approach indeed greatly simplifies practical computations. For example, the $S$-valued points of $\GL(p|q,\cplxs)$ (for $S$ a complex supermanifold resp.~a \emph{cs} manifold) are just the even $p|q\times p|q$-matrices with entries in the set $\Gamma(\sh O_S)$ of superfunctions on $S$ \cite{Lei80}. This makes it possible to perform calculations with these supergroups in terms of matrices. More details on $S$-valued points are summarised in \ref{app:functorpoints}. See also Ref.~\cite{Man88}.

\subsection{The relevant supergroups}

We now begin building the natural framework for the superbosonisation identity proper. Its right-hand side is the integral over a homogeneous superspace. The integrand turns out to be equivariant, \ie it transforms with respect to a character under a transitive supergroup action. We now introduce the supergroups which are relevant to its precise definition.

Consider the complex Lie supergroup $G'_\cplxs \defi \GL(2p|2q,\cplxs)$, with Lie superalgebra $\ger g' = \ger{gl}(2p|2q,\cplxs)$.\footnote{The notation comes from the fact that $G'_\cplxs$ arises as the Howe dual partner of the Lie group $G_\cplxs \defi \GL(n,\C)$ in the oscillator representation of $\ger{spo}(V)$, where $V\defi\cplxs^{p|q\times n}\oplus\cplxs^{n\times p|q}$.} The underlying Lie group of $G'_\cplxs$ is $G'_{\cplxs,0}=\GL(2p,\cplxs)\times\GL(2q,\cplxs)$. We will write $S$-valued points $g \in_S G'_\cplxs$ in the non-standard form 
\[
\kbordermatrix{
 & p|q & p|q \\
p|q & A & B \\
p|q & C & D},
\]
rather than in the more customary even-odd decomposition. Although this may seem unnatural at first, it will bear fruit below.

We define a Lie supergroup $K_\C:=\GL(p|q,\C) \times \GL(p|q,\C)$, which is a closed complex subsupergroup of $G'_\C$, embedded as
\[
\kbordermatrix{
 & p|q & p|q \\
p|q & A & 0 \\
p|q & 0 & D},
\]
with $A,D \in_S \GL(p|q,\C)$. Further, $\ger k \defi \mathrm{Lie}(K_\C)$.

We now define a \emph{cs} form $H$ of $K_\C$ that is of interest. We do so by specifying a real form of $K_{\cplxs,0}$ by 
\[
  H_0\defi\GL(p,\cplxs)\times\mathrm U(q)\times\mathrm U(q).
\]
Here, the latter is embedded into $K_{\cplxs,0}$ as the set of all matrices of the form 
\[
  \kbordermatrix{
     & p & q & p & q \\
    p & A & &  & \\
    q & & D & & \\
    p & & & (A^*)^{-1} & \\
    q& & & & D'}
\]
with $A\in \GL(p,\C)$, $D,D'\in \mathrm{U}(q)$.

To conclude this section, we define the homogeneous superspace $\Omega$. To that end, consider the complex super-vector space
$$
W \defi\ger{gl}(p|q,\cplxs)=\C^{p|q\times p|q}.
$$
We define a partial action of $G'_\C$ on the complex supermanifold $L(W)$: For 
$$
g' = \begin{Matrix}1A&B\\C&D\end{Matrix} \in_S G'_\C
$$
and $Z\in_S L(W)$, the action is by fractional linear transformations
$$
g'.Z = (AZ+B)(CZ+D)^{-1},
$$
whenever $CZ+D\in_S\GL(p|q,\cplxs)$. Observe that if $g\in_SK_\cplxs$, then this condition is always satisfied, so that $K_\cplxs$ acts on $L(W)$. This induces an action of the \emph{cs} form $H$ of $K_\cplxs$ on $L(W)$.

We define $\Omega\defi H.1$, the orbit through the ordinary point $1\in L(W)_0=W_\ev$. It is a \emph{cs} manifold, whose underlying manifold is the homogeneous space
$$
\Omega_0=\GL(p,\cplxs)/\mathrm U(p)\times(\mathrm U(q)\times\mathrm U(q))/\mathrm U(q)=\Herm^+(p) \times \mathrm{U}(q),
$$
where $\Herm^+(p)$ is the cone of positive definite Hermitian $p \times p$ matrices.

\subsection{\protect{The $Q$-morphism}}

The superbosonisation identity describes the transformation of an integral under a certain quadratic morphism $Q$. We introduce it in a form convenient for our purposes.

Consider the complex super-vector spaces 
$$
V\defi\cplxs^{p|q\times n}\oplus\cplxs^{n\times p|q}\nd U\defi \C^{(n+p|q)\times (n+p|q)}.
$$
We have an embedding $V \hookrightarrow U$, given by
\[
  \kbordermatrix{
    & n & p|q  \\
    n & 0 & a'   \\
    p|q & a & 0 },
\]
where $(a,a')\in V$. Further, we have an embedding of $W$ into $U$, given by
  \[
  \kbordermatrix{
 & n & p|q  \\
n & 0 & 0   \\
p|q & 0 & w},
  \]
  where $w\in W$. We obtain embeddings $L(V)\to L(U)$ and $L(W)\to L(U)$ as closed complex subsupermanifolds. 
  
  A real form of $U_\ev$ is given by $\csvr U=\Herm(n+p)\times\Herm(q)$, \ie the set of matrices
  \[  
\kbordermatrix{
 & n & p & q  \\
n & b_1 & b_2 & 0 \\
p & b_3 & b_4 & 0\\
q & 0 & 0 & b'}
  \]
  with $\begin{Matrix}0b_1&b_2\\b_3&b_4\end{Matrix}\in\Herm(n+p)$, $b'\in\Herm(q)$. We set 
  \[
    \csvr V\defi V\cap\csvr U\nd\csvr W\defi W\cap\csvr U.
  \]
  These are real forms of $V_\ev$ and $W_\ev$, respectively. We obtain \emph{cs} manifolds denoted by $L(\csv U)$, $L(\csv V)$, and $L(\csv W)$.
  
  For $x,y\in U$, we define
  \[
    P(x)y\defi xyx\in U.
  \]
  Since this expression is polynomial in $x$ and $y$, there is a unique morphism $L(U)\times L(U)\to L(U)$ given on $S$-valued points by $(x,y)\mapsto P(x)y$. We let $Q:L(V)\to L(W)$ denote the morphism given by $Q(x)=P(x)c_n$ where 
\[
c_n=\kbordermatrix{
 & n & p|q  \\
n & 1 & 0   \\
p|q & 0 & 0}.
\]

\medskip\noindent
For the remainder of the section, recall from Appendix \ref{app:schwartz} that for any \emph{cs} vector space $(\csv E)$, the \Define{Schwartz space} $\Sw0{\csv E}$ is defined to be the set of superfunctions $f\in\Gamma(\sh O_{\smash{\csv E}})$ \scth 
\[
  \forall N\in\nats,D\in S(E)\,:\,\sup\nolimits_{x\in\csvr E}(1+\Norm0x)^N\Abs0{f(D;x)}<\infty.
\]
Here, $S(E)$ is the supersymmetric algebra, considered as the set of constant coefficient differential operators on $L(\csv E)$, and $f(D;x)$ is the value of the superfunction $Df$ at the point $x$.

Moreover, $\TDi0{\csv E}$ denotes the space of continuous linear functionals on $\Sw0{\csv E}$; the elements thereof are called \Define{tempered superdistributions}. For more details, consult \ref{app:schwartz}.

\begin{Prop}
  The morphism $Q$ induces a morphism $L(\csv V)\to L(\csv W)$ of \emph{cs} manifolds, and the pullback along this morphism induces a continuous linear map $Q^\sharp:\Sw0{\csv W}\to\Sw0{\csv V}$. 
\end{Prop}

\begin{proof}
  For a \emph{cs} manifold $S$, let $x\in_SL(\csv V)$. Then 
  \[
    x=\kbordermatrix{
 & n & p|q  \\
n & 0 & a'   \\
p|q & a & 0 }
\]
implies
\[
x^2=P(x)1
=\kbordermatrix{
 & n & p|q  \\
n & a'a & 0   \\
p|q & 0 & aa' },\quad
Q(x)=\kbordermatrix{
 & n & p|q  \\
n & 0 & 0   \\
p|q & 0 & aa' }.
  \]

  The underlying values of the entries of $x$ satisfy $a_0=(a'_0)^*$, so $Q(x)_0$ is an $S_0$-valued point of $L(\csv W)_0=\csvr W$. Thus, $Q$ descends to a morphism of \emph{cs} manifolds $L(\csv V) \to L(\csv W)$, proving the first statement. 
  
  As for the second statement, the partial derivatives of the components of $Q$ are polynomials, so it is sufficient to prove that for any $f\in\Sw0{\csv W}$, $k\in\nats$, 
  \[
    \sup\nolimits_{x\in\csvr V}\Abs1{(1+\Norm0x^2)^kQ^\sharp(f)(x)}<\infty,
  \]
  and that this quantity, as a function of $f$, is a continuous seminorm on $\Sw0{\csv U}$. Here, $\Norm0x$ denotes the operator norm of $x\in\csvr V$. 

  Clearly, $\Norm0{x^2}\sle\Norm0x^2$. Since $x$ is a Hermitian matrix and therefore normal, we have the spectral identity 
  \[
      \Norm0x=\vrho(x)=\inf_n\Norm0{x^n}^{1/n}\sle\Norm0{x^2}^{1/2},
  \]
  where $\vrho$ denotes the spectral radius. Thus, with the notation above, it follows that 
  \[
    \Norm0x^2=\Norm0{x^2}=\Norm0{a'a}\Norm0{aa'}=\Norm0{aa'}^2, 
  \]
  where we have used the fact that $a'=a^*$ and thus $\Norm0{a'a}=\Norm0{aa'}=\Norm0a^2$ by the same argument as for $x$. Therefore, the statement is immediate. 
\end{proof}

\begin{Cor}[qsharp]
  There is a continuous linear map $Q_\sharp:\TDi0{\csv V}\to\TDi0{\csv W}$, given by 
  \[
    \Dual1{Q_\sharp(u)}f\defi\Dual1u{Q^\sharp(f)}\mathfa u\in\TDi0{\csv V}\,,\,f\in\Sw0{\csv W}.
  \]
\end{Cor}

\subsection{Recollections on Berezin integration}

Below, we make heavy use of Berezin integration for non-compactly supported integrands. In this brief interlude, we explain the formalism used to manipulate them. 

Let us first recall some basic facts. For a \emph{cs} manifold $X$ of $\dim_{cs}X=a|b$, we denote by $\sh Ber_X$ the Berezinian sheaf and by $\Abs0{\sh Ber}_X$ its twist by the orientation sheaf on $X_0$, called the sheaf of \Define{Berezinian densities}. 

If $X$ has global coordinates $(x,\xi)$, we have a distinguished basis $\Abs0{D(x,\xi)}$ of $\Abs0{\sh Ber}_X$. For $\omega=\Abs0{D(x,\xi)}\,f$ with $f=\sum_If_I\xi^I$, we let
\begin{equation}\label{eq:fibint-def}
  \fibint[X_0]X\omega\defi\Abs0{dx_0}\,f_{\{1,\dotsc,b\}}.  
\end{equation}
This is a density on $X_0$, and the Berezin integral of $\omega$ is defined by 
\begin{equation}\label{eq:berint-def}
  \int_X\omega\defi\int_{X_0}\Bracks3{\fibint[X_0]X\omega}  
\end{equation}
whenever this density is absolutely integrable. 

Unless $\omega$ is compactly supported, the definition depends on the chosen coordinates. Thus, to make sense of Berezin integrals for non-compactly supported $\omega$, one needs to fix some further information. As explained in Ref.~\cite{ahp}, this extra datum is that of a \Define{retraction}. We briefly recollect the basics. 

By definition, a \Define{retraction} of $X$ is a morphism $r:X\to X_0$ left-inverse to the canonical embedding $X_0\to X$---that is, the value of $r^\sharp(f)$ at any $x\in X_0$ is $f(x)$. Then, if one restricts attention to coordinate systems \Define{adapted to} $r$, the left-hand side of Equation \eqref{eq:berint-def} depends only on $r$. Here, the coordinate system $(x,\xi)$ is called adapted to $r$ if $x=r^\sharp(x_0)$. 

We say that $\omega$ is \Define{absolutely integrable} with respect to $r$ if the density defined locally by Equation \eqref{eq:fibint-def} is absolutely integrable over $X_0$. In this case, we define the integral of $\omega$, denoted $\int_X^r\omega$, by Equation \eqref{eq:berint-def}.

The theory of such integrals is fully developed in Ref.~\cite{ahp}. We give some details of the relative version thereof in \ref{app:super-int}. 

\subsection{Statement of the Theorem}

In this subsection, we state our Main Theorem as an identity of two generalised superfunctions and give a proof, up to the explicit computation of the Laplace transform of the right-hand side. For the reader's convenience, we give the necessary definitions here, but defer parts of the proofs to later sections of the paper. 

We begin by describing the left-hand side of the superbosonisation identity. On $U$ and the subspaces $V$ and $W$, we consider the supertrace form $\str(uu')$, which is positive definite on $\smash{\csvr V}$. We normalise the Lebesgue density $\Abs0{dv_0}$ on $\csvr V$ by fixing the volume of the unit cube with respect to an orthonormal basis to one. 

We choose an orientation on $V_\odd$ and identify the dual space $V^*$ with $V$ \via the supertrace form. The \Define{standard Berezinian density} $\Abs0{Dv}$ on $L(\csv V)$, introduced in \thmref{Def}{std-ber}, is then characterised by: 
\[
  \int_{L(\smash{\csv V})}\Abs0{Dv}\,f(v)=\int_{\smash{\csvr V}}\Abs0{dv_0}\,\frac\partial{\partial\nu_m}\dotsm\frac\partial{\partial\nu_1}f\mathfa f\in\Gamma_c(\sh O_{\csv V})
\]
where $\nu_1,\dotsc,\nu_m$ is an oriented symplectic basis of $V_\odd$. We consider $\Abs0{Dv}$ as a tempered superdistribution in $\TDi0{\csv V}$ \via
\[
    \Dual0{\Abs0{Dv}}f\defi\int_{L(\smash{\csv V})}\Abs0{Dv}\,f(v)\mathfa f\in\Sw0{\csv V}.
\]
Here, the integral is taken with respect to the standard retraction $r$. This is given by considering $V_{\smash\ev}^*\subseteq V^*$ \via the splitting $V=V_\ev\oplus V_\odd$ and letting 
\[
  r^\sharp(f)\defi f\in V^*\subseteq\Gamma(\sh O_{\csv V})
\]
\fa $f\in V_{\smash\ev}^*$. Compare the beginning of \ref{app:ft-schwartz} for further details.

The following proposition is immediate from the definitions and \thmref{Cor}{qsharp}.

\begin{Prop}[lhs-supp]
  We have $Q_\sharp(\Abs0{Dv})\in\TDi0{\csv W}$ and 
  \[
    \supp Q_\sharp(\Abs0{Dv})\subseteq\overline{\Herm^+(p)}\times\{0\}.
  \]
\end{Prop}

The tempered superdistribution thus defined is the \emph{left-hand side} of the superbosonisation identity. We now describe the right-hand side. 

On general grounds \cites{ah-berezin,a-hchom}, the homogeneous superspace $\Omega=H.1$ admits a non-zero $H$-invariant Berezinian density, unique up to a constant. Due to the special features of this example, we can give an explicit formula. 

Indeed, observe that $\Omega$ is a locally closed \emph{cs} submanifold of $\smash{L(W)}$ with purely even codimension. Consider the standard coordinates 
\[
  Z=\kbordermatrix{%
      & p & q \\
    p & z & \zeta\\
    q & \omega & w
  }
\]
on $\smash{L(W)}$. For any choice of even coordinates $x$ among the components of $z,w$ that define a local coordinate system on $\Omega_0$, $(x,\zeta,\omega)$ is a local coordinate system on $\Omega$. This defines a retraction of $\Omega$, which we also call \Define{standard}, and a system of adapted coordinates for this retraction. 

In particular, $D(\zeta,\omega)$ is a well-defined relative Berezinian (density) on $\Omega$ over $\Omega_0$, with respect to the standard retraction. Denote by $\Abs0{dz}$ the Lebesgue density on $\Herm(p)$, and by $\Abs0{dw}$ the normalised invariant density on $\mathrm U(q)$. We set 
\begin{equation}\label{eq:invber-def}
  D\mu(Z)\defi\frac{\Abs0{dz}\,\Abs0{dw}}{\Abs0{\det z}^p}\,D(\zeta,\omega)\,\det(z-\zeta w^{-1}\omega)^q\det(w-\omega z^{-1}\zeta)^p.  
\end{equation}
Then $\mu$ is $H$-invariant; the proof is deferred to \thmref{Prop}{inv-ber} below. In what follows, we write $y$ for $S$-valued points of $\Omega$, and set $\Abs0{Dy}\defi D\mu(Z)$.

In what follows, let $n\sge p$. We define $T_n$ by 
$$
\Dual0{T_n}f\defi\int_\Omega \Abs0{Dy}\,\Ber0y^n f(y)
$$
for any (holomorphic) superfunction $f$ on the open subspace $T(\gamma)$ of $L(W)$ whose underlying set is the tube
\begin{equation}\label{eq:tube}
  T(\gamma_0)\defi\Parens1{\Herm^+(p)\times\Herm(q)}+i\csvr W=T(\Herm^+(p))\times\cplxs^{q\times q}
\end{equation}
with Paley--Wiener estimates along $T(\Herm^+(p))\defi\Herm^+(p)+i\Herm(p)$, \ie
\begin{equation}\label{eq:pw-tube0}
  \sup\nolimits_{z\in T(\Herm^+(p))}\Abs1{e^{-R\Norm0{\Im z}}(1+\Norm0z)^Nf(D;z,w)}<\infty
\end{equation}
for any $D\in S(W)$, $N\in\nats$, $w\in\cplxs^{q\times q}$, and some $R>0$. Here, $2i\Im z\defi z-z^*$.

The following proposition is a direct consequence of \thmref{Prop}{conv}, whose proof is given below in Section \ref{sec:lap-conical}. 

\begin{Prop}[rhs]
  When $n\sge p$, the functional $T_n$ is well-defined and continuous on the space of all $f\in\Gamma(\sh O_{T(\gamma)})$ satisfying the estimate in Equation \eqref{eq:pw-tube0} \fs $R>0$ and any $D\in S(W)$, $w\in\cplxs^{q\times q}$, and $N\in\nats$.
\end{Prop}

In particular, $T_n$ may be considered as an element of $\PWFn0{\csv W}$, the topological dual space of the Paley--Wiener space $\PW0{\csv W}$. By \thmref{Def}{pw-def} below, the latter is given as the set of superfunctions $f\in\Gamma(\sh O_W)$ holomorphic on all of $L(W)$, satisfying estimates
\begin{equation}\label{eq:pw-tube}
  \sup\nolimits_{y\in W_\ev}\Abs1{e^{-R\Norm0{\Im y}}(1+\Norm0y)^Nf(D;y)}<\infty
\end{equation}
for any $D\in S(W)$, $N\in\nats$, and some $R>0$. Here, $2i\Im y\defi y-y^*$.

The functional $T_n$ is the \Define{right-hand side} of the superbosonisation identity. For $q=0$, it coincides with the unweighted \emph{Riesz distribution} for the parameter $n$, \vq Ref.~\cite{FK94}. Thus, it may also be called the \Define{Riesz superdistribution}. 

We now state the functional analytic version of the superbosonisation identity. 

\begin{Th}[sbos][superbosonisation identity]
  Assume that $n\sge p$. Then we have that $T_n\in\TDi0{\csv W}$, and 
  \[
    Q_\sharp(\Abs0{Dv})=\frac{\sqrt{\pi}^{np}}{\Gamma_\Omega(n\mathds1)}\cdot T_n,
  \]
  where the finite constant $\Gamma_\Omega(n\mathds1)>0$ is determined in \thmref{Th}{gamma}.

  Moreover, these generalised superfunctions extend as continuous functionals to the space of all superfunctions $f\in\Gamma(\sh O_{\csv W})$ that satisfy Schwartz estimates along $\smash{\Herm^+(p)}$, \ie 
  \[
    \sup\nolimits_{z\in\Herm^+(p)}\Abs1{(1+\Norm0z)^Nf(D;z,w)}<\infty
  \]
  for all $w\in\Herm(q)$, $N\in\nats$, and $D\in S(W)$.
\end{Th}

The first assertion of the theorem means that $T_n$ lies in the image of the continuous injection $\TDi0{\csv W}\to\PWFn0{\csv W}$. Such an injection exists, since the Paley--Wiener space $\PW0{\csv W}$ is a dense subspace of the Schwartz space $\Sw0{\csv W}$, as follows from the Paley--Wiener theorem (\thmref{Prop}{PW}).


\thmref{Th}{sbos} immediately implies the following explicit formula, which is a precise statement of the identity proved in Ref.~\cite{LSZ08}.

\begin{Cor}
  Let $f$ be a (holomorphic) superfunction on the open subspace of $L(W)$ whose underlying set is the tube $T(\gamma_0)$ from Equation \eqref{eq:tube}, satisfying the estimate in Equation \eqref{eq:pw-tube} \fs $R>0$ and any $D\in S(W)$, $w\in\cplxs^{q\times q}$, and $N\in\nats $. Then 
  \[
    \int_{L(\csv V)}\Abs0{Dv}\,Q^\sharp(f)(v)=\frac{\sqrt{\pi}^{np}}{\Gamma_\Omega(n\mathds1)}\int_\Omega \Abs0{Dy}\,\Ber0y^n f(y).
  \]
\end{Cor}

The proof of \thmref{Th}{sbos} relies on the theory of the Laplace transform of generalised superfunctions, which is developed in \ref{app:super-lap}. 
Before we embark upon the proof, let us briefly summarise the key features. 

The Laplace transform is defined for functionals on either of the super-vector spaces $\mathscr S=\Sw0{\csv W}$ and $\mathscr Z=\PW0{\csv W}$. We call these $E=\mathscr S,\mathscr Z$ \Define{test spaces} and their topological duals $E'$ spaces of \Define{generalised superfunctions}. 

As is the case classically, the Laplace transform of $\mu\in E'$ at $x+iy$ is the Fourier transform of $e^{-\Dual0x\cdot}\mu$, as a generalised superfunction of $y$. That is, we set 
\[
  \LT(\mu)(x)\defi\mathcal F\Parens1{e^{-\Dual0x\cdot}\mu}.
\]
In the Appendix, we make sense of this for $S$-valued points $x$. 

The Laplace transform at $x$ is a holomorphic superfunction of $x+iy$ on the tube $T(\gamma)\defi\gamma+iL(\csv W)$ where $\gamma=\gamma_\mathscr S^\circ(\mu)$ is the maximal domain on which $e^{-\Dual0x\cdot}\mu$ is a superfunction of $x$ with values in $\mathscr S'$, that is 
\[
  e^{-\Dual0x\cdot}\mu\in\Gamma(\sh O_S)\,{\Hat\otimes}_\pi\,\TDi0{\csv W}
\]
for every $S$-valued point of $\gamma$. Compare \ref{app:vectval}, Equation \eqref{eq:expins}, and the surrounding remarks for details. The thus defined holomorphic superfunction is called the \Define{Laplace transform of $\mu$} and denoted by $\LT(\mu)$. It entirely determines $\mu$ (\thmref{Th}{lap-hol}).

If $\mu$ is already tempered and supported in a pointed cone $\gamma$, then the domain $\gamma_\mathscr S^\circ(\mu)$ contains the open dual cone $\check\gamma$ (\thmref{Cor}{supp-lap}).

\begin{proof}[\protect{Proof of \thmref{Th}{sbos}}]
The idea of the proof is to compute the Laplace transforms of both sides of the equation and to compare the results. In view of the injectivity of the Laplace transform (\thmref{Th}{lap-hol}), this will show the equality of the two functionals. To make this argument rigorous, we have to show that the domains of holomorphy for the Laplace transforms intersect. Throughout the proof, we will consider the \emph{cs} manifolds $L(\csv V)$ and $L(\csv W)$ as embedded in $L(\csv U)$. 

We start by considering the left-hand side. Since the cone $\Herm^+(p)$ is self-dual, it is immediate from \thmref{Prop}{lhs-supp} and \thmref{Cor}{supp-lap} that
\[
  \gamma_{\mathscr S}^\circ\Parens1{Q_\sharp(\Abs0{Dv})}_0\supseteq\gamma_0=\Herm^+(p)\times\Herm(q).
\]
Therefore, the Laplace transform of $Q_\sharp(\Abs0{Dv})$ is defined and holomorphic on $T(\gamma)$, where $\gamma$ is the open subspace of $L(\csv W)$ whose underlying set is $\gamma_0$. For $x\in_ST(\gamma)$, we compute
$$
  \LT(Q_\sharp(\Abs0{Dv})(x)=\int_{L(\csv V)} \Abs0{Dv}\,Q^\sharp (e^{-\str(x\cdot)})(v) = \int_{L(\csv V)} \Abs0{Dv}\,e^{-\str(xQ(v))},
$$
arguing that if the integral converges absolutely, it must equal the Laplace transform.

For each $x\in W$, we define a linear map $\phi_x:V \to V$ by
$$
\phi_x\begin{Matrix}10&a'\\a&0\end{Matrix} \defi \begin{Matrix}10&a'x\\xa&0\end{Matrix},
$$
for $v=\begin{Matrix}00&a'\\a&0\end{Matrix}\in V$. By the cyclicity of the supertrace, notice that 
\[
  2\str(xQ(v))=2\str(xaa')=\str(a'xa)+\str(aa'x)=\str(v\phi_x(v)).
\]
Now, let $\gamma^+\subseteq\gamma$ be the open subspace corresponding to 
\[
  \gamma_0^+\defi\Herm^+(p)\times\Herm^+(q)
\]
and let $x\in_S\gamma^+$. Then we may choose $\smash{x^{1/2}\in_S\gamma^+}$ \scth $(x^{1/2})^2=x$ and make a change of coordinates $v \mapsto \smash{\phi_{x^{-1/2}}(v)}$. The Berezinian of this coordinate transformation is given by $\Ber0x^{-n}$ and $\str(v\phi_x(v)) \mapsto\str(v^2)$. So,
$$
\LT(Q_\sharp(\Abs0{Dv}))(x)= \Ber0x^{-n}\int_{L(\csv V)} \Abs0{Dv}\,e^{-\frac12\str(v^2)}.
$$
By holomorphicity, both sides of this equation coincide on $T(\gamma)\cap\GL(p|q,\cplxs)$.

To determine the Gaussian integral, pick coordinates $(x,\xi,\eta)$ on $L(\csv V)$:
$$
\kbordermatrix{
 & n & p & q  \\
n & 0 & \overline{x_{ji}} &  \eta_{ji} \\
p & x_{ij} & 0 & 0\\
q & \xi_{ij} & 0 & 0},
$$
so that 
$$
-\frac12\str(v^2)=-\tr(xx^*)+\tr(\xi\eta) =-\sum_{i=1}^p\sum_{j=1}^n |x_{ij}|^2 +\sum_{k=1}^q\sum_{j=1}^n\xi_{kj}\eta_{jk}.
$$
The Berezin integral is performed by picking the degree $nq$ term in the expansion of the exponential function for $e^{\mathrm{tr}(\xi\eta)}$, which is just $1$ (for a suitable choice of signs).

Further, the remaining integral is just $np$ copies of the Gaussian integral, which contributes $\smash{\sqrt{\pi}^{np}}$. Therefore, we find 
$$
\LT(Q_\sharp(\Abs0{Dv}))(x)= \sqrt{\pi}^{np}\,\Ber0x^{-n},
$$
for any $x\in_ST(\gamma)\cap\GL(p|q,\cplxs)$.

On the other hand, from the the definition of $\gamma_\mathscr S^\circ(T_n)$ and \thmref{Prop}{rhs}, we see that $\gamma\subseteq\gamma_\mathscr S^\circ(T_n)$. Hence, the Laplace transform of $T_n$ is defined and holomorphic on $T(\gamma)$, by \thmref{Th}{lap-hol}. Moreover, by \thmref{Cor}{sw-ext}, $T_n$ is already a tempered superdistribution. It remains to compute $\LT(T_n)$.

As a consequence of \thmref{Cor}{delta-lap}, whose proof is deferred to Section \ref{sec:lap-conical}, we have the identity 
\[
  \LT(T_n)(x)=\int_\Omega \Abs0{Dy}\,\Ber0y^n e^{-\str(xy)}=\Gamma_\Omega(n\mathds1)\,\Ber0x^{-n},
\]
for all $x\in_S T(\gamma)$ with all principal minors of $x$ invertible. Here, $\mathds1=(1,\dotsc,1)$ and $\Gamma_\Omega$ is a meromorphic function of $p+q$ indeterminates. It is explicitly computed in \thmref{Th}{gamma}, whose proof is also deferred. In particular, when $n\sge p$, one sees that the constant $\Gamma_\Omega(n\mathds1)$ is a non-zero positive number. 

Comparing the outcome of the two computations and the domains of holomorphy, the result follows from the injectivity of the Laplace transform as stated in \thmref{Th}{lap-hol}.
\end{proof}

\section{Laplace transforms of conical superfunctions}\label{sec:lap-conical}

As we have seen in Section \ref{sec:two}, the superbosonisation identity reduces to computing the Laplace transform of both sides. Whereas for the left hand side, this amounts to the evaluation of a standard Gaussian integral, the computation on the right hand side is more intricate. 

Following the procedure from the classical case, where $q=0$, we compute the Laplace transform of certain, more general conical superfunctions. 
The outcome for $q>0$ is more complicated than in the classical case, where the Laplace transform has poles, but no zeros; this is quite different for $q>0$. 

\subsection{The conical superfunctions} 

We introduce the basic objects of this section, the conical superfunctions. These are a natural generalisation, to the superspace $\Omega$, of the conical polynomials encountered in the theory of Riemannian symmetric spaces. Our main reference to the subject will be the book of Faraut--Kor\'anyi \cite{FK94}, which contains a beautiful, elementary, and self-contained account for the special class of such spaces consisting of the symmetric cones.

Let $N^+$ be the closed subsupergroup of $K_\cplxs$ whose functor of points is defined by 
\[
  N^+(S)\defi\left\{\begin{Matrix}1A&0\\0&D\end{Matrix}\in_SK_\cplxs\left|A=\trimat(*,1,1,0),D=\trimat(0,1,1,*)\right.\right\}.
\]
That this functor is represented by a closed, connected complex analytic subsupergroup of $K_\cplxs$ is immediate from the implicit function theorem.

Similarly, we define $T_\cplxs$ to be the complex supergroup representing the functor
\[
  T_\cplxs(S)\defi\Set4{\begin{Matrix}1A&0\\0&D\end{Matrix}\in_SK_\cplxs}{A=\diag(a),D=\diag(d), a,d\in_S(\cplxs^\times)^{p+q}}.
\]
Then $T_\cplxs$ normalises $N^+$, and the subsupergroup $B\defi T_\cplxs N^+$ of $K_\cplxs$ generated by $T_\cplxs$ and $N^+$ is again closed and connected. Its complex super-dimension is 
\[
  \dim_\cplxs B=(p+q)|0+\dim_\cplxs W.
\]

As a direct consequence of these definitions, we obtain the following lemma.

\begin{Lem}
  The orbit $K_\cplxs.1=\GL(p|q,\cplxs)$ is open in $W$, and the action of $B$ on $K_\cplxs$ also admits an open orbit, namely, $B.1$. Here, $1$ denotes the identity matrix, considered as an ordinary point of $W$.
\end{Lem}




Let $B_p$ and $B_q$ denote the Borel subgroups of $\GL(p,\cplxs)$ and $\GL(q,\cplxs)$, respectively, given by lower triangular matrices. Denote their opposite Borels with bars. Then the underlying Lie group of $B$ is $B_0=B_p\times B_q\times\overline B_p\times\overline B_q$. Hence, $B.1$ is the open subspace of $W$ corresponding to the open set $B_p\overline B_p\times B_q\overline B_q$. This justifies calling $B.1$ the \Define{big cell} of $K_\cplxs.1=\GL(p|q,\cplxs)$; it also shows that $(B.1)_0$ is Zariski open, so it makes sense to speak of regular superfunctions on $B.1$.

\bigskip\noindent
We will now define a family of rational superfunctions $\Delta_1,\dotsc,\Delta_{p+q}$ which in some sense are fundamental (relative) invariants. Here, the superalgebra $\cplxs(W)$ of \Define{rational superfunctions} is defined to be $\cplxs(W_\ev)[W_\odd]=\cplxs(W_\ev)\otimes\bigwedge(W_\odd)^*$, where $\cplxs(W_\ev)$ is the algebra of rational superfunctions on $W_\ev$. Each $f\in\cplxs(W)$ may be considered as a superfunction on an open subspace of $L(W)$ in an obvious fashion. 

For any $Z=(Z_{ij})\in_S W=\ger{gl}(p|q,\cplxs)$ and $1\sle k,\ell\sle p+q$, we consider 
\[
  [Z]_{k\ell}\defi(Z_{ij})_{1\sle i\sle k,1\sle j\sle \ell}\nd [Z]_k\defi[Z]_{kk},
\]
so that $[Z]_k$ is the $k$th principal minor of $Z$. Whenever $[Z]_k$ is invertible, we define
\[
  \Delta_k(Z)\defi\Ber0{[Z]_k}.
\]
This uniquely determines a rational superfunction $\Delta_k\in\cplxs(W)$. We also consider the (even) characters $\chi_k:T_\cplxs\to\cplxs^\times$, defined by 
\[
  \chi_k(t)\defi\frac{\prod_{j=1}^{\min(k,p)}a_j^{-1}d_j}{\prod_{j=p+1}^{\min(k,p+q)}a_j^{-1}d_j}\mathfa t=\begin{Matrix}1\diag(a)&0\\0&\diag(d)\end{Matrix}\in_ST_\cplxs.
\]
Note that $\chi_k(t)=\Delta_k(t^{-1}.1)$. In view of the isomorphism $T_\cplxs\cong B/N^+$, we may consider $\chi_k$ as a character of $B$.

\begin{Prop}[conical]
  Let $1\sle k\sle p+q$. Then the rational superfunction $\Delta_k$ is regular and invertible on the big cell $B.1$, and 
  \begin{equation}\label{eq:deltak-relinv}
    \Delta_k(b^{-1}.Z)=\chi_k(b)\Delta_k(Z),  
  \end{equation}
  for any $b\in_SB$ and $Z\in_SB.1$.
\end{Prop}

In the \emph{proof} of the proposition, we will need the following lemma.

\begin{Lem}[bigcell]
  Let $Z\in_SW$. Then $Z\in_SB.1$ if and only if $[Z]_k$ is invertible for every $1\sle k\sle p+q$.
\end{Lem}

\begin{proof}
  We write 
  \[
    Z=\kbordermatrix{%
          & \mathclap{p+q-1} & 1\\
      p+q-1   & A   & B\\
      1       & C   & D
    },
  \]
  where $A=[Z]_{p+q-1}$. If $A$ is invertible, we may perform a block decomposition of $Z$ as as follows:
  \[
    Z=
    \begin{Matrix}1
      1&0\\
      CA^{-1}&1
    \end{Matrix}
    \begin{Matrix}1
      A&0\\
      0&D-CA^{-1}B
    \end{Matrix}
    \begin{Matrix}1
      1&A^{-1}B\\
      0&1
    \end{Matrix}.
  \]
  In this decomposition, $D-CA^{-1}B$ is automatically an invertible element of $\Gamma(\sh O_{S,\ev})$ whenever $Z=[Z]_{p+q}$ is invertible.

  Provided that all principal minors are invertible, we may, by replacing $Z$ by $A$, continue with this procedure to arrive at a decomposition of the form $Z=ldu$ where $l$ is strictly lower triangular, $d$ is diagonal, and $u$ is strictly upper triangular. Then $b\defi (ld,u^{-1})\in_SB$, and $b.1=Z$. This shows that the set of all $Z\in_SW$ \scth all principal minor are invertible is contained in $(B.1)(S)$.

  Given any $Z\in_SW=\ger{gl}(p|q,\cplxs)$ \scth $[Z]_k$ is invertible, write
  \[
    Z=\kbordermatrix{%
          & k   & \mathrlap{p+q-k}\\
      k     & A   & *\\
      p+q-k & *   & *
    },
  \]
  where $A=[Z]_k$. Define integers 
  \[
    m\defi \min(p,k)\nd n\defi\min(0,p+q-k),
  \]
  and let $\alpha,\delta\in_S\GL(m|n,\cplxs)$. Then 
  \begin{equation}\label{eq:minor-trafo}
    \begin{Matrix}1
      \alpha  & 0\\
      *   & *
    \end{Matrix}
    \begin{Matrix}1
      A     & *\\
      *   & *
    \end{Matrix}
    \begin{Matrix}1
      \delta  & *\\
      0     & *
    \end{Matrix}
    =
    \begin{Matrix}1
      \alpha A\delta & *\\ *&*
    \end{Matrix}, 
  \end{equation}
  where the upper left block of the latter matrix is invertible. 

  In particular, the open subspace of $W$ whose $S$-valued points are the $Z\in_SW$ with all principal minors invertible is invariant under the action of $B$. But this already shows the equality. 
\end{proof}

\begin{proof}[\protect{Proof of \thmref{Prop}{conical}}]
  \thmref{Lem}{bigcell} shows that the functions $\Delta_k$ are regular on $B.1$. We have already noted the identity 
  \[
    \Delta_k(t^{-1}.1)=\chi_k(t)
  \]
  for $t\in_ST_\cplxs$. Since, as we have shown in the course of the proof of \thmref{Lem}{bigcell}, every $Z\in_SB.1$ may be written as $Z=b.1$ for some $b\in_SB$, it will be sufficient to show that $\Delta_k$ is $N^+$-invariant. 


  So, let $n=\begin{Matrix}0n'&0\\0&n''\end{Matrix}\in_SN^+$ and $1\sle k\sle p+q$. Denote the upper left $k\times k$ block of $n'$ resp.~$n''$ by $\alpha$ and $\delta$, respectively. Then $\alpha$ and $\delta$ are invertible. Moreover, Equation \eqref{eq:minor-trafo} shows that 
  \[
    [n^{-1}.Z]_k=\alpha^{-1}[Z]_k\delta,
  \]
  and hence, that 
  \[
    \Delta_k(n^{-1}.Z)=\Ber0\alpha^{-1}\Delta_k(Z)\Ber0\delta=\Delta_k(Z),
  \]
  due to the tridiagonal nature of $\alpha$ and $\delta$. This proves the claim.
\end{proof}

\begin{Rem}
  In case $q=0$, the functions $\Delta_1,\dotsc,\Delta_p$ are known from the theory of Jordan algebras. The statement corresponding to \thmref{Prop}{conical} is to be found, \eg in \cite{FK94}*{Proposition VI.3.10}. In this case, the $\Delta_k$ are all polynomials. 

  In general, this continues to hold for $\Delta_k$, $1\sle k\sle p$. As for the other $\Delta_k$, $k>p$, they are certainly regular on the larger domain where only $[Z]_{p+1},\dotsc,[Z]_{p+q}$ are invertible. By density of the big cell in $W$, Equation \eqref{eq:deltak-relinv} continues to hold there. 
\end{Rem}

The rational characters of $B$ (equivalently, of $T_\cplxs$) are given exactly by the superfunctions
\[
  \chi_\mathbf m\defi\chi_1^{m_1-m_2}\dotsm\chi_{p+q-1}^{m_{p+q-1}-m_{p+q}}\cdot\chi_{p+q}^{m_{p+q}},
\]
where $\mathbf m=(m_1,\dots,m_{p+q})\in\ints^{p+q}$. Explicitly, for $t=(\diag(a),\diag(d))\in_ST_\cplxs$:
\[
  \chi_\mathbf m(t)=\frac{\prod_{j=1}^p(a_j^{-1}d_j^{\vphantom{m_j}})^{m_j}}{\prod_{j=1}^q(a_{p+j}^{-1}d_{p+j}^{\vphantom{m_{p+j}}})^{m_{p+j}}}.
\]

If we define $\Delta_\mathbf m$ by 
\[
  \Delta_\mathbf m\defi\Delta_1^{m_1-m_2}\dotsm\Delta_{p+q-1}^{m_{p+q-1}-m_{p+q}}\Delta_{p+q}^{m_{p+q}},
\]
then $\Delta_\mathbf m$ is a regular superfunction on the big cell by \thmref{Prop}{conical}, and by the same token, we have
\[
  \Delta_\mathbf m(b^{-1}.Z)=\chi_\mathbf m(b)\Delta_\mathbf m(Z)
\]
\fa $b\in_SB$ and $Z\in_SB.1$.

\begin{Def}[conical][conical superfunctions]
A rational function $f\in\cplxs(W)$ is \Define{conical} if its domain of definition is $B$-invariant and there exists $\mathbf m\in\ints^{p+q}$ \scth 
\[
  f(b^{-1}.Z)=\chi_\mathbf m(b)f(Z)
\]
\fa $b\in_SB$ and all $S$-valued points $Z$ of the domain of definition of $f$.
\end{Def}

\begin{Lem}[conical-exhaustion]
  Let $f$ be a conical superfunction. Then $f$ is proportional to $\Delta_\mathbf m$ for some multi-index $\mathbf m\in\ints^{p+q}$.
\end{Lem}

\begin{proof}
  The domain of definition of $f$ is Zariski open and dense, as is the big cell $B.1$. Hence, $f$ is regular on $B.1$ and uniquely determined by its values $f(Z)$ for any $Z\in_SB.1$ (and any $S$). Then
  \[
    f(b.1)=\chi_\mathbf m(b^{-1})f(1)=f(1)\Delta_\mathbf m(b.1)
  \]
  for any $b\in_SB$, which already shows that $f=f(1)\Delta_\mathbf m$.
\end{proof}

In one instance below, it will be useful to have an alternative parametrisation of the `boson-boson sector' of $B$. To that end, define 
\[
  \tau_j(u)\defi  
  \kbordermatrix{
        &   j-1   &   1   &   p-j\\
      j-1 & 0     & 0   & 0\\
      1   &   0     &   0   &   u\\
      p-j & 0     &   0   & 0   
    }
\]
for any $1\sle j\sle p$ and $u\in\cplxs^{1\times(p-j)}$. Then 
\[
  e^{\tau_j(u)}=    
  \kbordermatrix{
        &   j-1   & 1   &   p-j\\
      j-1 & 1     &   0   & 0\\
      1   &   0     &   1   &   u\\
      p-j & 0     &   0   & 1   
    }\in\overline B_p.
\]

Let $u$ be a Hermitian matrix. We denote the diagonal entries by $u_j=u_{jj}$ and the rows of the upper triangle by 
\[
  u^{(j)}\defi(u_{j,j+1},\dotsc,u_{jp})\mathfa1\sle j<p.
\]
If $u_1,\dotsc,u_p>0$, the we define 
\[
  t(u)\defi\delta_1(u_1)e^{\tau_1(u^{(1)})^*}\delta_2(u_2)e^{\tau_2(u^{(2)})^*}\dotsm e^{\tau_{p-1}(u^{(p-1)})^*}\delta_p(u_p)\in\GL(p,\cplxs),
\]
where we set 
\[
  \delta_j(\lambda)\defi(\underbrace{1,\dotsc,1}_{j-1},\lambda,\underbrace{1,\dotsc,1}_{p-j}).
\]

Considering as usual $\GL(p,\cplxs)$ as a subgroup of $H$, \via the map 
\[
  a\mapsto\left(
  \begin{Matrix}[0]0
    a&0&0&0\\
    0&1&0&0\\
    0&0&(a^*)^{-1}&0\\
    0&0&0&1
  \end{Matrix}\right),
\]
we see immediately that for any $Z\in_S B.1$ and $\mathbf m\in\ints^{p+q}$, we have
\[
  \Delta_\mathbf m(t(u).Z)=\Delta_\mathbf m(Z).
\]

On the other hand, we obtain coordinates on $\Herm^+(p)$ in this fashion, as explained in the proposition below, which translates \cite{FK94}*{Proposition VI.3.8} to this special case. We give the direct proof using matrices both for the reader's convenience, and because we will need to refer to it later on.

\begin{Prop}[herm-param]
  Let $z\in\Herm^+(p)$. There is a unique $u\in\Herm^+(p)$, with diagonal entries $u_1,\dotsc,u_p>0$, \scth $z=t(u).1$, and it is determined by 
  \[
    z_{jj}=u_j^2+\sum_{k=1}^{j-1}\Abs0{u_{kj}}^2\nd z_{jk}=u_ju_{jk}+\sum_{\ell=1}^{j-1}u_{j\ell}u_{\ell k}
  \]
  for all $1\sle j\sle p$ and $j<k$, $k\sle p$.
\end{Prop}

\begin{proof}
  The proof is by induction on $p$. For $p=1$, the statement is trivial. For general $p$, 
  \[
    z=\delta_1(u_1)e^{\tau_1(u^{(1)})^*}
    \begin{Matrix}1
      1&0\\0&z'
    \end{Matrix}=
    \begin{Matrix}1
      u_1^2&u_1u^{(1)}\\
      *&z'+(u^{(1)})^*u^{(1)}
    \end{Matrix},
  \]
  where $u_1$ and $u^{(1)}$ are determined by 
  \[
    z_{11}=u_1^2\nd z_{1k}=u_1u_{1k}
  \]
  \fa $k>1$. By the inductive hypothesis, we have $z'=t(u').1=(z_{jk})_{2\sle j,k\sle p}$ for $u'=(u_{jk})_{2\sle j,k\sle p}$, where $z'$ and $u'$ are related by 
  \[
    z_{jj}'=u_j^2+\sum_{k=2}^{j-1}\Abs0{u_{kj}}^2\nd z_{jk}'=u_ju_{jk}+\sum_{\ell=2}^{j-1}u_{j\ell}u_{\ell k}
  \]
  for all $2\sle j\sle p$ and $j<k$, $k\sle p$. The assertion follows by noting simply that
  \[
    \Parens1{(u^{(1)})^*u^{(1)}}_{jk}=\overline{u_{1j}}u_{1k}=u_{j1}u_{1k},
  \]
  where we take the entries of this matrix to be indexed over the set $\{2,\dotsc,p\}$.
\end{proof}

\subsection{The invariant Berezinian}\label{subs:ber}

Recall the definition of the homogeneous superspace $\Omega=H.1$ from Section \ref{sec:two}. The underlying manifold is 
\[
  \Omega_0=\Herm^+(p)\times\mathrm U(q),
\]
where 
\[  
  \Herm^+(p)=\GL(p,\cplxs)/\mathrm U(p)
\]
is the cone of positive definite Hermitian $p\times p$ matrices. We observe that the dimension of the \emph{cs} manifold $\Omega$ coincides with the graded dimension of the complex supermanifold $L(W)$. In particular, since $\Omega_0$ is contained in $(B.1)_0$, it follows that $\Omega$ is a subspace of the \emph{cs} manifold associated with the complex supermanifold $B.1$.

Recall the definition of the Berezinian density $\mu$ from Equation \eqref{eq:invber-def}.

\begin{Prop}[inv-ber]
  The Berezinian density $\mu$ is $H$-invariant. 
\end{Prop}

We divide the \emph{proof} of this statement into several lemmata which will also be useful below when computing Laplace transforms. Recall the notion of nilpotent shifts from Appendix \ref{app:points}.

\begin{Lem}[uq-shift]
  Let $n$ be a nilpotent shift for $\mathrm U(q)=(\mathrm U(q)\times\mathrm U(q))/\diag\subseteq\cplxs^{q\times q}$. Then 
  \[
    \int_{\mathrm U(q)}\Abs0{dw}\,f(w+n)=\int_{\mathrm U(q)}\Abs0{dw}\,f(w)\frac{\det(w)^q}{\det(w-n)^q}
  \]
  for any smooth function $f$ on $\mathrm U(q)$.
\end{Lem}

\begin{proof}
  The invariant density on $\mathrm U(q)$ coincides, due to the invariance, with the Riemannian density for any invariant Riemannian metric. An application of \thmref{Lem}{shift} shows that we have 
  \[
    \int_{\mathrm U(q)}\Abs0{dw}\,f(w+n)=\int_{\mathrm U(q)}\Abs0{dw}\,fJ_1
  \]
  where $J_1$ is determined by $J_0=1$ and 
  \[
    \frac d{dt}\log J_t(w)=-\ddiv v_n(w-tn)=\tfrac12\tr_{\ger{gl}(q)}R_{(w-tn)^{-1}w}=q\tr((w-tn)^{-1}w),
  \]
  for $w\in\mathrm U(q)$, since for $u\in\ger u(q)$ and $w\in\cplxs^{q\times q}$, we have
  \[
    R_w(u)=\frac d{dt}\exp(tu)w\exp(tu)\Big|_{t=0}=uw+wu.
  \]
  Setting 
  \[
    J_t(w)\defi\frac{\det(w)^q}{\det(w-tn)^q}
  \]
  manifestly solves the equation.
\end{proof}

\begin{Lem}[herm-shift]
  Let $n$ be a nilpotent shift for $\Herm^+(p)=\GL(p,\cplxs)/U(p)\subseteq\cplxs^{p\times p}$. Then 
  \[
    \int_{\Herm^+(p)}\frac{\Abs0{dz}}{\Abs0{\det z}^p}\,f(z+n)=\int_{\Herm^+(p)}\frac{\Abs0{dz}}{\Abs0{\det z}^p}\,f(z)\frac{\det(z)^p}{\det(z-n)^p}
  \]
  for any compactly supported smooth function $f$ on $\Herm^+(p)$.
\end{Lem}

\begin{proof}
  Again \thmref{Lem}{shift} applies, since the invariant density is the Riemannian density, and we have 
  \[
    \int_{\Herm^+(p)}\frac{\Abs0{dz}}{\Abs0{\det z}^p}\,f(z+n)=\int_{\Herm^+(p)}\frac{\Abs0{dz}}{\Abs0{\det z}^p}\,fJ_1
  \]
  where $J_1$ is determined by $J_0=1$ and 
  \begin{align*}
    \frac d{dt}\log J_t(z)&=-\ddiv v_n(z-tn)=\tfrac12\tr_{\ger{gl}(p)}R_{(z-tn)^{-1/2}z(z-tn)^{-1/2}}\\
    &=p\tr((z-tn)^{-1/2}z(z-tn)^{-1/2})=p\tr((z-tn)^{-1}z),
  \end{align*}
  for $z\in\Herm^+(p)$, since for $u\in\Herm(p)$ and $z\in\cplxs^{p\times p}$, we have
  \[
    R_w(u)=\frac d{dt}\exp(tu)w\exp(tu)\Big|_{t=0}=uw+wu.
  \]
  Hence, all we have to do is to solve the same initial value problem as in the proof of \thmref{Lem}{uq-shift}, with $q$ replaced by $p$.
\end{proof}

Now, for $Z\in_SK_\cplxs.1$, set 
\[
  \vrho(Z)\defi\det(z-\zeta w^{-1}\omega)^q\det(w-\omega z^{-1}\zeta)^p.
\]
Observe that 
\begin{equation}\label{eq:vrho-ber}
  \vrho(Z)=\Ber0Z^{q-p}\det(z)^p\det(w)^q=\Delta_{q\mathds1'+(q-p)\mathds1''}(Z)\det(w)^q,
\end{equation}
where we define 
\[
  \mathds1'\defi(1_p,0_q)=(1,\dotsc,1,0,\dotsc,0),\quad\mathds1''\defi(0_p,1_q)=(0,\dotsc,0,1,\dotsc,1).
\]
For $h\in_SK_\cplxs$, we define
\[
  I_h^f(z,w)\defi\int_{\cplxs^{0|p\times q}\oplus\cplxs^{0|q\times p}}D(\zeta,\omega)\,\vrho(Z)f(h.Z)
\]
whenever $f(h.Z)$ is defined. Here, observe that for a purely odd super-vector space $V$ (say), the $\cplxs$-superspaces $L(V)$ and $L(V,0)$ (where $0$ is the unique real form of $V_\ev=0$) coincide. For this reason, when integrating over such a superspace, we will simply write $V$, instead of using the more cumbersome notation.

\begin{Lem}[ferm-invar]
  Let $f\in\Gamma(\sh O_\Omega)$ and $h=\begin{Matrix}0A&0\\0&D^{-1}\end{Matrix}$, where $A=\begin{Matrix}01&0\\\alpha&1\end{Matrix}\in_S\GL(p|q,\cplxs)$ and $D=\begin{Matrix}01&\delta\\0&1\end{Matrix}\in_S\GL(p|q,\cplxs)$. Then for any smooth function $\vphi$ on $\mathrm U(q)$, we have
  \[
    \int_{\mathrm U(q)}\Abs0{dw}\,\vphi(w)I_h^f(z,w)=\int_{\mathrm U(q)}\Abs0{dw}\,\int D(\zeta,\omega)\,\vrho(Z)f(Z)\vphi(w-\alpha\zeta-\omega\delta-\alpha z\delta).
  \]
  Here, the inner integral on the right-hand side is over $\cplxs^{0|p\times q}\oplus\cplxs^{0|q\times p}$.
\end{Lem}

\begin{proof}
  Firstly, note that the conditions set out above imply that $h\in_SH$, so that the statement of the lemma is meaningful. We have
  \[
    h.Z=AZD=
    \begin{Matrix}1
      z&\zeta+z\delta\\
      \omega+\alpha z&w+\alpha\zeta+\omega\delta+\alpha z\delta
    \end{Matrix}.
  \]
  By the use of the coordinate change $\zeta\mapsto\zeta+z\delta$, $\omega\mapsto\omega+\alpha z$, we find
  \[
    I_h^f(z,w)=\int D(\zeta,\omega)\,\vrho
    \begin{Matrix}1
      z&\zeta-z\delta\\
      \omega-\alpha z&w
    \end{Matrix}
    f
    \begin{Matrix}1
      z&\zeta\\ 
      \omega& w+\alpha\zeta+\omega\delta+\alpha z\delta
    \end{Matrix}.
  \]
  Applying \thmref{Lem}{uq-shift} with the nilpotent shift $n=\alpha\zeta+\omega\delta+\alpha z\delta$, we see that the left hand side $\smash{\int\Abs0{dw}\,\vphi(w)I_h^f(z,w)}$ equals
  \[
    \int\!D(\zeta,\omega)\!\int_{\mathrm U(q)}\Abs0{dw}\,\vrho(h^{-1}.Z)f(Z)\frac{\vphi(w-\alpha\zeta-\omega\delta-\alpha z\delta)\det(w)^q}{\det(w-\alpha\zeta-\omega\delta-\alpha z\delta)^q}.
  \]
  On appealing to Equation \eqref{eq:vrho-ber}, the claim follows. 
\end{proof}

\begin{Lem}[bos-invar]
  Let $f\in\Gamma_c(\sh O_\Omega)$ and $h=\begin{Matrix}0A&0\\0&D^{-1}\end{Matrix}$, where $A=\begin{Matrix}01&\alpha\\0&1\end{Matrix}\in_S\GL(p|q,\cplxs)$ and $D=\begin{Matrix}01&0\\\delta&1\end{Matrix}\in_S\GL(p|q,\cplxs)$. For any smooth function $\vphi$ on $\Herm^+(p)$, we have
  \[
    \int_{\Herm^+(p)}\!\frac{\Abs0{dz}}{\Abs0{\det z}^p}\,\vphi(z)I_h^f(z,w)=\!\int\!\!\frac{\Abs0{dz}}{\Abs0{\det z}^p}\!\int\!D(\zeta,\omega)\,\vrho(Z)f(Z)\vphi(z-\alpha\omega-\zeta\delta-\alpha w\delta).
  \]
  Here, on the right-hand side, the outer integral is over $\Herm^+(p)$ and the inner one is over $\cplxs^{0|p\times q}\oplus\cplxs^{0|q\times p}$.
\end{Lem}

\begin{proof}
  We may proceed as in the proof of \thmref{Lem}{ferm-invar}. Indeed,
  \[
    h.Z=AZD=
    \begin{Matrix}1
      z+\alpha\omega+\zeta\delta+\alpha w\delta&\zeta+\alpha w\\
      \omega+w\delta&w
    \end{Matrix}.
  \]
  The coordinate change $\zeta\mapsto\zeta+\alpha w$, $\omega\mapsto\omega+w\delta$ leads to
  \[
    I_h^f(z,w)=\int D(\zeta,\omega)\,\vrho
    \begin{Matrix}1
      z&\zeta-\alpha w\\
      \omega-w\delta&w
    \end{Matrix}
    f%
    \begin{Matrix}1
      z+\alpha\omega+\zeta\delta+\alpha w\delta&\zeta\\
      \omega&w
    \end{Matrix}.
  \]
  Applying \thmref{Lem}{herm-shift} with the nilpotent shift $n=\alpha\omega+\zeta\delta+\alpha w\delta$, we find that the left hand side $\int_{\Herm^+(p)}\frac{\Abs0{dz}}{\Abs0{\det(z)}^p}\vphi(z)I_h^f(z)$ equals
  \[
    \int\!D(\zeta,\omega)\!\int\!\!\frac{\Abs0{dz}}{\Abs0{\det(z)^p}}\,\vrho(h^{-1}.Z)f(Z)\frac{\vphi(z-\alpha\omega-\zeta\delta-\alpha w\delta)\det(z)^p}{\det(z-\alpha\omega-\zeta\delta-\alpha w\delta)^p}.
  \]
  As above, this proves the claim, by the use of Equation \eqref{eq:vrho-ber}.
\end{proof}

\begin{Rem}
  Observe that the Borel supergroups used in \thmref{Lem}{ferm-invar} and \thmref{Lem}{bos-invar} are opposite.
\end{Rem}

\begin{proof}[\protect{Proof of \thmref{Prop}{inv-ber}}]
  If $h\in_SH$, $h=\begin{Matrix}0A&0\\0&D^{-1}\end{Matrix}$, is of the form set out in \thmref{Lem}{ferm-invar} or \thmref{Lem}{bos-invar}, then by the token of these, $\mu$ is invariant under the action of $h$. Decomposing $A$ and $D$ in the general case into elements of block diagonal and block triangular form, we may thus assume
  \[
    A=
    \begin{Matrix}1
      a_1&0\\0&a_2
    \end{Matrix}\nd
    D=
    \begin{Matrix}1
      d_1&0\\0&d_2
    \end{Matrix}.
  \]
  Then 
  \[
    h.Z=AZD=
    \begin{Matrix}1
      a_1zd_1&a_1\zeta d_2\\
      a_2\omega d_1&a_2wd_2
    \end{Matrix}.
  \]
  The coordinate changes $\zeta\mapsto a_1\zeta d_2$ and $\omega\mapsto a_2\omega d_1$ show that $I_h^f(z,w)$ equals
  \[
    \det(a_1d_1)^q\det(a_2d_2)^p\int D(\zeta,\omega)\,\vrho
    \begin{Matrix}1
      z&a_1^{-1}\zeta d_2^{-1}\\
      a_2^{-1}\omega d_1^{-1}&w
    \end{Matrix}
    f
    \begin{Matrix}1
      a_1zd_1&\zeta\\
      \omega&a_2wd_2
    \end{Matrix},
  \]
  where as above, the integral is over $\cplxs^{0|p\times q}\oplus\cplxs^{0|q\times p}$. Since 
  \begin{align*}
    \vrho(A^{-1}ZD^{-1})&=\Ber0{A^{-1}ZD^{-1}}^{q-p}\det(a_1^{-1}zd_1^{-1})^p\det(a_2^{-1}wd_2^{-1})^q\\
    &=\det(a_1d_1)^{p-q}\det((a_1d_1)^{-1})^p\det(a_2d_2)^{q-p}\det((a_2d_2)^{-1})^q\vrho(Z)\\
    &=\det(a_1d_1)^{-q}\det(a_2d_2)^{-p}\vrho(Z),
  \end{align*}
  it follows, on applying the invariance of the densities on $\Herm^+(p)$ and $\mathrm U(q)$, that $\mu$ is invariant under the action of $h$. This proves the proposition.  
\end{proof}

\subsection{The Laplace transform of conical superfunctions}

We now come finally to the core of our paper, the explicit computation of the Laplace transforms of conical superfunctions. We will make heavy use of the facts and definitions laid down in Appendix \ref{app:super-int}.

Fix a superfunction $f\in\Gamma(\sh O_\Omega)$ and $x\in_SL(W)$. Whenever the integral converges, we define the \emph{Laplace transform} of $f$ at $x$ by 
\[
  \LT(f)(x)\defi\int_\Omega |Dy|\,e^{-\str(xy)}f(y),
\]
where we write $\Abs0{Dy}$ for the invariant Berezinian $\mu$ on $\Omega$. All integrals will be taken with respect to the standard retraction on $\Omega$.

\begin{Prop}[conv]
For $x\in_S B.1$, the integral
\[
  \LT(\Delta_\mathbf m)(x^{-1})=\int_\Omega |Dy|\,e^{-\str(x^{-1}y)}\Delta_\mathbf m(y)
\]
converges absolutely if and only if $m_j>j-1$ for $j=1,\dots,p$.
\end{Prop}

We make use of the following lemma. 

\begin{Lem}[berint-trans]
  Define 
  \[
    I_x\begin{Matrix}0z&0\\0&w\end{Matrix}\defi\int_{\cplxs^{0|p\times q}\oplus\cplxs^{0|q\times p}}D(\zeta,\omega)\,\vrho(y)e^{-\str(x^{-1}y)}\Delta_\mathbf m(y)
  \]
  for $x\in_S\GL(p|q,\cplxs)$. Let $h\in_SB_0$. Then 
  \[
    I_{h.x}\begin{Matrix}0z&0\\0&w\end{Matrix}=\chi_\mathbf m(h^{-1})I_x\Parens1{h^{-1}.\begin{Matrix}0z&0\\0&w\end{Matrix}}.
  \]
\end{Lem}

\begin{proof}
  If $h=\begin{Matrix}0A&0\\0&D\end{Matrix}\in_SK_\cplxs$, then 
  \begin{equation}\label{eq:str-equi}
    \str((h.x)^{-1}y)=\str((AxD^{-1})^{-1}y)=\str(x^{-1}A^{-1}yD)=\str(x^{-1}(h^{-1}.y)).
  \end{equation}
  Assume now that $h\in_SB_0$, where $A=\begin{Matrix}0a_1&0\\0&a_2\end{Matrix}$ and $D=\begin{Matrix}0d_1&0\\0&d_2\end{Matrix}$. Arguing as in the proof of \thmref{Prop}{inv-ber} and using Equation \eqref{eq:vrho-ber}, we obtain
  \[
    I_{h.x}\begin{Matrix}0z&0\\0&w\end{Matrix}=\frac{\det(a_1)^q\det(a_2)^p}{\det(d_1)^q\det(d_2)^p}\int\!D(\zeta,\omega)\,e^{-\str(x^{-1}y)}\Delta_{\mathbf m+q\mathds1'+(q-p)\mathds1''}(h.y')\det(w)^q
  \]
  with
  \[
    y'\defi h^{-1}.
    \begin{Matrix}1
      z&a_1^{\vphantom{-1}}\zeta d_2^{-1}\\
      a_2^{\vphantom{-1}}\omega d_1^{-1}&w
    \end{Matrix}=
    \begin{Matrix}1
      a_1^{-1}zd_2^{\vphantom{-1}}&\zeta\\
      \omega&a_2^{-1}wd_2^{\vphantom{-1}}
    \end{Matrix}.
  \]
  In view of 
  \[
    \chi_{q\mathds1'-p\mathds1''}(h^{-1})=\det(a_1^{-1}d_1^{\vphantom{-1}})^q\det(a_2^{\vphantom{-1}}d_2^{-1})^p,
  \]
  this leads to the desired conclusion immediately. 
\end{proof}

\begin{proof}[\protect{Proof of \thmref{Prop}{conv}}]
  Let us first show that the condition stated in the proposition is sufficient for the convergence of the integral.

  To that end, we perform the coordinate change $u\mapsto z=z(u)=t(u).1$ by the aid of \thmref{Prop}{herm-param}. The pullback of $\Abs0{dz}$ is $2^p\prod_{j=1}^pu_j^{2(p-j)+1}\,\Abs0{du}$, where, using the short-hand $\Abs0{du_{jk}}=\Abs0{d\Re u_{jk}}\Abs0{d\Im u_{jk}}$, we write $\Abs0{du}$ for the Lebesgue density $\prod_{j=1}^p\Abs0{du_j}\,\prod_{1\sle j<k\sle p}\Abs0{du_{jk}}$. Hence, $\Abs0{\det z}^{-p}\,\Abs0{dz}$ pulls back to $2^p\prod_{j=1}^pu_j^{-2j+1}\,\Abs0{du}$. 

  Now, let $x$ be the generic point of $B.1\cap G'.0$. Writing $x^{-1}=\begin{Matrix}0a&b\\c&d\end{Matrix}$, we see that
    \[
    (t(u)^{-1}.x)^{-1}=
    \begin{Matrix}1
      t(u)^*at(u)&t(u)^*b\\
      ct(u)&d
    \end{Matrix}.
  \]
  From \thmref{Lem}{berint-trans}, we obtain
  \begin{equation}\label{eq:zout}
    I_x\begin{Matrix}0z&0\\0&w\end{Matrix}=\chi_\mathbf m(t(u)^{-1})I_{t(u)^{-1}.x}(1,w)=\prod_{j=1}^pu_j^{2m_j}e^{-\tr(t(u)at(u)^*)}f(u,w),
  \end{equation}
  where
  \[
    f(u,w)\defi e^{\tr(dw)}\det(w)^q\int\!D(\zeta,\omega)\,\Delta_{\mathbf m+q\mathds1'+(q-p)\mathds1''}
    \begin{Matrix}0
      1&\zeta\\
      \omega&w
    \end{Matrix}\,e^{-\tr(t(u)^*b\zeta)+\tr(ct(u)\omega)}.
  \]
  Here and in the following, unless otherwise stated, it will be understood that $\cplxs^{0|p\times q}\oplus\cplxs^{0|q\times p}$ is the domain of integration for the fermionic integral $\int D(\zeta,\omega)$. Observe that $f\in\Gamma(\sh O_{\mathrm U(q)})[u]$. Here, for any locally convex vector space $E$, we denote by $E[u]$ the set of polynomials in $u_j,u_{jk},\overline{u_{jk}}$ with values in $E$. 

  Moreover, 
  \begin{align*}
    \tr\Parens1{t(u)^*at(u)}&=\tr\Parens1{at(u)t(u)^*}=\tr(az).
  \end{align*}
  To see that the integral $\int_\Omega\Abs0{Dy}\,e^{-\str(x^{-1}y)}\Delta_\mathbf m(y)$ converges absolutely, it will thus be sufficient to show that 
  \[
    J_\mathbf m(f,a)\defi\int_{\Herm^+(p)}\Abs0{du}\,\prod_{j=1}^pu_j^{2(m_j-j+1)-1}f(u)e^{-\tr(az(u))}
  \]
  converges absolutely for any $f\in\cplxs[u]$, uniformly on compact subsets with all derivatives in $a\in\Herm^+(p)+i\Herm(p)$. Since taking derivatives with respect to $a$ only introduces polynomials in $u$ into the integrand, it will be sufficient to show uniform convergence on compact subsets with respect to $a$.

  Thus, let $\alpha=(\alpha_j)_{1\sle j\sle p}\cup(\alpha_{jk},\bar\alpha_{jk})_{1\sle j<k\sle p}$, $\alpha_j,\alpha_{jk}\in\nats$, and consider 
  \[
    f(u)=u^\alpha\defi\prod_ju_j^{\alpha_j}\prod_{j<k}(u_{jk})^{\alpha_{jk}}(\overline{u_{jk}})^{\bar\alpha_{jk}}.
  \]
  Below, we will use the notation $\alpha^{(j)}=(\alpha_{jk},\bar\alpha_{jk})_{j<k\sle p}$ for any $1\sle j<p$. 

  We prove the convergence of the integral $J_\mathbf m^\alpha(a)\defi J_\mathbf m(u^\alpha,a)$ by induction on $p$. For $p=1$, we have
  \[
    J_{m_1}^{\alpha_1}(a_1)=\int_0^\infty\Abs0{du_1}\,u_1^{2m_1+\alpha_1-1}e^{-\Re(a_1)u_1^2}=\frac{\Gamma\Parens1{m_1+\frac{\alpha_1}2}}{2\Re(a_1)^{m_1+\alpha_1/2}},
  \]
  with uniform convergence on compact subsets of $\Re a_1>0$, provided that $m_1>0$.

  For $p\sge2$, by the proof of \thmref{Prop}{herm-param}, we may decompose $a$, $u$, and $z(u)$:
  \[
    a=\kbordermatrix{
        &   1       &   p-1\\
      1   &   a_1     & a^{(1)}\\
      p-1 &   (a^{(1)})^* & a'
    },\ 
    u=\kbordermatrix{
        &   1       &   p-1\\
      1   &   u_1     & u^{(1)}\\
      p-1 &   (u^{(1)})^* & u'
    },
  \]
  and 
  \[
    z(u)=\kbordermatrix{
        &   1         &   p-1\\
      1   &   u_1^2       & u_1u^{(1)}\\
      p-1 &   u_1(u^{(1)})^*  & z'(u')+(u^{(1)})^*u^{(1)}     
    }.
  \]
  Then, with $u_1$ running over $(0,\infty)$ and $u^{(1)}$ over $\cplxs^{1\times(p-1)}$, we have
  \[
    J_\mathbf m^\alpha(a)
    = J_{\mathbf m'}^{\alpha'}(a')
    \int\Abs0{du^{(1)}}\Abs0{du_1}\Abs0{u^{(1)}}^{\alpha^{(1)}}u_1^{2m_1+\alpha_1-1}e^{-a_1u_1^2-2u_1\Re a^{(1)}u^{(1)*}-\Re u^{(1)}a'u^{(1)*}}.
  \]

  We now write $\Re a'$ for the Hermitian part of the matrix $a'$, which is positive definite by assumption. Setting $\tilde a^{(1)}\defi(\Re a')^{-1/2}a^{(1)}$, we find
  \begin{align*}
    \int\Abs0{du^{(1)}}&\Abs0{u^{(1)}}^{\alpha^{(1)}}e^{-2u_1\Re a^{(1)}u^{(1)*}-\Re u^{(1)}a'u^{(1)*}}\\
    &=(\det\Re a')^{-1}\int\Abs0{du^{(1)}}\,\Abs1{(\Re a')^{-1/2}u^{(1)}}^{\alpha^{(1)}}e^{\Norm0{u^{(1)}-\tilde a^{(1)}}^2-u_1^2\Norm0{\tilde a^{(1)}}^2}\\
    &=\frac{e^{-u_1^2\Norm0{\tilde a^{(1)}}^2}}{\det\Re a'}
    \prod_{j=2}^p\int d\Abs0{u_{1j}}\,\Abs1{\Parens1{(\Re a')^{-1/2}u^{(1)}}_{1j}}^{\alpha_{1j}+\bar\alpha_{1j}}e^{-\Abs0{u_{1j}-(\tilde a^{(1)})_{1j}}^2},
  \end{align*}
  with the $u_{1j}$ running over $\cplxs$. The latter integral converges, uniformly on compact subsets as a function of $a'$ and $a^{(1)}$, since $\int_{-\infty}^\infty\Abs0{dt}\,\Abs0t^{2\beta-1} e^{-t^2}=\Gamma(\beta)$ for $\beta>0$. 

  Using the formula from the case $p=1$, we see that $J_\mathbf m^\alpha(a)$ converges uniformly on compact subsets as a function of $a$ for any $\alpha$, provided that $m_j>j-1$ for all $j=1,\dotsc,p$. This completes the proof of sufficiency. 

  Necessity follows from \cite{FK94}*{Theorem VII.1.1} by setting $x=1$ in the expression for the integral derived at the beginning of the proof. 
\end{proof}

Having established convergence, we may study the behaviour of $\LT(\Delta_\mathbf m)(x^{-1})$ as a function of $x$.

\begin{Prop}[laplace-conical]
  Assume that $m_j>j-1$ for all $j=1,\dotsc,p$. Then the function $F$ defined on $S$-valued points as $F(x)\defi\mathscr{L}(\Delta_{\bf m})(x^{-1})$ is conical.
\end{Prop}

\begin{proof}
  Let $v\in\ger b$ be an element of the Lie superalgebra of $B$. Denote by the same letter the vector field induced by $v$ on $B.1$. Following the exposition in Appendix \ref{app:vv-pt}, we may consider $v$ as an $S$-valued point of $B.1$ for $S\defi(B.1)[\eps,\tau]$; namely, after identifiying $\ger b$ with the $*[\eps,\tau]$-valued points of $B$ along $1_B:*\to B$, it is induced from $v\in\ger b$ by the action of $B$.

  Recall from Equation \eqref{eq:str-equi} that we have $\str((h^{-1}.x)^{-1}y)=\str(x^{-1}(h.y))$ for any $h\in_SB$. So we compute, with $h=v$ understood as above, that for $x\in_SB.1$
  \begin{align*}
    \sh L_v(F)(x)&=\LT(\Delta_\mathbf m)((v^{-1}.x)^{-1})=\int_\Omega\Abs0{Dy}\,e^{-\str(x^{-1}(v.y))}\Delta_\mathbf m(y)\\
    &=-\int_\Omega\Abs0{Dy}\,\sh L_v(e^{-\str(x^{-1}\cdot)})(y)\Delta_\mathbf m(y),
  \end{align*}
  where we have used Equation \eqref{eq:spec-int}.

  By \thmref{Prop}{inv-ber}, $\Abs0{Dy}$ is $H$-invariant, so $\sh L_u(\Abs0{Dy})=0$ for $u\in\ger h$. But $\ger b\subseteq\ger h$, so $\sh L_v(\Abs0{Dy})=0$, and we compute
  \begin{equation}\label{eq:lieder-ber}
    \sh L_v\Parens1{\Abs0{Dy}\,e^{-\str(x^{-1}\cdot)}\Delta_\mathbf m}
    =\Abs0{Dy}\,\Parens1{\sh L_v(e^{-\str(x^{-1}\cdot)})\Delta_\mathbf m+e^{-\str(x^{-1}\cdot)}\sh L_v(\Delta_\mathbf m)}.
  \end{equation}
  Since 
  \[
    \sh L_v(\Delta_\mathbf m)(y)=\Delta_\mathbf m(v^{-1}.y)=\chi_\mathbf m(v)\Delta_\mathbf m(y)=d\chi_\mathbf m(v)\Delta_\mathbf m(y),
  \]
  we see that the right hand side of Equation \eqref{eq:lieder-ber} is absolutely integrable, and hence, so is the left hand side. It follows that 
  \[
    \sh L_v(F)=d\chi_\mathbf m(v)F. 
  \]
  The differential equation 
  \[
    \dot\gamma=d\chi_\mathbf m(v)\gamma
  \]
  with initial condition $\gamma(0)=F(x)$, of which $\gamma(t)=F(\exp(-tv)x)$ is a solution, has values in the Fr\'echet space $\Gamma(\sh O_S)$. In general, the solutions of linear ODE with values in Fr\'echet spaces are not unique, \cf Ref.~\cite{neeb-monastir}. However, $d\chi_\mathbf m(v)$ is a scalar, so it induces a continuous endomorphism in the topology on $\Gamma(\sh O_S)$ generated by any fixed continuous norm from a defining set, and we may apply the uniqueness theorem from the Banach case. 

  Using the facts that $\exp:\ger b\to B$ is a local isomorphism and $B$ is connected, we deduce as in the Lie group case that $F(b^{-1}x)=\chi_\mathbf m(b)F(x)$ for any $b\in_SB$. 
\end{proof}

\begin{Def}[gamma][gamma function]
The gamma function of $\Omega$ is defined as:
\[
  \Gamma_\Omega(\mathbf m)\defi\LT(\Delta_\mathbf m)(1)=\int_\Omega |Dy| e^{-\str(y)}\Delta_\mathbf m(y),
\]
whenever $m_j>j-1$ \fa $j=1,\dotsc,p$.
\end{Def}

With this notation, the following is immediate.

\begin{Cor}[delta-lap]
  Assume $m_j>j-1$ \fa $j=1,\dotsc,p$. Then \fa $x\in_SB.1$, we have
  \[
    \LT(\Delta_\mathbf m)(x^{-1})=\Gamma_\Omega(\mathbf m)\Delta_\mathbf m(x).
  \]
\end{Cor}

Of course, the value of this corollary depends on the extent to which we have control over $\Gamma_\Omega$. In fact, we can give an entirely explicit expression, as follows.

\begin{Th}[gamma]
  Let $m_j>j-1$ \fa $j=1,\dotsc,p$. We have
  \begin{align*}
    \Gamma_\Omega({\mathbf m}) =&\; (2\pi)^{p(p-1)/2}\prod_{j=1}^p \Gamma(m_j-(j-1))\\
    & \times 
    \prod_{k=1}^q \frac{\Gamma(q-(k-1))}{\Gamma(m_{p+k}+q-(k-1))}\frac{\Gamma(m_{p+k}+k)}{\Gamma(m_{p+k}-p+k)}.
  \end{align*}
  In particular, $\Gamma_\Omega(\mathbf m)$ extends uniquely as a meromorphic function of $\mathbf m\in\cplxs^{p+q}$, and it has neither zeros nor poles if 
  \begin{eqnarray*}
    m_j>j-1& &j=1,\dotsc,p,\\
    m_{p+k}>p-k& &k=1,\dotsc,q.
  \end{eqnarray*}
  If, under this assumption, $\mathbf m$ is a double partition, then it is a hook partition. 
\end{Th}

\begin{proof}
  Our stategy of proof is to reduce the computation to three separate computations, which take place on the fermionic part $\cplxs^{0|p\times q}\oplus\cplxs^{0|q\times p}$, the `fermion-fermion sector' $\mathrm U(q)$, and on the `boson-boson sector' $\Herm^+(p)$, respectively. 

  Let us decompose $\mathbf m=(\mathbf m',\mathbf m'')$ where 
  \[
    \mathbf m'=(m_1,\dotsc,m_p)\nd\mathbf m''=(m_{p+1},\dotsc,m_{p+q}).
  \] 
  Then from Equation \eqref{eq:zout}, we have 
  \[
    \int D(\zeta,\omega)\,e^{-\str(y)}\Delta_{\mathbf m+q\mathds1'+(q-p)\mathds1'}(y)\det(w)^q=e^{-\tr(z)+\tr(w)}\Delta_{\mathbf m'}(z)\psi(w),
  \]
  where
  \begin{equation}\label{eq:psidef}
    \psi(w)\defi\int\!D(\zeta,\omega)\,\Delta_{\mathbf m+q\mathds1'+(q-p)\mathds1''}\begin{Matrix}11&\zeta\\\omega&w\end{Matrix}\det(w)^q.
  \end{equation}
  Appealing to \cite{FK94}*{Theorem VII.1.1}, we see that 
  \[
    \Gamma_\Omega(\mathbf m)=(2\pi)^{p(p-1)/2}\prod_{j=1}^p \Gamma(m_j-(j-1))\int_{\mathrm U(q)}\Abs0{dw}\,e^{\tr(w)}\psi(w).
  \]
  The final statement now follows from \thmref{Lem}{gamma-berezin} and \thmref{Lem}{gamma-unitary} below.
\end{proof}

\begin{Lem}[gamma-berezin]
  In the notation from Equation \eqref{eq:psidef}, we have 
  \[
    \psi(w)=(\Delta_{\mathbf m''}(w))^{-1}\prod_{k=1}^q\frac{\Gamma(m_{p+k}+k)}{\Gamma(m_{p+k}-p+k)},
  \]
  where $\mathbf m'':=(m_{p+1},\dotsc,m_{p+q})$.
\end{Lem}

\begin{proof}
  Notice that $\psi(w)$ is well-defined for $w$ in the open subset $B_q\overline B_q\subseteq\cplxs^{q\times q}$, in view of \thmref{Lem}{bigcell}. On applying \thmref{Lem}{berint-trans} for $x=0$, we see that 
  \[
    \psi(h^{-1}.w)=\chi_\mathbf m(h)\psi(w)
  \]
  for any $h=\begin{Matrix}0A&0\\0&D\end{Matrix}\in\mathrm U(q)\subseteq B_0$, $A=\begin{Matrix}01&0\\0&a\end{Matrix}$, $D=\begin{Matrix}01&0\\0&d\end{Matrix}$. Applying \thmref{Lem}{conical-exhaustion} for the case of $p=0$, we find
  \[
    \psi(w)=\Delta_\mathbf m\begin{Matrix}01&0\\0&w\end{Matrix}\psi(1)=(\Delta_{\mathbf m''}(w))^{-1}\psi(1).
  \]

  We have 
  \[
    \begin{Matrix}1
      1&\zeta\\\omega&1
    \end{Matrix}=
    \begin{Matrix}1
      1&0\\
      \omega&1
    \end{Matrix}
    \begin{Matrix}1
      1&0\\0&1-\omega\zeta
    \end{Matrix}
    \begin{Matrix}1
      1&\zeta\\0&1
    \end{Matrix},
  \]
  so 
  \[
    \Delta_{\mathbf m+q\mathds1'+(q-p)\mathds1''}
    \begin{Matrix}1
      1&\zeta\\\omega&1
    \end{Matrix}
    =(\Delta_{\mathbf m''+(q-p)\mathds1''}(1-\omega\zeta))^{-1}
  \]

  To calculate the resulting Berezin integral $\psi(1)=\gamma_{-(\mathbf m''+(q-p)\mathds1'')}$, where
  \[
    \gamma_{-\mathbf m''}\defi\int\!D(\zeta,\omega)\,(\Delta_{\mathbf m''}(1-\omega\zeta))^{-1},
  \]
  we decompose the matrices as
  \[
    \omega = \kbordermatrix{
      & p  \\
      1 & \omega_1\\
      q-1 & \omega^\prime} 
    \nd 
    \zeta=\kbordermatrix{
      & 1 & q-1 \\
      p & \zeta_1 & \zeta^\prime}.
  \]
  So, $1-\omega\zeta$ can be decomposed as
  \[
    \begin{Matrix}1
      1 & 0 \\
      \omega'\zeta_1 (1 - \omega_1\zeta_1)^{-1} & 1
    \end{Matrix}
    \begin{Matrix}1
      1 - \omega_1\zeta_1 & 0 \\
      0 & 1 - \omega'(1+N)\zeta'  
    \end{Matrix}
    \begin{Matrix}1
      1 & (1 - \omega_1\zeta_1)^{-1} \omega_1\zeta' \\
      0 & 1     
    \end{Matrix},
  \]
  where $N= \zeta_1 (1 - \omega_1\zeta_1)^{-1} \omega_1$. Then $\gamma_{-(\mathbf m''+(q-p)\mathds1'')}$ becomes
  \[
    \int\!\!D(\zeta_1,\omega_1)\!\int\!D(\zeta',\omega')\,\Delta_{\mathbf m''+(q-p)\mathds1''}\!\!%
    \begin{Matrix}1
      1-\omega_1\zeta_1&0\\
      0&1-\omega'(1+N)\zeta'
    \end{Matrix}^{-1}.
  \]
  Of these Berezin integrals, the outer is over $\cplxs^{0|p\times 1}\oplus\cplxs^{0|1\times p}$, while the inner is over the space $\cplxs^{0|p\times(q-1)}\oplus\cplxs^{0|(q-1)\times p}$.

  Observe that $\tr N=-(1-\omega_1\zeta_1)^{-1}\omega_1\zeta_1$ and hence 
  \[
    N^2=\zeta_1(1-\omega_1\zeta_1)^{-1}\omega_1\zeta_1(1-\omega_1\zeta_1)^{-1}\omega_1=-\tr N\cdot N.
  \]
  This implies 
  \[
    \tr(N^{k+1})=-\tr N\tr(N^k)=\dotsm=-(-\tr N)^{k+1}, 
  \]
  so that 
  \[
    \tr\log(1+N)=\sum_{k=1}^\infty\frac{(-1)^{k-1}}k\tr\Parens1{N^k}=-\log(1-\tr N).
  \]
  Then
  \begin{align*}
    \det(1+N)^{-1}&=\exp(-\tr\log(1+N))=1-\tr N\\
    &=1+(1-\omega_1\zeta_1)^{-1}\omega_1\zeta_1=(1-\omega_1\zeta_1)^{-1},
  \end{align*}
  and making a change of odd coordinates $\zeta^\prime \mapsto (1+N)\zeta^\prime$, the Berezin integral over $D(\zeta',\omega')$ simplifies to
  \[
    (1-\omega_1\zeta_1)^{q-1}\int\!D(\zeta',\omega')\,\Delta_{\mathbf m''+(q-p)\mathds1''}
    \begin{Matrix}1
      1-\omega_1\zeta_1&0\\
      0&1-\omega'\zeta'
    \end{Matrix}^{-1}.
  \]
  We obtain 
  \[
    \gamma_{-(\mathbf m''+(q-p)\mathds1'')}=\gamma_{-m_{p+1}-1+p}\gamma_{-(m_{p+2}+q-p,\dotsc,m_{p+q}+q-p)}=\dotsm=\prod_{j=1}^q\gamma_{-m_{p+j}-j+p}.    
  \]
  Finally, writing $a\defi1-\sum_{j=2}^p\omega_{1j}\zeta_{j1}$, we compute
  \[
    (1-\omega_1\zeta_1)^{-m}=a^{-m}(1-a^{-1}\omega_{11}\zeta_{11})^{-m}=a^{-m-1}(1+m\omega_{11}\zeta_{11})
  \]
  which recursively gives 
  \[
    \gamma_{-m}=m(m+1)\dotsm(m+p-1)=\frac{\Gamma(m+p)}{\Gamma(m)},
  \]
  and hence, our claim. 
\end{proof}

\begin{Lem}[gamma-unitary]
  We have
  \[
    \int_{{\rm U}(q)} |Dw| e^{{\rm tr}(w)}(\Delta_{\mathbf m''}(w))^{-1}=\prod_{k=1}^q \frac{\Gamma(q-(k-1))}{\Gamma(m_{p+k}+q-(k-1))},
  \]
  where $\mathbf m''\defi(m_{p+1},\dotsc,m_{p+q})$.
\end{Lem}

\begin{proof}
  We use spherical polynomials for this computation. They are defined as
  \[
    \Phi_\mathbf r(x):= \int_{{\rm U}(q)}\Abs0{Dw}\,\Delta_\mathbf r(wx),
  \]
  where $\mathbf r$ is an arbitrary multi-index of length $q$.

  Thanks to \cite{FK94}*{Proposition XII.1.3.(i)} we can write the exponential in this integral as an absolutely convergent series of spherical functions:
  \[
    e^{\tr(w)} = \sum_{\mathbf n\sge 0} \frac{d_\mathbf n}{q_\mathbf n}\Phi_\mathbf n(w),
  \]
  where $d_\mathbf n$ is the dimension of the finite-dimensional irreducible $\mathrm U(q)$-module of highest weight $\mathbf n$, and 
  \[
    q_\mathbf n\defi\prod_{k=1}^q\frac{\Gamma(n_k+q-(k-1))}{\Gamma(q-(k-1))}.
  \]
  Therefore, we can write our integral as
  \[
    \sum_{\mathbf n\sge 0} \frac{d_\mathbf n}{q_{\mathbf n}}\int_{\mathrm U(q)} |Dw|\,\Phi_\mathbf n(w) (\Delta_{\mathbf m''}(w))^{-1}.
  \]
  Notice that since $w\in\mathrm U(q)$, we have $(\Delta_{\mathbf m''}(w))^{-1} = \Delta_{\mathbf m''}(w^{-1})$. Due to the $\mathrm U(q)$-in\-va\-riance of $\Phi_\mathbf n$, we have
  \[
    \int_{\mathrm U(q)}|Du|\,\Phi_\mathbf n(xu) = \Phi_\mathbf n(x),
  \]
  so we obtain
  \begin{align*}
    \int_{\mathrm U(q)}|Dw|\,\Phi_{\bf n}(w) \Delta_{\mathbf m''}(w^{-1})
    &=\int_{\mathrm U(q)\times\mathrm U(q)} |Dw||Du|\,\Phi_\mathbf n(wu) \Delta_{\mathbf m''}(w^{-1}) \\
    &=\int_{\mathrm U(q)\times\mathrm U(q)} |Dw||Du|\,\Phi_\mathbf n(w) \Delta_{\mathbf m''}(uw^{-1}) \\
    &=\int_{\mathrm U(q)}|Dw|\,\Phi_\mathbf n(w) \Phi_{\mathbf m''}(w^{-1}).
  \end{align*}
  Further, it is easy to see that $\Phi_{\bf m^{\prime\prime}}(w^{-1})=\Phi_{\mathbf m''}(w^*) = \overline{\Phi_{\bf m^{\prime\prime}}(w)}$. 

  By the classical Schur orthogonality relations, this integral is only non-zero when $\mathbf m''=\mathbf n$, in which case the answer is $1/d_{\mathbf m''}$. Therefore,
  \[
    \int_{{\rm U}(q)} |Dw|\,e^{\tr(w)}(\Delta_{\mathbf m''}(w))^{-1} = 1/q_{\mathbf m''},
  \]
  which gives the claim.
\end{proof}

\appendix

\section{The functor of points}\label{app:super}

\subsection{The language of $S$-valued points}\label{app:functorpoints}

For a manifold $X$, a point can be thought of as a morphism $* \to X$, and this completely determines $X$. However, for a supermanifold $X$, such a morphism $* \to X$ is again only a point of the underlying manifold $X_0$, and therefore does not capture the supergeometric features of $X$. To deal with this, the notion of points has to be extended.

This idea is familiar in algebraic geometry. Here, it is common to talk about the $K$-rational points of a scheme $X$, which are nothing but morphisms $\mathrm{Spec}(K)\to X$. Grothendieck extended this idea and considered the scheme along with its $A$-points, for all commutative rings $A$. Then $X$ is completely recaptured by its collection of $A$-points, for all commutative rings $A$, along with admissible morphisms.

More generally, if $\cat C$ is any category, and $X$ is an object of $\cat C$, then an \Define{$S$-valued point} (where $S$ is another object of $\cat C$) is defined to be a morphism $x:S\to X$. One may view this as a `deformed' or `parametrised' point. Suggestively, one writes $x\in_SX$ in this case, and denotes the set of all $x\in_SX$ by $X(S)$.

For any morphism $f:X\to Y$, one may define a set-map $f_S:X(S)\to Y(S)$ by 
\[
  f_S(x)\defi f(x)\defi f\circ x\in_SY\mathfa x\in_SX.
\]
Clearly, the values $f(x)$ completely determine $f$, as can be seen by evaluating at the \Define{generic point} $x={\id}_X\in_XX$.

In fact, more is true. The following statement is known as Yoneda's Lemma \cite{maclane}: Given a collection of set-maps $f_S:X(S)\to Y(S)$, there exists a morphism $f:X\to Y$ \scth $f_S(x)=f(x)$ \fa $x\in_SX$ if and only if 
\[
  f_T(x(t))=f_S(x)(t)\mathfa t:T\to S.
\]
The points $x(t)$ are called \Define{specialisations of $x$}, so the condition states that the collection $(f_S)$ is invariant under specialisation.

The above facts are usually stated in the following more abstract form: For any object $X$, we have a set-valued functor $X(-):\cat C^{op}\to\Sets$, and the set of natural transformations $X(-)\to Y(-)$ is naturally bijective to the set of morphisms $X\to Y$. Thus, the functor $X\mapsto X(-)$ from $\cat C$ to $[\cat C^{op},\Sets]$, called the \Define{Yoneda embedding}, is fully faithful.

The Yoneda embedding preserves products \cite{maclane}, so if $\cat C$ admits finite products, it induces a fully faithful embedding of the category of group objects in $\cat C$ into the category $[\cat C^{op},\cat{Grp}]$ of group-valued functors. 

In other words, we have the following: Let $X$ be an object in $\cat C$. Then $X$ is a group object if and only if for any $S$, $X(S)$ admits a group law, which is invariant under specialisation. We will assiduously apply this point of view to the categories of complex supermanifolds and of \emph{cs} manifolds. 

\subsection{Vector fields and generalised points}\label{app:vv-pt}

We now show how vector fields can be understood in terms of generalised points. Among other things, this is a framework enabling us to exchange super-integration and differentiation under suitable assumptions.

Let $X$ and $S$ be \emph{cs} manifolds. Then $X$ is called a \Define{\emph{cs} manifold over $S$}, written $X/S$, if supplied with some morphism $X\to S$, which on some open cover $U_\alpha$ of $X$ fits into a commutative diagram 
\begin{center}
  \begin{tikzcd}
    U_\alpha\rar{}\dar{} & S\times Y\dar{p_1}\\
    V_\alpha\rar{}&S
  \end{tikzcd}
\end{center}
where the rows are open embeddings. Usually, we will consider only products, but the general language will be efficient nonetheless. There is an obvious notion of morphisms over $S$, which we denote $X/S\to Y/S$.

A system of (local) \Define{fibre coordinates} is given by the system $(x^a)=(x,\xi)$ of superfunctions on some trivialising open subspace $U\subseteq X$ obtained by pullback along a trivialisation from a coordinate system in the fibre $Y$.

If $X/S$ is a \emph{cs} manifold over $S$, then the \Define{relative tangent sheaf} is defined by 
\[
  \sh T_{X/S}\defi\ShGDer[_{p_0^{-1}\sh O_S}]0{\sh O_X,\sh O_X},
\]
the sheaf of superderivations of $\sh O_X$ which are linear over $\sh O_S$. Here, $p$ denotes the morphism $X\to S$. It is a basic fact that $\sh T_{X/S}$ is a locally free $\sh O_X$-module, with rank equal to the fibre dimension of $X/S$. 

More generally, let $\vphi:X/S\to Y/S$ be a morphisms over $S$. We let 
\[
  \sh T_{X/S\to Y/S}\defi\ShGDer[_{p_0^{-1}\sh O_S}]0{\vphi_0^{-1}\sh O_Y,\sh O_X}
\]
and call this the \Define{tangent sheaf along $\vphi$ over $S$}. Written out explicitly, the derivation property of a homogeneous element $\delta\in\sh T_{X/S\to Y/S}(U)$ is 
  \[
    \delta(fg)=\delta(f)\vphi^\sharp(g)+(-1)^{\Abs0f\Abs0\delta}\vphi^\sharp(f)\delta(g)
  \]
  \fa homogeneous $f,g\in\sh O_Y(V)$, $V\subseteq Y_0$ an open neighbourhood of $\vphi_0(U)$. The usual relative tangent bundle corresponds to $\vphi={\id}_X$.

We denote by $S[\eps,\tau]$ the $\cplxs$-superspace $(S_0,\sh O_S[\eps,\tau])$, where $\sh O_S[\eps,\tau]\defi\sh O_S\otimes_\cplxs A$, $A$ denoting the superalgebra $\cplxs[\eps,\tau]/(\eps^2,\eps\tau)$, where $\eps$ and $\tau$ are understood to be even and odd indeterminates, respectively. If $X/S$ and $Y/S$ are \emph{cs} manifolds over $S$ and $\vphi:X/S\to Y/S$ is a morphism over $S$, then there is a natural bijection 
  \[
    \Set1{\gamma\in\Hom[_S]1{X[\eps,\tau],Y}}{\gamma|_{\eps=\tau=0}=\vphi}\to\Gamma(\sh T_{X/S\to Y/S}):\gamma\mapsto\delta
  \]
  given by the equation $\gamma^\sharp(f)=\vphi^\sharp(f)+\eps\delta_\ev(f)+\tau\delta_\odd(f)$ for all local sections $f$ of $\sh O_Y$. 

In particular, consider the case of a \emph{cs} Lie supergroup $G$. The Lie superalgebra $\ger g$ is by definition the fibre over $1$ of the tangent bundle. Equivalently, elements of $\ger g$ may be seen as vector fields along the morphism $1_G:1\to G$, that is, as $*[\tau,\eps]$-valued points of $G$ along $1_G$.

\subsection{Nilpotent shifts of cycles in middle dimension}\label{app:points}

We will now show how the technique of nilpotent shifts common in physics can be understood in terms of $S$-valued points. 

Let $W\cong\cplxs^q$ be a complex vector space. Its associated $q$-dimensional complex manifold, with its sheaf $\sh O_W$ of complex-analytic functions, is naturally a complex supermanifold (namely, $L(W)$).

Abusing notation, we write $N$ for the \emph{cs} manifold associated with $\cplxs^{0|N}$. (This is also a complex supermanifold.) Assume 
\[
  n\in\Gamma(\sh N_{N,\ev})^q=\Gamma(\sh N_{N,\ev}\otimes W),
\]
where $\sh N_N=\bigwedge^+(\cplxs^N)^*$ is the ideal of $\sh O_N=\bigwedge(\cplxs^N)^*$ generated by $\sh O_{N,\odd}$. 

The generic point $w={\id}_W\in_WW$ of $W$ corresponds by Leites's Theorem \cite{Lei80} to the element $w=\sum_ie^i\otimes e_i\in\Gamma(\sh O_W\otimes W)$ where $e_i$, $i=1,\dotsc,2q$, is a real basis of $W$, and $e^i$ is its dual basis. The sum 
\[
  w+n\in\Gamma((\sh O_{W,\ev}\otimes1\oplus1\otimes\sh O_{N,\ev})\otimes W)\subseteq\Gamma(\sh O_{W\times N,\ev}\otimes W)
\]
corresponds by Leites's Theorem to a unique morphism $\phi:W\times N\to W$ of complex supermanifolds.

In particular, this gives a definite meaning to $f(w+n)=\phi^\sharp(f)$ for any complex analytic function $f$ defined on an open subset of $W$. It is known that 
\[
  f(w+n)=f(w)+\sum_{k=1}^N\frac1{k!}d^kf(w)(n,\dotsc,n)
\]
where the derivatives are extended multi-linearly over $\sh N_{N,\ev}$.

Let now $X$ be a closed real submanifold of $W$ of dimension $q$. We call such an $X$ a \emph{mid-dimensional cycle}. In this case, we call $n$ a \emph{nilpotent shift} for $X$.

By the use of the embedding $j:X\to W$, the real tangent space at any point of $X$ is naturally identified with a $q$-dimensional real subspace of $W$, and this gives a real vector bundle map $Tj:TX\to X\times W$. Thus, the complex tangent space at any point is naturally identified with $W$, which gives an isomorphism of complex vector bundles $T^\cplxs j:T^\cplxs X\to X\times W$.

For any superfunction $f\in\sh O_X(U_0)$, where $U\subseteq X$ is open, and for any $y\in_SU$, we define 
\[
  f(y+n)\defi f(y)+\sum_{k=1}^N\frac1{k!}T^kf(y)((T^\cplxs j)^{-1}(y,n),\dotsc,(T^\cplxs j)^{-1}(y,n))
\]
by multi-linear extension of the higher order tangent maps. 

Here, the left hand side lies in $\Gamma(\sh O_{S\times N})$. In particular, this defines a unique morphism $X\times N\to X$, which sends $f$ to $f(x+n)$, $x={\id}_X\in_XX$ denoting the generic point of $X$.

\section{Integration on supermanifolds}\label{app:super-int}

We will need to consider super-integrals depending on some parameters. An appropriate framework for this is that of relative Berezinians, paired with the understanding of parameter dependence in terms of $S$-valued points. 

\subsection{Relative Berezinians and fibre integrals}\label{sec:ber}

Let $X$ be a \emph{cs} manifold over $S$, and $\smash{\Omega^1_{X/S}}$ be the module of relative $1$-forms, by definition dual to $\sh T_{X/S}$. Then we define the sheaf of \Define{relative Berezinians} $\sh Ber_{X/S}$ to be the Berezinian sheaf associated to the locally free $\sh O_X$-module $\smash{\Pi\Omega^1_{X/S}}$ obtained by parity reversal. Furthermore, the sheaf of \Define{relative Berezinian densities} $\Abs0{\sh Ber}_{X/S}$ is the twist by the relative orientation sheaf, \ie 
\[
  \Abs0{\sh Ber}_{X/S}\defi\sh Ber_{X/S}\otimes_\ints or_{X_0/S_0}.
\]

Given a system of local fibre coordinates $(x^a)=(x,\xi)$ on $U$, their coordinate derivations $\frac\partial{\partial x^a}$ form an $\sh O_X|_{U_0}$-module basis of $\smash{\sh T_{X/S}|_{U_0}}$, with dual basis $dx^a$ of $\smash{\Omega^1_{X/S}|_{U_0}}$. One may thus consider the distinguished basis 
\[
  \Abs0{D(x^a)}=\Abs0{D(x,\xi)}=dx_1\dots dx_p\frac{\partial^\Pi}{\partial\xi^1}\dots\frac{\partial^\Pi}{\partial\xi^q}
\]
of the module of Berezinian densities $\Abs0{\sh Ber}_{X/S}$, \cf Ref.~\cite{Man88}.

If $X/S$ is a direct product $X=S\times Y$, then 
\[
  \Abs0{\sh Ber}_{X/S}=p_2^*\Parens1{\Abs0{\sh Ber}_Y}=\sh O_X\otimes_{p_{2,0}^{-1}\sh O_Y}p_{2,0}^{-1}\Abs0{\sh Ber}_Y.
\]
In particular, the integral over $Y$ of compactly supported Berezinian densities defines the integral over $X$ of a section of $(p_0)_!\Abs0{\sh Ber}_{X/S}$, where $(-)_!$ denotes the functor of direct image with compact supports \cite{iversen}. We denote the quantity thus obtained by 
\[
  \fibint[S]X\omega\in\Gamma(\sh O_S)\mathfa\omega\in\Gamma\Parens1{(p_0)_!\Abs0{\sh Ber}_{X/S}},
\]
and call this the \Define{fibre integral} of $\omega$.

We will, however, have to consider fibre integrals in a more general setting, beyond compact supports. Henceforth, we assume for simplicity that $X=S\times Y$. A \Define{fibre retraction} for $X$ is a morphism $r:Y\to Y_0$ which is left inverse to the canonical embedding $j:Y_0\to Y$, where $Y_0$ is the underlying manifold of $Y$.

A system of fibre coordinates $(x,\xi)$ of $X/S$ is called \Define{adapted} to $r$ if $x=r^\sharp(x_0)$. Given an adapted system of fibre coordinates, we may write $\omega=\Abs0{D(x,\xi)}\,f$ and 
\[
  f=\sum_{I\subseteq\{1,\dotsc,q\}}({\id}\times r)^\sharp(f_I)\,\xi^I
\]
for unique coefficients $f_I\in\Gamma(\sh O_{S\times Y_0})$, where $\dim Y=*|q$. Then one defines
\[
  \fibint[S\times Y_0]X\omega\defi\Abs0{dx_0}\,f_{\{1,\dotsc,q\}}\in\Gamma(\Abs0{\sh Ber}_{(S\times Y_0)/S}).
\]
Note that $\Abs0{\sh Ber}_{(S\times Y_0)/S}$ is $p_2^*$ of the sheaf of ordinary densities on the manifold $Y_0$, so we may write $\Abs0{dx_0}$.

This fibre integral only depends on $r$, and not on the choice of an adapted system of fibre coordinates. (See Ref.~\cite{ahp} for the absolute case.) If the resulting relative density is absolutely integrable along the fibre $Y_0$, then we say that $\omega$ is \Define{absolutely integrable} with respect to $r$, and define
\[
  \fibint[S]X\omega\defi\fibint[S]{S\times Y_0}\Bracks3{\fibint[S\times Y_0]X\omega}\in\Gamma(\sh O_S).
\]
Both this quantity and its existence depend heavily on $r$. 

Using the topology on $\Gamma(\sh O_S)$ introduced below, in Appendix \ref{app:super-lap}, the convergence may be understood in terms of vector-valued integrals. The relative density $\tfibint[S\times Y_0]X\omega$ may be viewed as a $\Gamma(\sh O_S)$-valued density on $Y_0$. It is absolutely integrable along $Y_0$ if and only if the corresponding vector-valued density is Bochner integrable. 

We shall use the language of $S$-points discussed above in Appendix \ref{app:functorpoints} to manipulate integrals of relative Berezinian densities in a hopefully more comprehensible formalism. This also gives a rigorous foundation for the super-integral notation common in the physics literature. 

If $f$ is a superfunction on $X=S\times Y$ and we are given some relative Berezinian density $\Abs0{Dy}$ on $X/S$, then we write
\[
  \int_Y\Abs0{Dy}\,f(s,y)\defi\fibint[S]X\Abs0{Dy}\,f.
\]
This is justified by the convention that the generic points of $S$ and $Y$ are denoted by $s$ and $y$, respectively. Moreover, it is easy to see that this notation behaves well under specialisation, since
\begin{equation}\label{eq:spec-int}
  \int_Y\Abs0{Dy}\,f(s(t),y)=\fibint[T]X(t\times{\id})^\sharp(\Abs0{Dy}\,f)=t^\sharp\Bracks3{\fibint[S]X\Abs0{Dy}\,f}
\end{equation}
for any $t\in_TS$. This follows from the fact that the fibre retractions are respected by the morphism $t\times{\id}$.

\subsection{Berezin integrals and nilpotent shifts}

We now return to the nilpotent shifts previously considered in Appendix \ref{app:points}, and apply them to certain integrals. 

Let $X$ be a mid-dimensional cycle in the complex vector space $W\cong\cplxs^q$. Let $X$ carry some pseudo-Riemannian metric $g$ and $\mu_g$ be the induced Riemannian density.

The following is a straightforward generalisation of \cite{LSZ08}*{Lemma 4.13}.

\begin{Lem}[shift]
  Let $n$ be a nilpotent shift for $X$. Then for any compactly supported smooth function $f$ on $X$, we have
  \[
    \int_Xf(x+n)\,d\mu_g(x)=\int_Xf(x)J_t(x)\,d\mu_g(x),
  \]
  where $J_t\in\Gamma(\sh O_{X\times N})$ is the solution of the ODE 
  \[
    \frac d{dt}\log J_t(y,s)=-\ddiv v_n(y-tn,s)
  \]
  with inital condition $J_0=1$ and $v_n$ is the vector field on $X\times N$ over $N$, defined by
  \[
    v_n(k)(y,s)\defi\frac d{dt}k(y+tn,s)\Big|_{t=0}
  \]
  for any smooth $k$ on some open subspace $U\subseteq X\times N$ and any $(y,s)\in_S(X\times N)$.
\end{Lem}

\begin{proof}
  Consider for fixed $t\in\reals$ the morphism $\phi_t:X\times N\to X\times N$, defined by 
  \[
    \phi_t(y,z)\defi (y-tn,z)
  \]
  for any $(y,z)\in_SX\times N$. 

  Then $\phi_t$ is an isomorphism over $N$, with inverse $\phi_{-t}$. We may consider $\mu_g$ as a Berezinian density of $X\times N$ over $N$. Thus, we have
  \[
    \int_Xf(x+tn)\,d\mu_g(x)=\fibint[N]{X\times N}\phi_{-t}^\sharp(p_1^\sharp(f))\,\mu_g=\fibint[N]{X\times N}p_1^\sharp(f)\,\phi_t^\sharp(\mu_g).
  \]

  Since $\sh Ber_{X\times N/N}$ is a free $\sh O_{X\times N}$-module with module basis $\mu_g$, there exists a unique even $J_t\in\Gamma(\sh O_{X\times N})$ \scth 
  \[
    \phi_t^\sharp(\mu_g)=J_t\,\mu_g.
  \]
  These superfunctions depend smoothly on $t$. Indeed, one may consider the morphism $\phi:\reals\times X\times N\to\reals\times X\times N$, given by 
  \[
    \phi(t,y,z)\defi(t,y-tn,z)
  \]
  \fa $(t,y,z)\in_S\reals\times X\times N$. Then $\phi$ is an isomorphism over $\reals\times N$, and $\phi(t,y,z)=(t,\phi_t(y,z))$ if $t$ is the specialisation of an ordinary point of $\reals$.

  Now, $\phi_t$ is the flow of the vector field $-v_n$, considered as vector field on $(X\times N)/N$. Indeed, if $h$ is a smooth function on $\reals\times X\times N$, then
  \[
    \frac d{dt}\phi^\sharp(h)(t,y,z)=\frac d{ds}h(t,y-(s+t)n,z)\Big|_{s=0}=-\phi^\sharp\Parens1{({\id}\otimes v_n)h}(t,y,z).
  \]
  We have 
  \[
    \frac d{dt}J_t=\frac d{d\tau}\phi_t^\sharp(J_\tau)\Big|_{\tau=0}=-\phi^\sharp_t(\ddiv_g v_n)\cdot J_t,
  \]
  so
  \[
    \frac d{dt}\log J_t(y,s)=-\phi_t^\sharp(\ddiv v_n)(y,s)=-\ddiv v_n(y-tn,s)
  \]
  with $J_0=1$, proving the lemma.
\end{proof}

Let $G$ be a Lie group acting linearly on $W$ and $o\in W$ \scth $K\defi G_o$ is a compact subgroup, open in the fixed point set of an involutive automorphism $\theta$ of $G$. Then the orbit $X\defi G.o=G/K$ is a Riemannian symmetric space. 

Also denote by $\theta$ the involutive automorphism of the Lie algebra $\ger g$ of $G$ induced by $\theta$. Then $\ger g=\ger k\oplus\ger p$ where $\ger k$, the $+1$ eigenspace of $\theta$, is the Lie algebra of $K$, and $\ger p$, the $-1$ eigenspace of $\theta$, identifies as a $K$-module with $T_oX=\ger g/\ger k$.

We will denote the action of $G$ on $W$ by $g.w$, and the derived action of $\ger g$ by $u.w$. Since $\ger p_\cplxs=W$ if $X$ is mid-dimensional in $W\cong\cplxs^q$ as above, we obtain for any $w\in W$ endomorphisms $R_w$ of $\ger p_\cplxs$ by $R_w(u)\defi u.w$.

In this setting, we have the following generalisation of \cite{LSZ08}*{Lemma 4.12}.

\begin{Lem}[shift-div]
  Assume that $X=G/K$ is a mid-dimensional cycle. Then 
  \[
    (\ddiv v_n)(gK)=-\tfrac12\tr_{\ger p_\cplxs}(R_{g^{-1}.n})
  \]
  for any nilpotent shift $n$.
\end{Lem}

\begin{proof}
  Let $\nabla$ denote the Levi-Civita connection of $X$. For any vector field $v$ on $X$, one has the identity
  \[
    \ddiv v=\tr\nabla v
  \]
  where $\nabla v$ denotes the endomorphism of the tangent sheaf $\sh T_X$ given by $u\mapsto\nabla_uv$. This statement extends immediately to the case of relative vector fields. In order to compute this explicitly, we use the transitive $G$-action to reduce everything to a Lie algebra computation, as follows. 

  Consider the vector bundle $G\times^K\ger p$ associated to the principal $K$-bundle $G\to X$ \via the adjoint action of $K$ on $\ger p$. By definition, this is  the quotient $(G\times\ger p)/K$, where the right action of $K$ on $G\times\ger p$ is given by $(g,u)\cdot k\defi(gk,\Ad(k^{-1})(u))$. We denote the equivalence class of $(g,u)$ in this quotient by $[g,u]$.

  There is a natural vector bundle map $G\times^K\ger p\to TX$ mapping $[g,u]$ to the the tangent vector $\dot\beta_u(0)$ at $gK$, given as the derivative of the geodesic $\beta_u(t)\defi g\exp(tu)$ \cite{paradan}*{Corollary 5.8}. This map is a $G$-equivariant isomorphism. Here, on $TX$, $g\in G$ acts by the derivative $TL_g$ of the left multiplication $L_g$ of $X$, while it acts on $G\times^K\ger p$ by sending $[g',u]$ to $[gg',u]$. Similarly, we identify $\End0{TX}$ with $G\times^K\End0{\ger p}$.

  Let $u\in\ger g$. The value at $gK$ of the fundamental vector field $u_X$ associated with $u$ is the tangent vector $\dot\alpha_u(0)$, where $\alpha_u(t)\defi\exp(tu)gK$. By \cite{paradan}*{Equation (4.45)}, the connection $\nabla$ is given by  
  \[
    \nabla_{u_X}=\ad u_X-\Lambda(u)
  \] 
  where $\Lambda:\ger g\to\Gamma(TX)$ is a $G$-equivariant map. By \cite{paradan}*{Propositions 5.2 and 5.9}, $\Lambda$ is determined by $\Lambda(u)(o)=[1,\lambda(u)]$, where 
  \[
    \lambda(x+y)\defi\ad x\mathfa x\in\ger k,y\in\ger p,
  \]
  in view of \cite{paradan}*{Proposition 5.9}. By $G$-equivariance, we have 
  \[
    \Lambda(u)(gK)=T_oL_g\circ\lambda\Parens1{\Ad(g^{-1})(u)}\circ T_oL_{g^{-1}}.
  \]

  Thus, for $[g,u]\in T_{gK}X$, $v\in\Gamma(\sh T_X)$, and any smooth function $f$ on an open neighbourhood of $gK$, we have 
  \[
    (\nabla_{[g,u]}v)_{gK}f=([u_X,v]f)(gK)-[g,\lambda(\Ad(g^{-1})(u))(T_oL_g)^{-1}v_{gK}]f.
  \]
  Here, $[u_X,v]$ is the bracket of vector fields, whereas in the second term on the right-hand side, the brackets denote an equivalence class in $G\times^K\ger p$.

  By the invariance of the trace, we may now compute 
  \[
    \tr_{T_{gK}X}(\nabla v)=\tr_{T_oX}\Parens1{[1,u]\mapsto T_{gK}L_{g^{-1}}(\nabla_{[g,u]}v)_{gK}}.
  \]
  Denoting the action of $G$ on vector fields by $g\cdot v$, we compute, using the $G$-invariance of the Levi-Civita connection, that 
  \[
    \nabla_{[g,u]}v_{gK}=T_oL_g\Parens1{\nabla_{[1,u]}T_gL_{g^{-1}}v_{gK}}=T_oL_g\Parens1{\nabla_{[1,u]}(g^{-1}\cdot v)_o}.
  \]
  In particular, since $\lambda(u)=0$ for all $u\in\ger p$, we find
  \[
    (\ddiv v)(gK)=\tr_{\ger p}\Parens1{u\mapsto[u_X,g^{-1}\cdot v]_o}
  \]
  for any vector field $v$. We may replace $\tr_\ger p$ by $\frac12\tr_{\ger p_\cplxs}$.

  Using the identification $T^\cplxs j$ of $X\times W$ with $T^\cplxs X=G\times^K\ger p_\cplxs$ given by the $G$-action on $W$, we see $(v_n)_{gK}=[g,g^{-1}.n]$. Hence, one finds that 
  \begin{align*}
    (g^{-1}\cdot v_n)_{hK}&=T_{ghK}L_{g^{-1}}(v_n)_{ghK}=T_{ghK}L_{g^{-1}}[gh,(gh)^{-1}.n]\\
    &=[h,h^{-1}.g^{-1}.n]=(v_{g^{-1}.n})_{hK}.
  \end{align*}
  Computing the commutator 
  \[
    [u_X,v_n]=\frac d{ds}\exp(-su)\cdot v_n\Big|_{s=0}=\frac d{ds}v_{\exp(-su).n}\Big|_{s=0}=-v_{u.n},
  \]
  we see that 
  \[
    (\ddiv v_n)(gK)=\tfrac12\tr_{\ger p_\cplxs}\Parens1{u\mapsto[u_X,v_{g^{-1}.n}]_o}=-\tfrac12\tr_{\ger p_\cplxs}R_{g^{-1}.n},
  \]
  thus completing the proof of the lemma.
\end{proof}

\section{Superdistributions and Laplace transforms}\label{app:super-lap}

In this appendix, we develop some basic Euclidean Fourier analysis for superdistributions. The facts about Fourier inversion on the Schwartz space are well-known, but we are not aware of a convenient reference. The account we give of the Laplace transform is to our knowledge new. 

\subsection{Superdistributions}\label{app:schwartz}

In this subsection, we give a self-contained development of the basic functional analytic properties of the spaces of superfunctions and superdistributions we encounter in this article. 

All the results can also be derived quickly from the classical case by invoking Batchelor's theorem. However, we deliberately avoid this point of view on grounds that it is generally useful to have definitions of the relevant topologies at hand, which do not appeal to coordinates in their definition. 

We use some well-established functional analysis terminology at liberty. Basic texts include Refs.~\cites{meise-vogt,schaefer,treves}.

\begin{Def}[temp-schwartz][temperedness and Schwartz class]
  Let $\csvr V$ be a real vector space and $\Norm0\cdot$ a norm function on $\csvr V$. A complex-valued function $f$ on $\csvr V$ is of \Define{moderate growth} if there exists an $N\in\nats$ \scth 
  \[
    \sup\nolimits_{x\in \csvr V}(1+\Norm0x)^{-N}\Abs0{f(x)}<\infty.
  \]  
  Now, let $(\csv V)$ be a \emph{cs} vector space. Consider the associated \emph{cs} manifold $L(\csv V)$, with sheaf of superfunctions $\sh O_{\csv V}$. For $f\in\Gamma(\sh O_{\csv V})$ and $D\in S(V)$, considered as a differential operator, we write 
  \[
    f(D;x)\defi(Df)(x)
  \]
  \fa $x\in \csvr V$. Then $f$ is called \Define{tempered} if for any $D\in S(V)$, the function $f(D;\cdot)$ is of moderate growth; it is of \Define{Schwartz class} if for any $D\in S(V)$, and any tempered superfunction $h$, the function $(hf)(D;\cdot)$ is bounded. In the latter case, we set
  \[
    p_{h,D}(f)\defi\sup\nolimits_{x\in \csvr V}\Abs0{(hf)(D;x)}.
  \]
  
  Then $f$ is tempered if any only if for every $D\in S(V)$, we have 
  \[
    \sup\nolimits_{x\in \csvr V}(1+\Norm0x)^{-N}\Abs0{f(D;x)}<\infty
  \]
  \fs $N>0$; $f$ is of Schwartz class if and only if for $D\in S(V)$, we have   
  \[
    p_{N,D}(f)\defi\sup\nolimits_{x\in \csvr V}(1+\Norm0x)^N\Abs0{f(D;x)}<\infty
  \]
  \fa $N>0$. The totality of all tempered superfunctions (resp.~superfunctions of Schwartz class) is denoted by $\mathscr T(\csv V)$ (resp.~$\Sw0{\csv V}$).  We endow $\Sw0{\csv V}$ with the locally convex topology defined by the seminorms $p_{h,D}$ (or, equivalently, $p_{N,D}$), and let $\TDi0{\csv V}$ denote the topological dual space of $\Sw0{\csv V}$, with the strong topology. The elements of $\TDi0{\csv V}$ are called \Define{tempered superdistributions}. 
\end{Def}

\begin{Lem}
  The space $\Sw0{\csv V}$ is a Fr\'echet space and in particular, barrelled.
\end{Lem}

\begin{proof}
  In view of the above discussion, $\Sw0{\csv V}\cong\Sw0{\csvr V}^m$ where $m=2^{\dim V_\odd}$. The assertion follows from \cite{gelfand-shilov-vol2}*{Chapter II, Section 2.2, Theorem 2}.
\end{proof}

If $X$ is a \emph{cs} manifold, then we endow $\sh O_X(U)$, for any open set $U\subseteq X_0$, with the locally convex topology induced by the seminorms 
\[
  p_{K,D}(f)\defi\sup_{x\in K}\Abs1{(Df)(x)},
\]
where $D$ runs throught the set $\sh D_X(U)$ of superdifferential operators (of finite order) on $X_U$, and $K\subseteq U$ is compact. In what follows, we require $X_0$ to be \emph{metrisable} (or equivalently, paracompact \cite{spivak}*{Appendix}). 

\begin{Prop}[cover-projlim]
  Let $(U_\alpha)$ be an open cover of $U$. Then $\sh O_X(U)$ is the locally convex projective limit of the $\sh O_X(U_\alpha)$, with respect to the restriction morphisms. In particular, $\sh O_X(U)$ is complete, and if $U$ is $\sigma$-compact, then $\sh O_X(U)$ is Fr\'echet. 
\end{Prop}

\begin{proof}
  Since $\sh D_X$ is an $\sh O_X$-module and $\sh O_X$ is $c$-soft, so is $\sh D_X$. This readily implies that the restriction maps are continuous. Hence, we have that the linear map $\sh O_X(U)\to\varprojlim_\alpha\sh O_X(U_\alpha)$ is continuous, and is bijective by the sheaf property. Conversely, to see that it is open, one may pass to a locally finite refinement, and then argue similarly using partitions of unity. The remaining statements then carry over from the case of coordinate neighbourhoods, which is easily dealt with. 
\end{proof}

We recall that if $A$ is an algebra (not necessarily unital or associative) endowed with a locally convex topology, then this topology is called \Define{locally $m$-convex} if it is generated by a system of submultiplicative seminorms. 

\begin{Cor}[smooth-mconvex]
  The topology on $\sh O_X(U)$ is locally $m$-convex. 
\end{Cor}

\begin{proof}
  In view of \thmref{Prop}{cover-projlim}, it is sufficient to prove this in a coordinate neighbourhood. Then, as in the even case \cite{michor-cinftyalgs}*{2.2}, a locally $m$-convex topology is generated by the seminorms
  \[
    p_{K,k}(f)=2^k\cdot \max_{\alpha\in\nats^p\times\{0,1\}^q\Abs0\alpha\sle k}\sup_{x\in K}\Abs0{f(\partial^\alpha,x)},
  \]
  for $k\in\nats$ and $K\subseteq U$ compact, where we agree to write
  \[
    \partial^\alpha\defi\frac{\partial^{\alpha_{p+q}}}{\partial (x^{p+q})^{\alpha_{p+q}}}\dotsm\frac{\partial^{\alpha_1}}{\partial (x^1)^{\alpha_1}},
  \]
  and $(x^a)$ is some local coordinate system on $U$.
\end{proof}

By similar arguments as the proof of \thmref{Prop}{cover-projlim}, one proves the following. 

\begin{Prop}[mor-cont]
  Let $\phi:X\to Y$ be a morphism of \emph{cs} manifolds. Then the even linear pullback map $\phi^\sharp:\Gamma(\sh O_Y)\to\Gamma(\sh O_X)$ is continuous. 
\end{Prop}

Let $X$ be a \emph{cs} manifold where $X_0$ is $\sigma$-compact.\footnote{By assumption, $X_0$ is metrisable, and it is locally path-connected as a manifold. Hence, $X_0$ is second countable if and only if it has countably many connected components \cite{spivak}*{Appendix}.} For any compact $K\subseteq X_0$, let $\Gamma_K(\sh O_X)$ denote the set of all global sections of $\sh O_X$ with support in $K$. Endowed with the relative topology from $\Gamma(\sh O_X)=\sh O_X(X_0)$, it is a Fr\'echet space. 

Let $\Gamma_c(\sh O_X)=\bigcup_K\Gamma_K(\sh O_X)$, where the union extends over all compact subsets $K\subseteq X_0$, be the set of all compactly supported sections of $\sh O_X$, equipped with the locally convex inductive limit topology.

\begin{Prop}[cpt-nuclear]
  The locally convex space $\Gamma_c(\sh O_X)$ has the following properties: 
  \begin{enumerate}
    \item It is LF, and in particular, complete, barrelled, and bornological. 
    \item It is nuclear, and in particular, reflexive and Montel.
  \end{enumerate}
  The latter statement also holds for $\Gamma(\sh O_X)$.
\end{Prop}

\begin{proof}
  (i). If $K'\supseteq K$, then $\Gamma_K(\sh O_X)\to\Gamma_{K'}(\sh O_X)$ is by definition a topological embedding. Since $X_0$ is $\sigma$-compact, the limit topology is computed by taking any countable exhaustive filtration of $X_0$ by compact subsets.  
  
  (ii). It is sufficient to prove the nuclearity of $\Gamma_K(\sh O_X)$, since this property is preserved under countable locally convex inductive limits \cite{treves}*{Proposition 50.1}. The same holds true for locally convex projective limits (\loccit), so the question is reduced to the case of \emph{cs} domain, in view of \thmref{Prop}{cover-projlim}. In this case, $\Gamma_K(\sh O_X)\cong\sh C^\infty_K(U)^N$ where $U\subseteq\reals^p$ is open, $K\subseteq U$ is compact, and $N=2^q$ is some non-negative integer. The claim then follows from \loccit{} and the Corollary to Theorem 51.4 (\opcit). For $\Gamma(\sh O_X)$, we argue analogously.
\end{proof}

\begin{Cor}[sw-nuclear]
  The locally convex spaces $\Sw0{\csv V}$ and $\TDi0{\csv V}$ are nuclear, barrelled, reflexive, and Montel. 
\end{Cor}

The \emph{proof} makes use of the following lemma.

\begin{Lem}
  Both of the natural even linear maps $\Gamma_c(\sh O_{\csv V})\to\Sw0{\csv V}$ and $\TDi0{\csv V}\to\Gamma_c(\sh O_{\csv V})'$ are continuous, injective, and have dense image.
\end{Lem}

\begin{proof}
  For any compact $K\subseteq \csvr V$, we have an injection $\Gamma_K(\sh O_{\csv V})\to\Sw0{\csv V}$, which is continuous by definition of the topologies. Thus, so is $\Gamma_c(\sh O_{\csv V})\to\Sw0{\csv V}$. This map has dense image, by \cite{gelfand-shilov-vol2}*{Chapter II, Section 2.5}. Hence, its transpose defines an injection $\TDi0{\csv V}\to\Gamma_c(\sh O_{\csv V})'$ with dense image. 
\end{proof}

\begin{proof}[\protect{Proof of \thmref{Cor}{sw-nuclear}}]
  By \thmref{Prop}{cpt-nuclear}, $\Gamma_c(\sh O_{\csv V})$ is a nuclear LF space. Therefore, by \cite{treves}*{Proposition 50.6}, $\Gamma_c(\sh O_{\csv V})'$ is the locally convex projective limit of nuclear spaces, and thus itself nuclear, in view of Theorem 50.1 (\opcit). As a subspace of a nuclear space, $\TDi0{\csv V}$ is nuclear (\loccit). Hence, so is $\Sw0{\csv V}$, by Proposition 50.6 (\opcit). 
  
  Any barrelled nuclear space is Montel, any nuclear space is reflexive, and the strong dual of a Montel space is Montel (thus, barrelled). Hence the claim.
\end{proof}

Using nuclearity, we derive along the lines of \cite{treves}*{proof of Theorem 51.6} the following corollary.

\begin{Cor}[prod-space]
  Let $X$ and $Y$ be \emph{cs} manifolds. There is a natural isomorphism of locally convex super-vector spaces $\Gamma(\sh O_X)\widehat\otimes_\pi\Gamma(\sh O_Y)\to\Gamma(\sh O_{X\times Y})$ where $\widehat\otimes_\pi$ denotes the completed projective tensor product.
\end{Cor}

Let $X$ be a \emph{cs} manifold. The assignment $U\mapsto\Gamma_c(\sh O_X|_U)$ is a cosheaf \cite{bredon}, and its extension maps $\Gamma_c(\sh O_X|U)\to\Gamma_c(\sh O_X|_V)$ for open subsets $U\subseteq V\subseteq X_0$ are continuous, as follows from the definition of the topologies. 

Thus, we have a presheaf $\Db_X$ on $X_0$, defined by 
\[
  \Db_X(U)\defi\Gamma_c(\sh O_X|_U)',
\]
the topological dual space of $\Gamma_c(\sh O_X|_U)$. Because $\sh O_X$ is $c$-soft, it follows easily that $\Db_X$ is a sheaf. Sections of this sheaf are called \Define{superdistributions} on $X$. 

In view of \thmref{Prop}{cpt-nuclear}, when equipped with the strong topology, $\Db_X(U)$ is nuclear and Montel, and in particular, reflexive and barrelled. 

\subsection{Vector-valued superfunctions}\label{app:vectval}

In this subsection, we generalise the notion of a function with values in a locally convex space to the super case. Our motivation is Laurent Schwartz's approach to the study of the Laplace transform \cite{schwartz-laplace}, which we will need to super-extend in order to prove the main result of this paper. However, the notion of vector-valued superfunctions is also useful in other contexts.

Rather than giving the most general definition, which would appeal to some category of infinite-dimensional supermanifolds, we define vector-valued superfunctions \via completed tensor products. This becomes tractable by a suitable extension of the formalism of $S$-valued points.

In what follows, let $E$ denote a locally convex super-vector spaces and $E'$ its strong continuous linear dual space. For any \emph{cs} manifold $S$, we define
\[
  \sh O(S,E)\defi\Gamma(\sh O_S)\widehat\otimes_\pi E,
\]
where $\widehat\otimes_\pi$ denotes the completed projective tensor product, endowed with the standard grading. The elements of $\sh O(S,E)$ are called \Define{$E$-valued superfunctions on $S$}. Observe that since $\Gamma(\sh O_S)$ is nuclear by \thmref{Prop}{cpt-nuclear}, we might have taken any other locally convex tensor product topology in the definition \cite{treves}. 

\begin{Prop}
  Let $E$ be a locally convex super-vector space. The assignment $S\mapsto\sh O(S,E)$ is a functor from \emph{cs} manifolds to the category of locally convex super-vector spaces with even continuous linear maps. 

  For any \emph{cs} manifolds $S$ and $T$, there is a natural isomorphism
  \[
    \sh O(S\times T,E)=\sh O(S,\sh O(T,E))
  \]
  of locally convex super-vector spaces. 
\end{Prop}

\begin{proof}
  The functoriality of $\sh O(-,E)$ follows from the definitions: Given a morphism $\phi:T\to S$ of \emph{cs} manifolds, we may form 
  \[
    \sh O(\phi,E)\defi\phi^\sharp\,\widehat\otimes_\pi\id_E:\sh O(S,E)\to\sh O(T,E),
  \]
  by \thmref{Prop}{mor-cont}. The second assertion is a consequence of \thmref{Cor}{prod-space}.
\end{proof}

\begin{Def}[spt-val-vectval][values of vector-valued superfunctions]
  Let $X$ be a \emph{cs} manifold and $f\in\sh O(X,E)$. For any $x\in_SE$, we define 
  \[
    f(x)\defi\sh O(x,E)(f)=\Parens1{x^\sharp\,\widehat\otimes_\pi\id_E}(f)\in\sh O(S,E),
  \]
  and call this the \Define{value of $f$} at the $S$-valued point $x$.
\end{Def}
 
The following is immediate from the definitions.

\begin{Prop}[vectval-nat]
  Let $X_1,\dotsc,X_n$ be \emph{cs} manifolds, $E_1,\dotsc,E_n,F$ be locally convex spaces and $b:\prod_jE_j\to F$ an even continuous $n$-linear map. The assignment 
  \[
    (f_1,e_1,\dotsc,f_n,e_n):\prod\nolimits_j\Parens1{\Gamma(\sh O_{X_j})\times E_j}\mapsto f_1\otimes\dotsm\otimes f_b\otimes b(e_1,\dotsc,e_n)
  \]
  extends uniquely to a continuous linear map
  \[
    b:\sh O\Parens1{\textstyle\prod_jX_j,\widehat\bigotimes_{\pi,j}E_j}\to\sh O\Parens1{\textstyle\prod_jX_j,F}
  \]
  which satisfies
  \[
    b(f)(x_1,\dotsc,x_n)=b(f(x_1,\dotsc,x_n))\in\sh O(S,F)
  \]
  for any $f\in\sh O\Parens1{\prod_jX_j,\widehat\bigotimes_{\pi,j} E_j}$ and $(x_1,\dotsc,x_n)\in_S\prod_jX_j$.
\end{Prop}

\begin{Rem}[vectval-sep]
  The conclusion of \thmref{Prop}{vectval-nat} continues to hold if there is some integer $k\sle n$ \scth 
  \begin{enumerate}[wide]
    \item the spaces $E_1,\dotsc,E_k$ are nuclear,
    \item $b$ is separately continuous, and
    \item for any $(e_1,\dotsc,e_k)\in\prod_{j=1}^kE_j$, the following map is continuous:
    \[
      b(e_1,\dotsc,e_k,\cdot):E_{k+1}\times\dotsm\times E_n\to F.
    \]
  \end{enumerate} 
Thus, if $E$ is nuclear and $\Dual0\cdot\cdot_E$ denotes the canonical pairing $E'\times E\to\cplxs$, then 
\[
  \Dual0\mu f(x,y)=\Dual0{\mu(x)}{f(y)}\in\sh O(S,\cplxs)=\Gamma(\sh O_S)
\]
for any $\mu\in\sh O(X,E')$, $f\in\sh O(Y,E)$, and $(x,y)\in_SX\times Y$.
\end{Rem}

\subsection{Fourier transform on the Schwartz space $\mathscr S$}\label{app:ft-schwartz}

In this subsection, we extend the classical theory of the Fourier transform on the Schwartz space to the super setting. Everything is more or less straightforward. However, we do not know a convenient reference, and consider it worthwhile to supply one. 

In what follows, recall the facts and definitions from Appendix \ref{app:super-int}.

\begin{Def}[std-ber]
  Let $(\csv V)$ be a \emph{cs} vector space of $\dim V=p|q$, endowed with a homogeneous basis $(v_a,\nu_b)$, where we assume $v_a\in \csvr V$. Let $\smash{(v^a,\nu^b)}$ be the dual basis. The Lebesgue density $\Abs0{dv_0}$ is the unique translation invariant density on $\csvr V$ \scth the cube with side one spanned by $(v_a)$ has volume one. 
  
  Moreover, there is a unique Berezinian density $\Abs0{Dv}$ on the \emph{cs} manifold $L(\csv V)$ associated with $(\csv V)$, \scth 
  \[
    \fibint[\csvr V]{L(\csv V)}\Abs0{Dv}\,f=\Abs0{dv_0}\,\frac\partial{\partial\nu^q}\dotsm\frac\partial{\partial\nu^1}f\mathfa f\in\Gamma(\sh O_{\csv V}).
  \]
  We let $(\csv{V^*})$ be the dual \emph{cs} vector space, with densities $\Abs0{dv_0^*}$ and $\Abs0{Dv^*}$ associated with the dual basis $\smash{(v^a,\nu^b)}$.
  
  The \Define{standard retraction} $r=\smash{r_{\csv V}}$ of $L(\smash{\csv V})$ is defined by $r^\sharp(v^a)=v^a$. The definition is in fact independent of the choice of basis. 

  For $f\in\Sw0{\csv V}$, we define the \Define{Fourier transform} $\sh F(f)\in\Sw0{\csvst V}$ by 
  \[
    \sh F(f)\defi\frac1{(2\pi)^{p/2}}\int_{L(\csv V)}\Abs0{Dv}\,e^{-i\Dual0\cdot v}f(v),
  \]
  where $\Dual0\cdot\cdot:V^*\times V\to\cplxs$ denotes the canonical pairing. For $f\in\Sw0{\csvr V}$, we normalise the ordinary Fourier transform $\sh F_0(f)\in\Sw0{\csvrst V}$ by
  \[
    \sh F_0(f)\defi\frac1{(2\pi)^{p/2}}\int_{\csvr V}\Abs0{dv_0}\,e^{-i\Dual0\cdot{v_0}}f(v_0).
  \]
  
  For $f\in\Sw0{\csvst V}$, we define the \Define{Fourier cotransform} $\check{\sh F}(f)\in\Sw0{\csv V}$:
  \[
    \check{\sh F}(f)\defi\frac{(-i)^q(-1)^{\frac12q(q+1)}}{(2\pi)^{p/2}}\int_{L(\csv{V^*})}\Abs0{Dv^*}\,e^{i\Dual0{v^*}\cdot}f(v^*).
  \]
  The ordinary Fourier cotransform $\check{\sh F}_0(f)\in\Sw0{\csvr V}$ of $f\in\Sw0{\csvrst V}$ is  
  \[
    \check{\sh F}_0(f)=\frac1{(2\pi)^{p/2}}\int_{\csvr{V^*}}\Abs0{dv_0^*}\,e^{i\Dual0{v_0^*}\cdot}f(v_0^*),
  \]
  so that $\check{\sh F_0}=\sh F_0^{-1}$ by the classical Fourier inversion theorem \cite{dieudonne-vol6}. 
  
  From the classical theory, one deduces easily that $\sh F$ and $\check{\sh F}$ are continuous. 
\end{Def}

\begin{Lem}[pairing-basis]
  Let $(U,U_{\ev,\reals})$ be a \emph{cs} vector space and $\Dual0\cdot\cdot$ denote the canonical pairing $U^*\times U\to\cplxs$. If $(u_a)$ is a homogeneous basis of $U$ with dual basis $(u^a)$, so that $u^b(u_a)=\delta_{ab}$, then as a superfunction on $L(U^*,U^*_{\ev,\reals})\times L(U,U_{\ev,\reals})$, we have
  \[
    \Dual0\cdot\cdot=\sum_a u_a\otimes u^a.
  \]
\end{Lem}

\begin{proof}
  Let $S$ be any \emph{cs} manifold, and let $u^*\in_SL(U^*,U_{\ev,\reals}^*)$, $u\in_SL(U,U_{\ev,\reals})$, where 
  \[
    u^*=\sum_af_au^a\nd u=\sum_ag^au_a
  \]
  under the identification $U(S)=(\Gamma(\sh O_S)\otimes U)_{\ev,\reals}$ (and similarly for $U^*$). Then 
  \[
    \Dual0{u^*}u=\sum_{ab}(-1)^{\Abs0{u^a}\Abs0{u_b}}f_ag^bu^a(u_b)=\sum_a(-1)^{\Abs0{u_a}}f_ag^a=\str(f_ag^b).
  \]
  On the other hand,
  \[
    \Parens2{\sum\nolimits_au_a\otimes u^a}(u^*,u)=\sum_{abc}(-1)^{\Abs0{u_a}\Abs0{u^b}+\Abs0{u^a}\Abs0{u_c}}f_bg^cu_a(u^b)u^a(u_c)=\str(f_ag^b),
  \]
  since $u_a(u^b)=(-1)^{\Abs0{u_a}}\delta_{ab}$. 
\end{proof}

\begin{Lem}[exp-deriv]
  For any $f\in\Sw0{\csv V}$, $g\in\Sw0{\csvst V}$, and $a=1,\dotsc,p+q$, we have
  \begin{align*}
    \sh F(\nu^af)&=(-1)^qi\frac\partial{\partial\nu_a}\sh F(f),&\check{\sh F}(\nu_a g)&=(-1)^qi\frac\partial{\partial\nu^a}\check{\sh F}(g),\\
    \sh F\Parens2{\frac\partial{\partial\nu^a}f}&=-(-1)^qi\nu_a\sh F(f),&\check{\sh F}\Parens2{\frac\partial{\partial\nu_a}g}&=-(-1)^qi\nu^a\check{\sh F}(g).
  \end{align*}
\end{Lem}

\begin{proof}
  By \thmref{Lem}{pairing-basis}, we have 
  \[
    \frac\partial{\partial\nu_a}\Dual0\cdot\cdot=\nu^a\nd\frac\partial{\partial\nu^a}\Dual0\cdot\cdot=-\nu_a,
  \]
  so that 
  \[
    \frac\partial{\partial\nu_a}e^{\pm i\Dual0\cdot\cdot}=\sum_{n=0}^\infty\frac{(\pm i)^{n+1}}{n!}\nu^a\Dual0\cdot\cdot^n=\pm i\nu^a e^{\pm i\Dual0\cdot\cdot}
  \]
  and $\frac\partial{\partial\nu^a}e^{\pm i\Dual0\cdot\cdot}=\mp i\nu_ae^{\pm i\Dual0\cdot\cdot}$. Then 
  \begin{align*}
    \sh F(\nu^af)&=\int_{L(\csv V)}\Abs0{Dv}\,\nu^ae^{-i\Dual0\cdot v}f(v)\\
    &=i\int_{L(\csv V)}\Abs0{Dv}\,\frac\partial{\partial\nu_a}e^{-i\Dual0\cdot v}f(v)=(-1)^qi\frac\partial{\partial\nu_a}\sh F(f),
  \end{align*}
  and similarly $\check{\sh F}(\nu_a g)=(-1)^qi\frac\partial{\partial\nu^a}\check{\sh F}(g)$. Using integration by parts, we see that
  \[
    \sh F\Parens2{\frac\partial{\partial\nu^a}f}=-\int_{L(\csv V)}\Abs0{Dv}\,\frac\partial{\partial\nu^a}e^{-i\Dual0\cdot v}f(v)=-(-1)^qi\nu_a\sh F(f),
  \]
  and similarly $\check{\sh F}\Parens1{\tfrac\partial{\partial\nu_a}g}=-(-1)^qi\nu^a\check{\sh F}(g)$. 
\end{proof}

\begin{Prop}[fourier-inv]
  The linear map $\sh F:\smash{\Sw0{\csv V}\to\Sw0{\csvst V}}$ is an isomorphism of locally convex vector spaces of parity $\equiv q\,(2)$, with inverse $\check{\sh F}$. 
\end{Prop}

Note that when considered on the level of Berezinian densities instead of functions, the Fourier transform is an even map. 

\begin{proof}[\protect{Proof of \thmref{Prop}{fourier-inv}}]
  The idea is to reduce to the classical case by taking derivatives, as is done for $\dim V=0|q$ in Ref.~\cite{guillemin-sternberg-susyequivar}*{Chapter 7}.

  Let $f\in\Sw0{\csvr V}$, considered as an element of $\Sw0{\csv V}$ \via the standard retraction of $V$. Letting $v^*$ denote the generic point of $L(\csvst V)$, the proof of \thmref{Lem}{exp-deriv} shows that 
  \begin{align*}
    \sh F(f)(v^*)&=\frac1{(2\pi)^{p/2}}\int_{\csvr V}\Abs0{dv_0}\,j_{L(\csv V)_0}^\sharp\Parens2{\frac\partial{\partial\nu^q}\dotsm\frac\partial{\partial\nu^1}e^{-i\Dual0{v^*}\cdot}}(v_0)f(v_0)\\
    &=i^q(-1)^{\frac12q(q-1)}\nu_1\dotsm\nu_q\,\sh F_0(f),
  \end{align*}
  where $\sh F_0(f)$ is considered as a superfunction on $L(\csvst V)$ \via the standard retraction. It follows that  
  \begin{align*}
    \check{\sh F}\sh F(f)=\frac1{(2\pi)^{p/2}}\!\int_{\csvr V^*}\!\!\Abs0{dv_0^*}\,e^{i\Dual0{v_0^*}{v_0}}j_{L(V^*,\csvr V^*)_0}^\sharp\Parens2{\frac\partial{\partial\nu_q}\dotsm\frac\partial{\partial\nu_1}\nu_1\dotsm\nu_q\,\sh F_0(f)}(v_0^*)=f,
  \end{align*}
  since $e^{i\Dual0\cdot\cdot}$ is even, and by the classical Fourier inversion formula \cite{dieudonne-vol6}. 
  
  If now $f\in\Sw0{\csv V}$ is arbitrary, then by \thmref{Lem}{exp-deriv},
  \[
    \check{\sh F}\sh F(\nu^af)=(-1)^qi\,\check{\sh F}\Parens2{\frac\partial{\partial\nu_a}\sh F(f)}=\nu^a\check{\sh F}\sh F(f).
  \]
  This reduces the proof of the equation $\check{\sh F}\sh F=\id$ to the subspace $\Sw0{\csvr V}$, which was treated above. For the converse composition, one proceeds analogously. 
\end{proof}

\begin{Def}
  For $f,g\in\Sw0{\csv V}$, define the \Define{convolution} $f*g\in\Sw0{\csv V}$ by demanding that 
  \[
    \int_{L(\smash{\csv V})}\Abs0{Dv}\,(f*g)(v)h(v)=\int_{L(V\times V,\smash{\csvr V\times\csvr V})}\,\Abs0{Dv_1}\Abs0{Dv_2}\,f(v_1)g(v_2)h(v_1+v_2)
  \]
  for any tempered superfunction $h$. Then, for any \emph{cs} manifold $S$, and any $S$-valued point $x\in_SL(\csv V)$, 
  \[
    (f*g)(x)=\int_{L(\csv V)}\Abs0{Dv}\,f(v)g(x-v).
  \]
  The following lemma, which is an easy consequence of the Leibniz rule and H\"older's inequality, shows that indeed $f*g\in\Sw0{\csv V}$. 
\end{Def}

\begin{Lem}
  \Fa $I\subseteq\{1,\dotsc,q\}$, there exist $D_{I1},D_{I2}\in S(V^*)$, \scth
  \[
    \Abs3{\int_{L(\csv V)}\Abs0{Dv}\,(fg)(v)}\sle\sum_I\sup_{v_0\in \csvr V}\Abs1{f(D_{I1};v_0)}\int_{\csvr V}\Abs0{dv_0}\,\Abs1{g(D_{I2};v_0)},
  \]
  for any superfunctions $f,g$ \scth all integrals in question converge absolutely. 
\end{Lem}

The behaviour of convolution products under Fourier transform carries over to the super case. The non-trivial signs are again an artefact introduced by considering the Fourier transform on the level of functions rather than of Berezinian densities.

\begin{Lem}[ft-conv]
  For any $f,g\in\Sw0{\csv V}$, we have 
  \[
    \sh F(f*g)=(-1)^{q\Abs0f}(2\pi)^{p/2}\sh F(f)\sh F(g).
  \]
\end{Lem}

\begin{proof}
  By the definition of $f*g$, and writing $\int_V$ for $\smash{\int_{L(\csv V)}}$, we compute
  \begin{align*}
    \sh F(f*g)&=\frac1{(2\pi)^{p/2}}\int_{V\times V}\Abs0{Dv_1}\Abs0{Dv_2}\,e^{-i\Dual0\cdot{v_1+v_2}}f(v_1)g(v_2)\\
    &=\frac{(-1)^{q\Abs0f}}{(2\pi)^{p/2}}\int_V\Abs0{Dv}\,e^{-i\Dual0\cdot v}f(v)\cdot\int_V\Abs0{Dv}\,e^{-i\Dual0\cdot v}g(v)\\
    &=(-1)^{q\Abs0f}(2\pi)^{p/2}\sh F(f)\sh F(g).
  \end{align*}
  Here, the equality 
  \[
    e^{-i\Dual0u{v_1+v_2}}=e^{-i\Dual0u{v_1}}e^{-i\Dual0u{v_2}}
  \]
  \fa $u,v_j\in_SL(\csv V)$ follows as usual using Cauchy summation, since the exponential series in question converge absolutely in $\Gamma(\sh O_S)$, in view of \thmref{Cor}{smooth-mconvex} and the (elementary) fact that complete locally $m$-convex algebras admit a functional calculus for entire functions, \cf Ref.~\cite{mitiagin-rolewicz-zelazko}. This proves the claim. 
\end{proof}

\begin{Prop}
  Let $f,g\in\Sw0{\csv V}$. Then 
  \[
    \int_{V^*}\Abs0{Dv^*}\,\sh F(f)(v^*)\sh F(g)(v^*)=(-1)^{q\Abs0f}i^q(-1)^{\frac12q(q-1)}\int_V\Abs0{Dv}\,f(v)g(-v),
  \]
  where we write $\int_{V^*}$ for $\smash{\int_{L(\smash{\csvst V})}}$ and $\int_V$ for $\smash{\int_{L(\smash{\csv V})}}$.
\end{Prop}

\begin{proof}
  We have 
  \[
    \check{\sh F}(h)(0)=\frac{(-i)^q(-1)^{\frac12q(q-1)}}{(2\pi)^{p/2}}\int_{V^*}\Abs0{Dv^*}\,h(v^*)
  \]
  \fa $h\in\Sw0{\csvst V}$, so by \thmref{Lem}{ft-conv} and \thmref{Prop}{fourier-inv},
  \begin{align*}
    \int_{V^*}\Abs0{Dv^*}\,\sh F(f)(v^*)\sh F(g)(v^*)&=\frac{(-1)^{q\Abs0f}}{(2\pi)^{p/2}}\int_{V^*}\Abs0{Dv^*}\,\sh F(f*g)(v^*)\\
    &=(-1)^{q\Abs0f}i^q(-1)^{\frac12q(q-1)}\check{\sh F}\sh F(f*g)(0)\\
    &=(-1)^{q\Abs0f}i^q(-1)^{\frac12q(q-1)}\int_V\Abs0{Dv}\,f(v)g(-v),
  \end{align*}
  which was our assertion. 
\end{proof}

Finally, we define the Fourier (co)transform on $\TDi0{\csv V}$ by duality.

\begin{Def}
  For any homogeneous $\mu\in\TDi0{\csv V}$, we define the distributional \Define{Fourier transform} $\sh F(\mu)\in\TDi0{\csvst V}$ by 
  \[
    \Dual0{\sh F(\mu)}{\sh F(f)}\defi (-1)^{q\Abs0\mu}i^q(-1)^{\frac12q(q-1)}\Dual0\mu{\check f}\mathfa f\in\Sw0{\csv V}
  \]
  where $\check f(v)=f(-v)$ \fa $v\in_SL(\csv V)$ and any \emph{cs} manifold $S$. Similarly, we define for any homogeneous tempered superdistribution $\nu\in\TDi0{\csvst V}$, the \Define{Fourier cotransform} $\check{\sh F}(\nu)\in\TDi0{\csv V}$ by 
  \[
    \Dual0{\check{\sh F}(\nu)}{\check{\sh F}(g)}\defi (-1)^{q\Abs0\nu}(-i)^q(-1)^{\frac12q(q-1)}\Dual0\nu{\check g}\mathfa g\in\Sw0{\csvst V}.
  \]
  The following is immediate.
\end{Def}

\begin{Cor}
  The Fourier transform $\sh F:\TDi0{\smash{\csv V}}\to\TDi0{\smash{\csvst V}}$ is an isomorphism of locally convex vector spaces, of parity $\equiv q\,(2)$, with inverse $\check{\sh F}$. 
\end{Cor}

Again, the parity problems disappear if instead we consider the Fourier transform on the level of tempered generalised superfunctions (\ie the dual of the space of Berezinian densities of Schwartz class).

\subsection{The Paley--Wiener space $\mathscr Z$} 

We will also need to consider the Fourier transform in situations which do not exhibit the same self-duality as the case of the Schwartz space. A useful receptacle will be the so-called Paley--Wiener space, which arises by Fourier transform of compactly supported smooth functions. 

\begin{Def}[pw-def]
  Fix a positive inner product $\Scp0\cdot\cdot_V$ on $\csvr V$. Write $\Norm0\cdot_V$ for the associated norm function. 
  
  Denote by $L(V)$ the complex supermanifold associated with the complex super-vector space $V$, with sheaf of superfunctions $\sh O_V$. We denote by $\PW0{\csv V}$ the following subspace of $\Gamma(\sh O_V)$, 
  \[
    \PW0{\csv V}\defi\Set1{f\in\Gamma(\sh O_V)}{\exists R>0\,\forall D\in S(V)\,,\,p\in S(V^*)\colon z_{R,D,p}(f)<\infty},
  \]
  where 
  \[
    z_{R,D,p}(f)\defi\sup\nolimits_{v\in V_\ev}\Abs1{e^{-R\Norm0{\Im v}_V}(p\cdot f)(D;v)}
  \]
  \fa $R>0$, $p\in S(V^*)$, $D\in S(V)$ and $f\in\Gamma(\sh O_V)$. Here, $\Im v\defi\frac1{2i}(v-\bar v)$ where $\bar v$ is the complex conjugate of $v$ with respect to the real form $\csvr V$ of $V_\ev$. 
  
  Thus, $\PW0{\csv V}$ consists of holomorphic superfunctions of exponential type; it is called the \Define{Paley--Wiener space} of $(\csv V)$. We also consider for fixed $R>0$ the subspace $\PW[_R]0{\csv V}\subseteq\Gamma(\sh O_V)$ defined by the requirement that all $z_{R,D,p}$, for $p$ and $D$ arbitrary, are finite. With the topology induced by these seminorms, $\PW[_R]0{\csv V}$ is a Fr\'echet space. We endow $\PW0{\csv V}$ with the locally convex inductive limit topology of the spaces $\PW[_R]0{\csv V}$. Obviously, the restriction map $\Gamma(\sh O_V)\to\Gamma(\sh O_{\csv V})$ induces a continuous injection $\PW0{\csv V}\to\Sw0{\csv V}$.
\end{Def}

\begin{Prop}[PW][Paley--Wiener]
  Let $\check{\sh F}:\Sw0{\csv V}\to\Sw0{\smash{\csvst V}}$ be the Fourier cotransform of $(\csvst V)$ and $f\in\PW0{\csv V}$. Then $\smash{\check{\sh F}}(f)\in\Gamma_c(\sh O_{\csvst V})$, and $f\in\PW[_R]0{\csv V}$ if and only if $\supp\smash{\check{\sh F}}(f)\subseteq B_{V^*}(0,R)$, where 
  \[
    B_{V^*}(0,R)=\Set1{v^*\in\csvr V^*}{\Norm0{v^*}_{V^*}\sle R},
  \]
  and $\Norm0\cdot_{V^*}$ denotes the norm dual to $\Norm0\cdot_V$. Moreover, this sets up an isomorphism of locally convex spaces $\Gamma_{B_{V^*}(0,R)}\Parens1{\sh O_{\smash{\csvst V}}}\cong\PW[_R]0{\csv V}$. 
\end{Prop}

\begin{proof}
  In view of \thmref{Lem}{exp-deriv}, the proof is reduced to the case of $f\in\PW0{\csvr V}$. Then by the proof of \thmref{Prop}{fourier-inv}, $\smash{\check{\sh F}}(f)=\nu_1\dotsm\nu_q\smash{\check{\sh F_0}}(f)$ which has support in $B_{V^*}(0,R)$ if and only if this is the case for $\check{\sh F}_0(f)$. Hence, the statement reduces to the classical Paley--Wiener theorem, \vq Refs.~\cites{dieudonne-vol6,gelfand-shilov-vol1,ehrenpreis-analytic1,schwartz-vol2}.
\end{proof}

\begin{Cor}
  The locally convex spaces $\PW[_R]0{\csv V}$ are nuclear Fr\'echet, and $\PW0{\csv V}$ is nuclear LF. In particular, both spaces are complete, reflexive, barrelled, Montel, and bornological. 
\end{Cor}

\begin{proof}
  The topology on $\PW[_R]0{\smash{\csv V}}$ (resp.~$\PW0{\smash{\csv V}}$) is the one induced \via the Fourier transform from $\Gamma_{B_{V^*}(0,R)}(\sh O_{V^*})$ (resp.~$\Gamma_c(\sh O_{V^*})$), and the latter is nuclear Fr\'echet (resp.~nuclear LF) and Montel by \thmref{Prop}{cpt-nuclear}. 
\end{proof}

\begin{Def}
  For any homogeneous $\mu\in\PWFn0{\csv V}$, we define the distributional \Define{Fourier transform} $\sh F(\mu)\in\Gamma(\smash{\Db_{\csvst V}})$ by 
  \[
    \Dual0{\sh F(\mu)}{\sh F(f)}\defi (-1)^{q\Abs0\mu}i^q(-1)^{\frac12q(q-1)}\Dual0\mu{\check f}\mathfa f\in\PW0{\csv V}.
  \]
  Here, we write $\smash{\Db_{\csvst V}}\defi\smash{\Db_{L(\csvst V)}}$. Similarly, we define, for $\nu\in\Gamma(\Db_{\csvst V})$, the \Define{Fourier cotransform} $\smash{\check{\sh F}}(\nu)\in\smash{\PWFn0{\csv V}}$ by 
  \[
    \Dual0{\check{\sh F}(\mu)}{\check{\sh F}(g)}\defi (-1)^{q\Abs0\nu}(-i)^q(-1)^{\frac12q(q-1)}\Dual0\nu{\check g}\mathfa g\in\smash{\Gamma_c(\sh O_{\csvst V})}.
  \]
\end{Def}

\begin{Cor}
  The Fourier transform $\smash{\sh F:\PWFn0{\csv V}\to\smash{\Gamma(\Db_{\csv{V^*}}})}$ is an isomorphism of locally convex vector spaces, of parity $\equiv q\,(2)$, with inverse $\check{\sh F}$. 
\end{Cor}

\subsection{Laplace transforms}\label{app:laplace} 

In this subsection, we give an account of the basics of the Laplace transform. The two main results are that the Laplace transform is injective, and that under mild conditions it can be computed as an integral. Another point, which we discuss at some length, is the extension of generalised superfunctions as functionals to certain larger spaces of test superfunctions, depending on the domains of definition of their Laplace transforms. 

Essentially, we follow the classical exposition by Schwartz \cite{schwartz-laplace} (see also the Exercises 4 and 6 in \cite{dieudonne-vol6}*{Chapter XXII.18}), although we also need to consider the Laplace transform for functionals on $\mathscr Z$ (as in Ref.~\cite{ishihara-laplace}). Moreover, a rigorous account of this theory for superspaces needs to use $S$-valued points; in this, we follow the exposition given in Appendix \ref{app:vectval}. 

\begin{Def}
  A locally convex super-vector space $E$ is called a \Define{test space} for $(\csv V)$ if it is one of $\Gamma_c(\sh O_{\csv V})$, $\PW0{\csv V}$, or $\Sw0{\csv V}$. In this case, $\check E\defi\sh F(E)$ (where $\sh F$ is the Fourier transform on $(\csv V)$) is called the \Define{the dual test space} (for $(\csvst V)$); one also has $\check E=\check{\sh F}(E)$, where $\check{\sh F}$ is the Fourier cotransform on $(\csvst V)$. Notice that all test spaces are contained as dense subspaces in $\Sw0{\csv V}$

  The strong dual $E'$ of a test space is called a space of \Define{generalised functions} for $(\csv V)$. Notice that $\TDi0{\csv V}$ is contained as a dense subspace in any space of generalised functions $E'$. The set $\sh M(E)$ of all functions $f\in\Gamma(\sh O_{\csv V})$ such that $f\cdot E\subseteq E$ in $\Gamma(\sh O_{\csv V})$ is called the \Define{multiplier space} of $E$. We write $\smash{\check E'}$ for the strong dual of the dual test space $\smash{\check E}$.
\end{Def}

For $E=\Gamma_c(\sh O_{\csv V})$, we have $\sh M(E)=\Gamma(\sh O_{\csv V})$. For $E=\Sw0{\csv V}$, we have $\sh M(E)=\mathscr T(\csv V)$, the space of tempered superfunctions \cite{dieudonne-vol6}. Finally, for $E=\PW0{\csv V}$, we have \cite{gelfand-shilov-vol1}
\[
  \sh M(E)=\Set1{f\in\Gamma(\sh O_V)}{\forall D\in S(V)\,\exists R,N>0:z_{R,D,-N}(f)<\infty},
\]
where 
\[
  z_{R,D,-N}(f)\defi\sup_{v\in V_\ev}(1+\Norm0v_V)^{-N}e^{-R\Norm0{\Im v}_V}\Abs0{f(D;v)}.
\]
Notice that $\sh M(E)$ contains the space $S(V^*)$ of superpolynomials for any $E$.

Let $E\subseteq F$ be test spaces for $(\csv V)$ and $\mu\in E'$. For $z\in_SL(V^*)$, we write 
\begin{equation}\label{eq:expins}
  e^{-\Dual0z\cdot}\mu\in_SF'  
\end{equation}
if the series 
\begin{equation}\label{eq:lapdef-series}
  \sum_{k=0}^\infty\frac{(-1)^k}{k!}\Dual0zv^k\mu
\end{equation}
converges in $\sh O(S,E')$, and its limit, denoted by $e^{-\Dual0z\cdot}\mu$, lies in the subspace $\sh O(S,F')$. Here, $v$ denotes the generic point of $L(\csv V)$, so that $\Dual0zv^k$ can be interpreted as an element of 
\[
  \sh O\Parens1{S,\Gamma(\sh O_{\smash{\csv V}})}=\Gamma(\sh O_S)\Hat\otimes_\pi\Gamma(\sh O_{\csv V})=\Gamma\Parens1{\sh O_{S\times L(\csv V)}}.
\]
Hence, we have $\Dual0zv^k\mu\in\sh O(S)\Hat\otimes_\pi E'=\sh O(S,E')$. 

We may thus define
\[
  \gamma_F(\mu)_S\defi\Set1{x\in_SL(\csvst V)}{e^{-\Dual0x\cdot}\mu\in_SF'}.
\]
We also write $\gamma_{\mathscr D}$, $\gamma_{\mathscr S}$, and $\gamma_{\mathscr Z}$, for the case of $E=\Gamma_c(\sh O_{\csv V})$, $E=\Sw0{\csv V}$, and $E=\PW0{\csv V}$, respectively. 

We remark that for any $z=x+iy\in_SL(V^*)$, $x,y\in_SL(\csvst V)$, we have 
\begin{equation}\label{eq:lap-tube}
  e^{-\Dual0z\cdot}\mu\in_SF'\quad\Leftrightarrow\quad x\in\gamma_F(\mu)_S.
\end{equation}

\begin{Def}
  Let $E\subseteq F$ be test spaces for $(\csv V)$ and $\mu\in E'$. For any $x\in\gamma_F(\mu)_S$, we define the \Define{Laplace transform}
  \[
    \LT(\mu)(x)\defi\sh F(e^{-\Dual0x\cdot}\mu)\in\sh O(S,\check F').
  \]
\end{Def}

\begin{Lem}[lap-nat]
  Let $E\subseteq F$ be test spaces for $(\csv V)$ and $\mu\in E'$. If $x\in\gamma_F(\mu)_S$ and $t:T\to S$, then $x(t)\in\gamma_F(\mu)_T$ and 
  \[
    \LT(\mu)(x(t))=(\LT(\mu)(x))(t).
  \]
\end{Lem}

\begin{proof}
  Applying \thmref{Prop}{mor-cont} to exchange $t^\sharp\,\widehat\otimes_\pi\id$ with the limit of the series, we find $x(t)\in\gamma_F(\mu)_T$. The second assertion follows from \thmref{Prop}{vectval-nat} and the continuity of the Fourier transform. 
\end{proof}

We now relate the thus defined Laplace transform of generalised superfunctions to the Laplace transform of ordinary generalised functions. First, observe that the following is immediate by the Leibniz rule.

\begin{Lem}
  Let $E\subseteq F$ be test spaces and $\mu\in E'$. Then 
  \[
    \gamma_F(\nu^b\mu)_S\supseteq\gamma_F(\mu)_S\nd\gamma_F\Parens1{\tfrac\partial{\partial\nu_b}\mu}_S\supseteq\gamma_F(\mu)_S.
  \]
\end{Lem}

Denote the standard retraction of $L(\csv V)$ by $\smash{r_{\csv V}:L(\csv V)\to\csvr V}$. Let $E$ be a test space for $(\csv V)$ and $E_0$ the corresponding test space for $\csvr V$; that is, $E_0=\PW0{\csvr V}$ if $E=\PW0{\csv V}$, \etc{} Then $\smash{r_{\csv V}^\sharp(E_0)\subseteq E}$, so 
\[
  \Dual0{r_{\csv V,\sharp}(\mu)}f\defi\Dual0\mu{r_{\csv V}^\sharp(f)}\mathfa\mu\in E'\,,\,f\in E_0.
\]
defines a continuous even linear map $\smash{r_{\csv V,\sharp}:E'\to E_0'}$. 

\begin{Prop}[super-laplace-even]
  Let $E\subseteq F$ be test spaces and $\mu\in E'$. Then 
  \[
    \gamma_F(\mu)_S=\bigcap_{k=0}^q\bigcap_{1\sle b_1<\dotsm<b_k\sle q}\gamma_{F_0}\Parens1{r_{\csv V,\sharp}(\nu^{b_1}\dotsm\nu^{b_k}\mu)}_S.
  \]
\end{Prop}

\begin{proof}
  The statement is immediate from the Taylor expansion
  \begin{equation}\label{eq:genfn-taylor}
    \Dual0\mu f=\sum_{k=0}^q\sum_{B=(1\sle b_1<\dotsm<b_k\sle q)}\pm\Dual2{r_{\csv V,\sharp}(\nu^B\mu)}{j_{L(\csv V)}^\sharp\Parens2{\frac{\partial^k}{\partial\nu_B}f}},    
  \end{equation}
  which can be applied to $\Dual0{e^{-\Dual0x\cdot}\cdot\mu}f$. Here, we use \thmref{Prop}{vectval-nat} in conjunction with \thmref{Rem}{vectval-sep}.
\end{proof}

If $X$ is any \emph{cs} manifold and $j:Y\to X$ is an embedding, then if a morphism $S\to X$ factors through $j$, then it does so uniquely. Hence, for any open subspace $U$ of $X$, $U(S)$ may be considered a as subset of $X(S)$ for any $S$, and the totality of these subsets form a topology on $X(S)$. For this topology, any morphism $f:X\to Y$ induces continuous maps $X(S)\to Y(S)$ on $S$-valued points; moreover, if $j:Y\to X$ is an embedding, then $j$ is open if and only if $Y(S)\subseteq X(S)$ is open for every $S$.

In what follows, given an open subspace $\gamma\subseteq L(\smash{\csvst V})$, we denote by $T(\gamma)$ the open subspace $\gamma+iL(\smash{\csvst V})$ of $L(V^*)$.

\begin{Th}[lap-hol]
  Let $E\subseteq F$ be test spaces for $(\smash{\csv V})$ and $\mu\in E'$. 
  \begin{enumerate}[wide]
    \item There exists a unique open subspace $\gamma_F^\circ(\mu)$ of $L(\smash{\csvst V})$ \scth for any \emph{cs} manifold $S$, $\gamma_F^\circ(\mu)(S)$ is the interior of $\gamma_F(\mu)_S$ in $L(\csvst V)(S)$.
    \item Assume that $\gamma_{\mathscr S}^\circ(\mu)\neq\vvoid$ and let $T(\mu)\defi T(\gamma_\mathscr S^\circ(\mu))$. Then there is a unique $f\in\Gamma(\sh O_{T(\mu)})$ \scth 
    \begin{equation}\label{eq:lap-fn}
      \int_{L(\csvst V)}\Abs0{Dy}\,f(x+iy)\vphi(y)=\Dual0{\LT(\mu)(x)}\vphi
    \end{equation}
    for any $\vphi\in\Sw0{\smash{\csvst V}}$ and $x\in_S\gamma_\mathscr S^\circ(\mu)$, in the sense that the integral on the left hand side converges absolutely and the equality holds. Moreover, for any $D\in S(V^*)$, the function $f(D;z)$ is bounded on $T(K)$ by some polynomial in $\Im z$, for every compact subset $K\subseteq\gamma_\mathscr S^\circ(\mu)_0$.
    \item Given an open subspace $\gamma\subseteq L(\smash{\csvst V})$ and $f\in\Gamma(\sh O_{T(\gamma)})$, there is a tempered superdistribution $\mu\in\TDi0{\smash{\csv V}}$ \scth $\gamma\subseteq\gamma_\mathscr S^\circ(\mu)$ and Equation \ref{eq:lap-fn} holds, if and only if for any $D\in S(V^*)$, the function $f(D;z)$ is bounded on $T(K)$ by a polynomial in $\Im z$, for every compact $K\subseteq\gamma_0$. In this case, $\mu$ is unique. 
  \end{enumerate}
\end{Th}

\begin{proof}
  (i). Unicity is obvious, and so we check existence. In view of \thmref{Prop}{super-laplace-even}, we may assume that $V=V_\ev$ and $\mu\in E'=E_0'$ is an ordinary generalised function. In this case, we let $\gamma_F^\circ(\mu)$ be the set of all $x\in\smash{\csvr V}$ \scth $e^{-\Dual0y\cdot}\cdot\mu\in F'$ for all $y$ a some neighbourhood of $x$. In general, $\gamma_F^\circ(\mu)$ will be defined by 
  \[
    \gamma_F^\circ(\mu)\defi\bigcap_{k=0}^q\bigcap_{1\sle b_1<\dotsm<b_k\sle q}\gamma_{F_0}^\circ\Parens1{r_{\csv V,\sharp}(\nu^{b_1}\dotsm\nu^{b_k}\mu)}.
  \]

  If $\gamma_F(\mu)_S^\circ\neq\vvoid$ where $S_0\neq\vvoid$, then there exists an open subset $U\subseteq\smash{\csvr V}$ \scth for any $y\in_SU$, we have $e^{-\Dual0y\cdot}\cdot\mu\in\sh O(S,F')$. In particular, $e^{-\Dual0u\cdot}\cdot\mu\in F'$ for any $u\in U$ that appears as the value of some $y\in_SU$. But since the image of $S_0$ is non-empty, any $u\in U$ appears in this way (by considering constant morphisms). 

  Then by Equations \eqref{eq:lap-tube} and \eqref{eq:lapdef-series}, the map $z\mapsto e^{-\Dual0z\cdot}\cdot\mu:U+i\smash{\csvr V}\to F'$ is (strongly) holomorphic, and in particular, if $u$ denotes the generic point of $U$, we have $e^{-\Dual0z\cdot}\cdot\mu\in\sh O(U,F')=\Ct[^\infty]0U\widehat\otimes_\pi F'$. In other words, $U\subseteq\gamma_F^\circ(\mu)$. 

  Thus, any $S$-valued point of $L(\csv V)$ in the interior of $\gamma(\mu)_S$ is the specialisation of the identity of some open subset of $\gamma_F^\circ(\mu)$, and this proves the claim. 

  (ii)--(iii). Using \thmref{Prop}{vectval-nat} in conjunction with \thmref{Rem}{vectval-sep}, the statements reduce by the token of Equation \eqref{eq:genfn-taylor} to the classical case \cite{schwartz-laplace}*{Proposition 6}.
\end{proof}

The following is immediate by combining items (ii) and (iii) of \thmref{Th}{lap-hol}.

\begin{Cor}[sw-ext]
  Let $E$ be a test space for $(\smash{\csv V})$ and $\mu\in E'$. If $\gamma_\mathscr S^\circ(\mu)\neq\vvoid$, then $\mu\in\TDi0{\smash{\csv V}}$, in the sense that $\mu$ extends continuously to $\Sw0{\smash{\csv V}}$.
\end{Cor}

\begin{Def}
  Let $E$ be a test space for $(\smash{\csv V})$ and $\mu\in E'$, where we assume $\gamma_\mathscr S^\circ(\mu)\neq\vvoid$. The holomorphic superfunction $f$ on $T(\mu)$ defined in \thmref{Th}{lap-hol} will be denoted by $\LT(\mu)$ and called the \Define{Laplace transform} of $\mu$.
\end{Def}

As above, using Taylor expansion in odd directions, the following two statements are immediate from the classical case \cite{schwartz-laplace}*{Proposition 8, Corollaire, Remarque}.

\begin{Prop}
  Let $E$ be a test space for $(\smash{\csv V})$ and $\mu\in E'$, where we assume $\gamma\defi\gamma_\mathscr S^\circ(\mu)\neq\vvoid$. Let $\xi\in\smash{\csvr V^*}$ and $C\in\reals$. Then the support of $\mu\in\TDi0{\smash{\csv V}}$ is contained in the half-space
  \[
    H_{\xi,C}\defi\Set1{v\in\csvr V}{\Dual0\xi v\sge C}
  \]
  if and only if for every $c<C$ and for some, or equivalently, for all $\xi_0\in\gamma_0$, we have $\xi_0+\reals_{\sge0}\xi\subseteq\gamma_0$ and each of the functions
  \[
    e^{tc}\LT(\mu)(D;\xi_0+t\xi+i\eta)
  \]
  for $D\in S(V^*)$, is bounded for all $t\sge0$ by some polynomial in $\eta$ independent of $t$.
\end{Prop}

\begin{Cor}[supp-lap]
  Let $\mu\in\TDi0{\smash{\csv V}}$ and assume $\supp\mu\subseteq\gamma$ where the latter is a closed convex cone with $\gamma\cap(-\gamma)=0$. Then 
  \[
    \gamma_\mathscr S^\circ(\mu)_0\supseteq\check\gamma\defi\Set1{\xi\in\csvr V^*}{\forall v\in\gamma\setminus\{0\}:\Dual0\xi v>0}.
  \]
\end{Cor}

\begin{Rem}
  A nice and more elementary proof of the latter result (for the classical case) is given in Ref.~\cite{hn-riesz}.
\end{Rem}

\end{document}